\documentclass[a4paper]{article}

\usepackage{import} 

\usepackage{iftex} 
\usepackage{xparse} 
\usepackage{xifthen} 

\usepackage{newsec} 
\usepackage{parskip} 
\usepackage{verbatim} 
\usepackage{todonotes} 
\usepackage{comment} 
\usepackage{enumerate} 
\usepackage{multicol} 
\usepackage{authblk} 
\ifpdftex
	\relax
\else
	\usepackage{fontspec}
	\defaultfontfeatures{Mapping = tex-text , Scale = MatchUppercase}
\fi
\usepackage{amssymb}
\usepackage{amsmath}
	\numberwithin{equation}{section}
\usepackage{mathtools}
\usepackage{stackengine} 
\usepackage{proof}
\ifpdftex
	\usepackage{latexsym} 
	\usepackage{lmodern} 
	\usepackage{frenchmath} 
	\usepackage{upgreek} 
	\usepackage{mathabx} 
	\usepackage{mathbbol} 
\else
	\usepackage{unicode-math} 
	\unimathsetup
	{
		math-style = french , 
		partial = upright ,
		nabla = upright
	}
\fi

\ifpdftex
	\usepackage{newunicodechar} 
	\newunicodechar{′}{{^{\prime}}}
	\newunicodechar{−}{-}
	\newunicodechar{×}{\times}
	\newunicodechar{∧}{\wedge}
	\newunicodechar{∨}{\vee}
	\newunicodechar{⊕}{\oplus}
	\newunicodechar{⊖}{\ominus}
	\newunicodechar{⊗}{\otimes}
	\newunicodechar{⊙}{\odot}
	\newunicodechar{⋅}{\cdot}
	\newunicodechar{∘}{\circ}
	\newunicodechar{∙}{\bullet}
	\newunicodechar{°}{{}^\circ}
	\newunicodechar{⋯}{\cdots}
	\newunicodechar{…}{\ldots}
	\newunicodechar{⋆}{\star}
	\newunicodechar{∞}{\infty}
	\newunicodechar{→}{\rightarrow}
	\newunicodechar{←}{\leftarrow}
	\newunicodechar{⇸}{\stackinset{c}{-0.125ex}{c}{0.25ex}{$\shortmid$}{$\rightarrow$}}
	\newunicodechar{⇷}{\stackinset{c}{0.125ex}{c}{0.25ex}{$\shortmid$}{$\leftarrow$}}
	\newunicodechar{⤍}{\dashrightarrow}
	\newunicodechar{⤌}{\dashleftarrow}
	\newunicodechar{⇒}{\Rightarrow}
	\newunicodechar{⇐}{\Leftarrow}
	\newunicodechar{↓}{\downarrow}
	\newunicodechar{↑}{\uparrow}
	\newunicodechar{⇑}{\Uparrow}
	\newunicodechar{⇓}{\Downarrow}
	\newunicodechar{◇}{\Diamond}
	\newunicodechar{□}{\square}
	\newunicodechar{⟨}{\langle}
	\newunicodechar{⟩}{\rangle}
	\newunicodechar{⌜}{\ulcorner}
	\newunicodechar{⌝}{\urcorner}
	\newunicodechar{⌞}{\llcorner}
	\newunicodechar{⌟}{\lrcorner}
	\newunicodechar{∫}{\int}
	\newunicodechar{∏}{\prod}
	\newunicodechar{∐}{\coprod}
	\newunicodechar{∑}{\sum}
	\newunicodechar{≅}{\cong}
	\newunicodechar{≃}{\simeq}
	\newunicodechar{∼}{\sim}
	\newunicodechar{≔}{\coloneqq}
	\newunicodechar{≕}{\eqqcolon}
	\newunicodechar{∈}{\in}
	\newunicodechar{∋}{\ni}
	\newunicodechar{∀}{\forall}
	\newunicodechar{∃}{\exists}
	\newunicodechar{⊤}{\top}
	\newunicodechar{⊥}{\bot}
	\newunicodechar{⊢}{\vdash}
	\newunicodechar{⊣}{\dashv}
	\newunicodechar{∆}{\Delta} 
	\newunicodechar{∇}{\nabla}
	\newunicodechar{∂}{\partial}
	\newunicodechar{α}{\alpha}
	\newunicodechar{β}{\beta}
	\newunicodechar{γ}{\gamma}
	\newunicodechar{δ}{\delta}
	\newunicodechar{ε}{\varepsilon}
	\newunicodechar{η}{\eta}
	\newunicodechar{κ}{\kappa}
	\newunicodechar{λ}{\lambda}
	\newunicodechar{ρ}{\rho}
	\newunicodechar{μ}{\mu}
	\newunicodechar{ν}{\nu}
	\newunicodechar{π}{\pi}
	\newunicodechar{φ}{\varphi}
	\newunicodechar{ψ}{\psi}
	\newunicodechar{θ}{\theta}
	\newunicodechar{σ}{\sigma}
	\newunicodechar{τ}{\tau}
	\newunicodechar{χ}{\chi}
	\newunicodechar{ξ}{\xi}
	\newunicodechar{ζ}{\zeta}
	\newunicodechar{ω}{\omega}
	\newunicodechar{Γ}{\Gamma}
	\newunicodechar{Δ}{\Delta}
	\newunicodechar{Λ}{\Lambda}
	\newunicodechar{Φ}{\Phi}
	\newunicodechar{Ψ}{\Psi}
	\newunicodechar{Ω}{\Omega}
\else
\fi
\usepackage{amsthm}
	\numberwithin{equation}{section}

\newtheoremstyle
	{plainer} 
	{2ex} 
	{2ex} 
	{} 
	{} 
	{\bfseries} 
	{} 
	{\newline} 
	{} 

\theoremstyle{plainer} 
\newtheorem{theorem}{Theorem}

\newtheorem{definition}[theorem]{Definition}

\numberwithin{theorem}{section}
\numberwithin{example}{section} 
\usepackage{xcolor} 
\usepackage{listings} 

\usepackage{tikz}

\usetikzlibrary
{
	calc ,
	matrix ,
	arrows ,
	angles ,
	quotes ,
	intersections ,
	positioning ,
	fit ,
	graphs ,
	shapes.geometric ,
	shapes.misc ,
	decorations ,
	decorations.markings ,
	decorations.pathmorphing ,
	decorations.pathreplacing ,
	decorations.text
}

\ifluatex
	\usetikzlibrary{graphdrawing}
	\usegdlibrary{trees , force}
\else
	\relax
\fi

\usepgflibrary
{
	arrows.meta
}

\usepackage{xcolor}

\tikzset
{
	diagram/.style =
	{
		line cap = butt ,
		line join = bevel ,
		arrows = -> ,
		> = angle 60 ,
		auto = left ,
		font = \small ,
		text height = 1.5ex , 
		text depth = 0.25ex , 
		inner sep = 0.25em ,
	} ,
	smaller/.append style =
	{
		font = \small
	} ,
	code node/.append style =
	{
		every node/.append style =
		{
			execute at begin node = \begin{texttt} ,
			execute at end node = \end{texttt}
		}
	} ,
	code diagram/.style =
	{
		diagram ,
		code node
	} ,
	math node/.append style =
	{
		every node/.append style =
		{
			execute at begin node = \begin{math} ,
			execute at end node = \end{math}
		}
	} ,
	math diagram/.style =
	{
		diagram ,
		math node
	} ,
	string diagram/.style =
	{
		math diagram ,
		arrows = - ,
	} ,
	online/.style =
	{
		shape = rectangle ,
		rounded corners = 2mm ,
		inner sep = 1.5pt ,
		fill = white ,
		opacity = 1.0 ,
		anchor = center	
	} ,
	double arrow/.style =
	{
		double distance between line centers = 2pt ,
		-implies 
	} ,
	equals/.style =
	{
		double distance between line centers = 2pt ,
		-implies , 
		arrows = -
	} ,
	mapto/.append style =
	{
		arrows = |-> ,
	} ,
	paph/.append style =
	{
		decorate ,
		decoration = {snake , segment length = 5mm , amplitude = 0.5mm} ,
	} ,
	includel/.append style =
	{
		arrows = Hooks[right]-> ,
	} ,
	includer/.append style =
	{
		arrows = Hooks[left]-> ,
	} ,
	monomorphism/.append style =
	{
		arrows = >-> ,
	} ,
	epimorphism/.append style =
	{
		arrows = ->> ,
	} ,
	cofibration/.append style =
	{
		arrows = >-> ,
	} ,
	fibration/.append style =
	{
		arrows = ->> ,
	} ,
	equivalence/.append style =
	{
		decoration = {markings , mark = at position 1/2 with {\node[online , transform shape] {∼};}} ,
		postaction = {decorate} ,
	} ,
	mark pos/.store in = \markpos ,
	mark pos = 1/2 ,
	proarrow/.append style =
	{
		postaction = {decorate} ,
		decoration =
		{
			markings ,
			mark = at position \markpos with {\draw [solid , -] (0 , -0.6ex) to (0 , 0.6ex);}
		} ,
	} ,
	proarrows/.append style =
	{
		postaction = {decorate} ,
		decoration =
		{
			markings ,
			mark = at position #1with {\draw [solid , -] (0 , -0.6ex) to (0 , 0.6ex);}
		} ,
	} ,
	heteromorphism/.append style =
	{
		decoration = {markings , mark = at position \markpos with {\draw [white , solid , semithick , -] (-1.0ex , 0) to (1.0ex , 0);}} ,
		postaction = {decorate} ,
	} ,
	arrow/.append style =
	{
		arrows = -> ,
	} ,
	string/.append style =
	{
		arrows = - ,
	} ,
	cofibration string/.append style =
	{
		decoration = {markings , mark = at position 1/4 with {\arrow{>}} , mark = at position 3/4 with {\arrow{>}}} ,
		postaction = {decorate} ,
	} ,
	fibration string/.append style =
	{
		decoration = {markings , mark = at position 4/8 with {\arrow{>}} , mark = at position 5/8 with {\arrow{>}}} ,
		postaction = {decorate} ,
	} ,
	overcross/.append style =
	{
		preaction = {draw = white , - , line width = 4pt , line cap = round}
	} ,
	label/.append style =
	{
		font = \scriptsize
	} ,
	incoming/.append style =
	{
		pos = 0 ,
		anchor = south
	} ,
	outgoing/.append style =
	{
		pos = 1 ,
		anchor = north
	} ,
	outcoming/.append style =
	{
		pos = 0 ,
		anchor = north
	} ,
	ingoing/.append style =
	{
		pos = 1 ,
		anchor = south
	} ,
	fromabove/.append style =
	{
		pos = 0 ,
		anchor = south
	} ,
	tobelow/.append style =
	{
		pos = 1 ,
		anchor = north
	} ,
	frombelow/.append style =
	{
		pos = 0 ,
		anchor = north
	} ,
	toabove/.append style =
	{
		pos = 1 ,
		anchor = south
	} ,
	fromleft/.append style =
	{
		pos = 0 ,
		anchor = east
	} ,
	toright/.append style =
	{
		pos = 1 ,
		anchor = west
	} ,
	fromright/.append style =
	{
		pos = 0 ,
		anchor = west
	} ,
	toleft/.append style =
	{
		pos = 1 ,
		anchor = east
	} ,
	forall/.append style =
	{
		line width = 0.5pt ,
	} ,
	exists/.append style =
	{
		densely dashed
	} ,
	structural/.append style =
	{
		opacity = 2/3 ,
	} ,
	0-cell/.style =
	{
		shape = rectangle ,
		rounded corners = 2mm ,
		inner sep = 2pt
	} ,
	1-cell/.style =
	{
		shape = rectangle ,
		rounded corners = 2mm ,
		inner sep = 1.5pt ,
		fill = white ,
		opacity = 1.0 ,
		anchor = center	
	} ,
	2-cell/.style =
	{
		draw = black!75 ,
		fill = white ,
		opacity = 0.9 ,
		inner sep = 0.75mm ,
		minimum size = 4.0mm
	} ,
	bead/.style =
	{
		2-cell ,
		shape = rectangle ,
		rounded corners = 2.0mm
	} ,
	diamond/.style =
	{
		2-cell ,
		shape = diamond
	} ,
	bra/.style =
	{
		2-cell ,
		shape = isosceles triangle ,
		isosceles triangle apex angle = 70 ,
		shape border rotate = -90 ,
		minimum size = 5mm
	} ,
	ket/.style =
	{
		2-cell ,
		shape = isosceles triangle ,
		isosceles triangle apex angle = 70 ,
		shape border rotate = 90 ,
		minimum size = 6mm
	} ,
	box/.style =
	{
		2-cell ,
		shape = rectangle
	} ,
	white dot/.style =
	{
		2-cell ,
		shape = circle ,
		minimum size = 1.5mm ,
		fill = white
	} ,
	black dot/.style =
	{
		2-cell ,
		shape = circle ,
		minimum size = 1.5mm ,
		fill = gray
	} ,
	dot/.style =
	{
		2-cell ,
		shape = circle ,
		inner sep = 0mm ,
		minimum size = 0.25mm ,
		fill = gray
	} ,
	sheet/.style =
	{
		draw = black!50 ,
		fill = gray!10 ,
		opacity = 1.0 ,
		fill opacity = 0.5
	} ,
	torn/.append style =
	{
		decoration = {random steps , segment length = 1pt , amplitude = 0.5pt}
	} ,
	on sheet/.append style =
	{
		semithick
	} ,
	edge annotation/.append style =
	{
		anchor = west
	} ,
	compositor/.style =
	{
		draw = black!50 ,
		fill = gray!10 ,
		draw opacity = 0.75 ,
		fill opacity = 0.5 ,
		dashed
	} ,
	callout/.append style =
	{
		inner sep = 0mm ,
		minimum height = 2mm ,
		minimum width = 4mm ,
		draw ,
		draw opacity = 0.25
	} ,
	callout-pre/.append style =
	{
		callout ,
		inner sep = 1.2mm ,
		solid
	} ,
	callout-post/.append style =
	{
		callout ,
		inner sep = 0.8mm ,
		dashed
	} ,
	pullback/.pic =
	{
		\draw [semithick , arrows = -] (-0.1 , 0.1) -- (0.1 , 0.1) -- (0.1 , -0.1) ;
	} ,
	pushout/.pic =
	{
		\draw [semithick , arrows = -] (-0.1 , 0.1) -- (-0.1 , -0.1) -- (0.1 , -0.1) ;
	} ,
}

\pgfdeclaredecoration{simple line}{start}
{
  \state{start}[width = +0pt,
                next state=step]{
    \pgfpathmoveto{\pgfpoint{0pt}{0pt}}
  }
  \state{step}[auto end on length    = 3pt,
               auto corner on length = 3pt,               
               width=+1pt]
  {
    \pgfpathlineto{\pgfpoint{1pt}{0pt}}
  }
  \state{final}
  {}
}
\usepackage
[
	backend = biber ,
	style = alphabetic ,
	maxnames = 6 ,
	url = false ,
]{biblatex} 

\bibliography{bibliography} 
\usepackage{hyperref}
\hypersetup
{
	unicode = true ,  
	bookmarks = true , 
	bookmarksopen = true , 
	bookmarksopenlevel = 2 , 
	breaklinks = true , 
	colorlinks = false , 
	pdfborder = {0 0 0} , 
	pdfborderstyle={/S/U/W 0} ,
}

\NewDocumentCommand \define {s o m}
{%
	\hypertarget%
	{\IfValueTF{#2}{#2}{#3}}%
	{\IfBooleanTF{#1}{#3}{\emph{#3}}}%
}

\NewDocumentCommand \refer {s o m}
{%
	\hyperlink%
	{\IfValueTF {#2}{#2}{#3}}%
	{\IfBooleanTF{#1}{#3}{#3}}%
}

\NewDocumentCommand \ensuretext {m} {\textrm{#1}}

\NewDocumentCommand \bbb {m} {\mathbb{#1}}    
\NewDocumentCommand \scp {m} {\ensuretext{\textsc{#1}}}  

\NewDocumentCommand \Nat {} {\bbb{N}}

\RenewDocumentCommand \arg {o}
{
	\IfNoValueTF{#1}
	{{−}}
	{\overset{#1}{−}}
}

\NewDocumentCommand \infarg {o}
{
	\IfNoValueTF{#1}
	{{\_}}
	{\overset{#1}{\_}}
}

\NewDocumentCommand \bind {o m m} 
{
	\IfNoValueTF{#1}
	{} {\mathop{{#1}}}
	{#2} \mathrel{.} {#3}
}

\NewDocumentCommand \unaryminus {} {\scalebox{0.5}[1.0]{\( - \)}}
\NewDocumentCommand \inv {m} {{#1} ^{\unaryminus1}} 

\NewDocumentCommand \cat {m} {\bbb{#1}} 
\NewDocumentCommand \Cat {m} {\scp{#1}} 
\NewDocumentCommand \celldim {m m} {{{#2}_{#1}}} 

\makeatletter
\newcommand*{\length}[1]{%
    \@tempcnta\z@
    \@for\@tempa:=#1\do{\advance\@tempcnta\@ne}%
    \the\@tempcnta%
}
\makeatother

\makeatletter
\newcount\my@repeat@count
\newcommand{\natrec}[2]{%
  \begingroup
  \my@repeat@count=\z@
  \@whilenum\my@repeat@count<#1\do{#2\advance\my@repeat@count\@ne}%
  \endgroup
}
\makeatother

\DeclarePairedDelimiter \hombrackets {(} {)}
\NewDocumentCommand \makehom {m} 
{
	\DeclareDocumentCommand \f {o m m} 
	{
		\IfNoValueTF{##1}
		{{##2} \mathrel{#1} {##3}}
		{{##1} \, \hombrackets{{##2} \mathrel{#1} {##3}}}
	}
}

\NewDocumentCommand \monoid {m m} 
{
	\DeclareDocumentCommand \f {m} 
	{
		\ifthenelse{\isempty{##1}}
		{
			{#2}
		}
		{
			\foreach \item [count=\index] in {##1}
			{
				\ifthenelse{\equal{\index}{1}}
				{} 
				{#1} 
				{\item}
			}
		}
	}
}

\DeclarePairedDelimiter \setbrackets {\{} {\}}
\NewDocumentCommand \set {m}
{
	\ifthenelse{\isempty{#1}}
	{∅}
	{
		\monoid{\mathpunct{,}}{\,}
		\setbrackets{\f{#1}}
	}
}

\DeclarePairedDelimiter \sequencebrackets {[} {]}
\NewDocumentCommand \sequence {m}
{
	\ifthenelse{\isempty{#1}}
	{∅}
	{
		\monoid{\mathpunct{,}}{\,}
		\sequencebrackets{\f{#1}}
	}
}

\RenewDocumentCommand \hom {o m m}
{
	\makehom{→}
	\f[#1]{#2}{#3}
}

\NewDocumentCommand \comp {O{1} m}
{
	\monoid
	{\mathbin{\natrec{#1}{\mathord{⋅}}}}
	{\mathrm{id} \ifthenelse{#1 > 1}{^{#1}}{}}
	\f{#2}
}

\NewDocumentCommand \idty {o m}
{
	\comp[#1]{} ({#2})
}

\NewDocumentCommand \pmoc {O{1} m}
{
	\monoid
	{\mathbin{\natrec{#1}{\mathord{∘}}}}
	{\mathrm{id} \ifthenelse{#1 > 1}{^{#1}}{}}
	\f{#2}
}

\NewDocumentCommand \id {m}
{
	\mathrm{id}
	\ifthenelse
	{\equal{#1}{}}
	{}
	{({#1})}
}

\NewDocumentCommand \adjoint {m}
{
	\monoid{\mathrel{⊣}}{\mathrm{id}}
	\f{#1}
}

\NewDocumentCommand \product {m}
{
	\monoid{×}{1}
	\f{#1}
}

\DeclarePairedDelimiter \tuplebrackets {⟨} {⟩}
\DeclarePairedDelimiter \globaltuplebrackets {(} {)}
\NewDocumentCommand \tuple {s m}
{
	\IfBooleanTF{#1}
	{
		\monoid{\mathbin{,}}{\,}
		\globaltuplebrackets{\f{#2}}
	}
	{
		\ifthenelse{\isempty{#2}}
		{!}
		{
			\monoid{\mathbin{,}}{\,}
			\tuplebrackets{\f{#2}}
		}
	}
}

\DeclarePairedDelimiter \cotuplebrackets {[} {]}
\NewDocumentCommand \cotuple {m}
{
	\ifthenelse{\isempty{#1}}
	{¡}
	{
		\monoid{\mathbin{,}}{\,}
		\cotuplebrackets{\f{#1}}
	}
}

\NewDocumentCommand \boundary {o}
{
	∂
	{
		\IfNoValueTF{#1}
		{}
		{^{#1}}
	}
}

\NewDocumentCommand \dom {}
{
	\boundary[−]
}

\NewDocumentCommand \cod {}
{
	\boundary[+]
}

\NewDocumentCommand \rightclass {m m} 
{
	{#2}^{#1}
}
\NewDocumentCommand \leftclass {m m} 
{
	{}^{#1}{#2}
}

\NewDocumentCommand \since {m} {\ensuremath{[\text{#1}]}}
\NewDocumentEnvironment{relationalreasoning}{}{\begin{array}{cl}}{\end{array}}
\NewDocumentCommand \term {o o m}
{
	\IfNoValueTF{#2}
	{{#1} & {#3}\\}
	{{#1} & \since{#2} \\ & {#3} \\}
}

\NewDocumentEnvironment{maptable}{}{\begin{array}{rcl}}{\end{array}}

\NewDocumentCommand \sequent {o o m m}
{
	{#3} \mathrel{⊢\IfNoValueTF{#1}{}{^{#1}}\IfNoValueTF{#2}{}{_{#2}}}{#4}
}

\NewDocumentCommand \cancel {o m} 
{
	\IfNoValueTF{#1}
	{[{#2}]}
	{{[{#2}] ^{#1}}}
}

\NewDocumentCommand \interchange {m}
{
	χ_{\tuple*{#1}}
}

\RenewDocumentCommand \lim {} {\underset{⟵}{\mathrm {lim}}}

\NewDocumentCommand \push {o m m} 
{
	\IfNoValueTF{#1}
	{
		{#2}_!
	}
	{
		Σ_{#1} ({#2})
	}
	\ifthenelse{\isempty{#3}}{}{({#3})}
}
\NewDocumentCommand \forward {o m} 
{
	\overrightarrow{#2}
}
\NewDocumentCommand \pull {o m m} 
{
	\IfNoValueTF{#1}
	{
		{#2}^*
	}
	{
		\inv{#1} ({#2})
	}
	\ifthenelse{\isempty{#3}}{}{({#3})}
}
\NewDocumentCommand \backward {o m} 
{
	\overleftarrow{#2}
}

\NewDocumentCommand \cartesianlift {o m m} 
{
	\overset{\hom{∙}{#3}}{#2}
}
\NewDocumentCommand \opcartesianlift {o m m} 
{
	\overset{\hom{#3}{∙}}{#2}
}

\NewDocumentCommand \prohom {o m m}
{
	\makehom{⇸}
	\f[#1]{#2}{#3}
}

\NewDocumentCommand \het {o m m}
{
	\makehom{⤍}
	\f[#1]{#2}{#3}
}

\NewDocumentCommand \pathhom {o m m}
{
	\makehom{⇝}
	\f[#1]{#2}{#3}
}

\NewDocumentCommand \lifts {}
{⧄}
\NewDocumentCommand \orthogonal {s m m}
{
	\IfBooleanTF{#1}
		{{#2} \mathrel{\overline{\lifts}} {#3}}
		{{#2} \mathrel{\lifts} {#3}}
}

\NewDocumentCommand \homotopy {o m m}
{
	\IfNoValueTF{#1}
	{
		\makehom{↝}
		\f{#2}{#3}
	}
	{
		\makehom{\overset{#1} {↝}}
		\f{#2}{#3}
	}
}

\NewDocumentCommand \hcomp {O{1} m}
{
	\monoid
	{\mathbin{\natrec{#1}{\mathord{+}}}}
	{0 \ifthenelse{#1 > 1}{^{#1}}{}}
	\f{#2}
}

\NewDocumentCommand \doublehom {o m m m m} 
{
	\IfNoValueTF{#1}
	{\prescript{#4\vphantom{#3}}{#2\vphantom{#5}}{}{◇}{}^{#3\vphantom{#4}}_{#5\vphantom{#2}}}
	{{#1} \, \hombrackets{\prescript{#4\vphantom{#3}}{#2\vphantom{#5}}{}{◇}{}^{#3\vphantom{#4}}_{#5\vphantom{#2}}}}
}

\NewDocumentCommand \fix {o m}
{
	\mathop{\mathrm{fix} {\IfValueTF {#1}{_{#1}}{}}} \ifthenelse{\isempty{#2}}{}{({#2})}
}

\NewDocumentCommand \boundrel {o m m}
{
	{#2} \mathbin{⪯} {#3}
	\IfNoValueTF{#1}{}{ : #1}
}

\NewDocumentCommand \costaction{o m m} 
{
	{#2} \mathbin{+_{\IfNoValueTF{#1}{⋅}{#1}}} {#3}
}

\NewDocumentCommand \evaluatesto {o m m}  
{
	{#2} \mathrel{\mathord{⇓}{\IfNoValueTF{#1}{}{^{#1}}}} {#3}
}

\NewDocumentCommand \diag {} {∆}

\hypersetup
{
	pdfauthor = {Edward Morehouse} ,
	pdftitle = {2-Categories from a Gray Perspective}
}

\title{$2$-Categories from a Gray Perspective}

\author
{
	\small
	Edward Morehouse%
	\thanks
	{
	This research was supported by
	the ESF funded Estonian IT Academy research measure
	(project 2014-2020.4.05.19-0001).
	}
}
\affil
{
	\small
	Tallinn University of Technology
}
\date{}

\begin{document}
	\maketitle
%

\begin{abstract}

In this paper we present $2$-category theory from the perspective of Gray-categories
using the graphical calculus of separated surface diagrams.
As an extended example we consider cones and limits of $2$-functors.
Then we use the canonical adjunction between $2$-computads and $2$-categories
to interpret the comparison structure of lax functors
and extend the surface diagram calculus with compositor sheets
in order to represent and reason about them.

\end{abstract}

\begin{sect}{A Gloss on Natural Transformations} \label{section: natural transformations}
	\subimport*{parts/}{strict_naturality}
\end{sect}

\begin{sect}{Picturing Gray-Categories} \label{section: gray categories}
	The globular $2$-dimensional categorical structure formed by categories, functors and natural transformations is known as a $2$-category.
$2$-categories, in turn, are related by morphisms of three different dimensions that come in several variations.
The structure of $2$-categories and their hierarchy of morphisms was studied by Gray in \cite{gray-1974-formal_category_theory}.
This structure has come to be known as a Gray-category.
We build it up in stages.

\begin{definition}[sesquicategory]
	A \define{sesquicategory} $\cat{C}$ consists of a category of $0$-cells (objects) and $1$-cells (arrows)
	together with the following additional structure:
	\begin{description}
		\item[$2$-cells:]
			for each parallel pair of $1$-cells $f , f′ : \hom[\cat{C}]{A}{B}$
			a collection of $2$-cells\footnote
			{\
				For brevity we may omit any well-formed prefix of a boundary specification
				if it is either clear from context or irrelevant.
				This lets us write $\hom[\hom[\cat{C}]{A}{B}]{f}{f′}$
				as either $\hom[\hom[]{A}{B}]{f}{f′}$ or $\hom{f}{f′}$,
				dropping the brackets as well when we drop the entire prefix.
			}
			or ``disks''
			 $\hom[\hom[\cat{C}]{A}{B}]{f}{f′}$,
		\item[nullary $2$-cell composition:]
			for each $1$-cell $f : \hom{A}{B}$
			a $2$-cell $\comp{} \, f : \hom{f}{f}$,
		\item[binary $2$-cell composition:]
			for consecutive $2$-cells $α : \hom[\hom[]{A}{B}]{f}{f′}$ and $β : \hom[\hom[]{A}{B}]{f′}{f′′}$
			a $2$-cell $\comp{α , β} : \hom[\hom[]{A}{B}]{f}{f′′}$,
		\item[$2$-cell right-whiskering:]
			for a $2$-cell $α : \hom[\hom[]{A}{B}]{f}{f′}$ and an adjacent $1$-cell $g : \hom{B}{C}$
			a $2$-cell $\comp[2]{α , g} : \hom[\hom[]{A}{C}]{\comp{f , g}}{\comp{f′ , g}}$,
		\item[$2$-cell left-whiskering:]
			for a $1$-cell $f : \hom{A}{B}$ and an adjacent $2$-cell $β : \hom[\hom[]{B}{C}]{g}{g′}$ 
			a $2$-cell $\comp[2]{f , β} : \hom[\hom[]{A}{C}]{\comp{f , g}}{\comp{f , g′}}$.
	\end{description}
	
	In addition to the composition laws for the underlying category we have the following relations:
	\begin{description}
		\item[strict $2$-cell composition:]
			for $0$-cells $A , B : \cat{C}$,
			the $1$-cells and $2$-cells of $\cat{C}$ in $\hom[\cat{C}]{A}{B}$ form a category:
			\begin{equation} \label{strict 2-cell composition}
				\comp{\comp{} \, f , α}  =  α  =  \comp{α , \comp{} \, f′}
				\; , \;
				\comp{(\comp{α , β}) , γ}  =  \comp{α , (\comp{β , γ})}
			\end{equation}
		\item[$2$-whiskering functoriality:]
			the whiskering operations are functorial in their $2$-cell argument:
			\begin{equation} \label{2-whiskering functoriality}
				\begin{array}{c}
					\comp[2]{\comp{} \, f , g}  =  \comp{} (\comp{f , g})  =  \comp[2]{f , \comp{} \, g}
					\; ,
					\\
					\comp[2]{(\comp{α , β}) , g}  =  \comp{(\comp[2]{α , g}) , (\comp[2]{β , g})}
					\; , \;
					\comp[2]{f , (\comp{γ , δ})}  =  \comp{(\comp[2]{f , γ}) , (\comp[2]{f , δ})}
					\\
				\end{array}
			\end{equation}
		\item[$2$-whiskering algebra laws:]
			the whiskering operations are compatible with $1$-cell composition:
			\begin{equation} \label{2-cell whiskering algebra laws}
				\begin{array}{c}
					\comp[2]{α , \comp{} \, B}  =  α  =  \comp[2]{\comp{} \, A , α}
					\; ,
					\\
					\comp[2]{α , (\comp{g , h})}  =  \comp[2]{(\comp[2]{α , g}) , h}
					\; , \;
					\comp[2]{(\comp{f , g}) , γ}  =  \comp[2]{f , (\comp[2]{g , γ})}
				\end{array}
			\end{equation}
		\item[$2$-whiskering bialgebra law:]
			the whiskering operations are compatible with each other:
			\begin{equation} \label{2-cell whiskering bialgebra law}
				\comp[2]{(\comp[2]{f , β}) , h}  =  \comp[2]{f , (\comp[2]{β , h})}
			\end{equation}
	\end{description}
\end{definition}

In light of the strict associativity of $2$-cell composition
we adopt an unbracketed notation for it,
letting us write expressions like $\comp{α , β , γ}$ without ambiguity.
Because whiskering a $2$-cell with any number of $1$-cells on either side is also unambiguous
we adopt an unbracketed notation for this as well,
letting us write the unique $2$-cells appearing in the binary composition instances of the algebra laws and in the bialgebra law
as $\comp[2]{α , g , h}$, $\comp[2]{f , g , γ}$ and $\comp[2]{f , β , h}$.

Because of the unit laws for $2$-cell composition in \eqref{strict 2-cell composition} and for whiskering in \eqref{2-cell whiskering algebra laws}
the following string diagram represents a unique $2$-cell.
$$
	
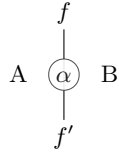
\begin{tikzpicture}[string diagram , x = {(12mm , 0mm)} , y = {(0mm , -12mm)} , baseline=(align.base)]
		\coordinate (f) at (1/2 , 1/2) ;
		\draw (f |- 0 , 0) to [out = south , in = north] node [fromabove] {f} coordinate [pos = 1/2] (a) node [tobelow] {f′} (f |- 1 , 1) ;
		\node [bead] at (a) (align) {α} ;
		\node at ( $(a) + (-1/2 , 0) $) {A} ;
		\node at ( $(a) + (1/2 , 0) $) {B} ;
	\end{tikzpicture}
$$
Because of the associative law for $2$-cell composition in \eqref{strict 2-cell composition}
the following string diagram represents a unique $2$-cell.
$$
	
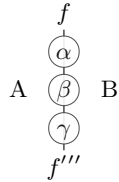
\begin{tikzpicture}[string diagram , x = {(12mm , 0mm)} , y = {(0mm , -16mm)} , baseline=(align.base)]
		\coordinate (f) at (1/2 , 1/2) ;
		\draw (f |- 0 , 0) to [out = south , in = north] node [fromabove] {f}
			coordinate [pos = 1/6] (a) coordinate [pos = 1/2] (b) coordinate [pos = 5/6] (c) 
			node [tobelow] {f′′′} (f |- 1 , 1) ;
		\node [bead] at (a) {α} ;
		\node [bead] at (b) (align) {β} ;
		\node [bead] at (c) {γ} ;
		\node at ( $(b) + (-1/2 , 0) $) {A} ;
		\node at ( $(b) + (1/2 , 0) $) {B} ;
	\end{tikzpicture}
$$
Because of whiskering functoriality \eqref{2-whiskering functoriality}
each of the following string diagrams represents a unique $2$-cell.
$$
	
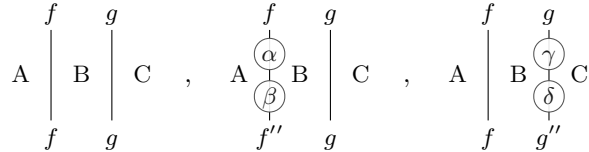
\begin{tikzpicture}[string diagram , x = {(16mm , 0mm)} , y = {(0mm , -12mm)} , baseline=(align.base)]
		\coordinate (f) at (1/4 , 1/2) ;
		\coordinate (g) at (3/4 , 1/2) ;
		\draw (f |- 0 , 0) to [out = south , in = north] node [fromabove] {f} coordinate [pos = 1/4] (a) coordinate [pos = 3/4] (b) node [tobelow] {f} (f |- 1 , 1) ;
		\draw (g |- 0 , 0) to [out = south , in = north] node [fromabove] {g} coordinate [pos = 1/4] (c) coordinate [pos = 3/4] (d) node [tobelow] {g} (g |- 1 , 1) ;
		\node at ($ (f) ! -1/2 ! (g) $) {A} ;
		\node at ($ (f) ! 1/2 ! (g) $) (align) {B} ;
		\node at ($ (f) ! 3/2 ! (g) $) {C} ;
	\end{tikzpicture}
	\quad , \quad
	\begin{tikzpicture}[string diagram , x = {(16mm , 0mm)} , y = {(0mm , -12mm)} , baseline=(align.base)]
		\coordinate (f) at (1/4 , 1/2) ;
		\coordinate (g) at (3/4 , 1/2) ;
		\draw (f |- 0 , 0) to [out = south , in = north] node [fromabove] {f} coordinate [pos = 1/4] (a) coordinate [pos = 3/4] (b) node [tobelow] {f′′} (f |- 1 , 1) ;
		\draw (g |- 0 , 0) to [out = south , in = north] node [fromabove] {g} coordinate [pos = 1/4] (c) coordinate [pos = 3/4] (d) node [tobelow] {g} (g |- 1 , 1) ;
		\node [bead] at (a) {α} ;
		\node [bead] at (b) {β} ;
		\node at ($ (f) ! -1/2 ! (g) $) {A} ;
		\node at ($ (f) ! 1/2 ! (g) $) (align) {B} ;
		\node at ($ (f) ! 3/2 ! (g) $) {C} ;
	\end{tikzpicture}
	\quad , \quad
	\begin{tikzpicture}[string diagram , x = {(16mm , 0mm)} , y = {(0mm , -12mm)} , baseline=(align.base)]
		\coordinate (f) at (1/4 , 1/2) ;
		\coordinate (g) at (3/4 , 1/2) ;
		\draw (f |- 0 , 0) to [out = south , in = north] node [fromabove] {f} coordinate [pos = 1/4] (a) coordinate [pos = 3/4] (b) node [tobelow] {f} (f |- 1 , 1) ;
		\draw (g |- 0 , 0) to [out = south , in = north] node [fromabove] {g} coordinate [pos = 1/4] (c) coordinate [pos = 3/4] (d) node [tobelow] {g′′} (g |- 1 , 1) ;
		\node [bead] at (c) {γ} ;
		\node [bead] at (d) {δ} ;
		\node at ($ (f) ! -1/2 ! (g) $) {A} ;
		\node at ($ (f) ! 1/2 ! (g) $) (align) {B} ;
		\node at ($ (f) ! 3/2 ! (g) $) {C} ;
	\end{tikzpicture}
$$
Because of the whiskering algebra laws for binary composites in \eqref{2-cell whiskering algebra laws}
and the biaglebra law \eqref{2-cell whiskering bialgebra law}
each of the following string diagrams represents a unique $2$-cell.
$$
	
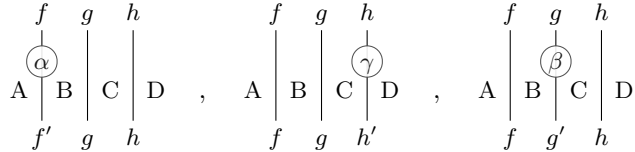
\begin{tikzpicture}[string diagram , x = {(24mm , 0mm)} , y = {(0mm , -12mm)} , baseline=(align.base)]
		\coordinate (f) at (1/4 , 2/3) ;
		\coordinate (g) at (2/4 , 2/3) ;
		\coordinate (h) at (3/4 , 2/3) ;
		\draw (f |- 0 , 0) to [out = south , in = north] node [fromabove] {f} coordinate [pos = 1/3] (a) node [tobelow] {f′} (f |- 1 , 1) ;
		\draw (g |- 0 , 0) to [out = south , in = north] node [fromabove] {g} coordinate [pos = 1/3] (b) node [tobelow] {g} (g |- 1 , 1) ;
		\draw (h |- 0 , 0) to [out = south , in = north] node [fromabove] {h} coordinate [pos = 1/3] (c) node [tobelow] {h} (h |- 1 , 1) ;
		\node [bead] at (a) {α} ;
		\node at ($ (f) ! -1/2 ! (g) $) {A} ;
		\node at ($ (f) ! 1/2 ! (g) $) {B} ;
		\node at ($ (g) ! 1/2 ! (h) $) {C} ;
		\node at ($ (g) ! 3/2 ! (h) $) (align) {D} ;
	\end{tikzpicture}
	\quad , \quad
	\begin{tikzpicture}[string diagram , x = {(24mm , 0mm)} , y = {(0mm , -12mm)} , baseline=(align.base)]
		\coordinate (f) at (1/4 , 2/3) ;
		\coordinate (g) at (2/4 , 2/3) ;
		\coordinate (h) at (3/4 , 2/3) ;
		\draw (f |- 0 , 0) to [out = south , in = north] node [fromabove] {f} coordinate [pos = 1/3] (a) node [tobelow] {f} (f |- 1 , 1) ;
		\draw (g |- 0 , 0) to [out = south , in = north] node [fromabove] {g} coordinate [pos = 1/3] (b) node [tobelow] {g} (g |- 1 , 1) ;
		\draw (h |- 0 , 0) to [out = south , in = north] node [fromabove] {h} coordinate [pos = 1/3] (c) node [tobelow] {h′} (h |- 1 , 1) ;
		\node [bead] at (c) {γ} ;
		\node at ($ (f) ! -1/2 ! (g) $) {A} ;
		\node at ($ (f) ! 1/2 ! (g) $) {B} ;
		\node at ($ (g) ! 1/2 ! (h) $) {C} ;
		\node at ($ (g) ! 3/2 ! (h) $) (align) {D} ;
	\end{tikzpicture}
	\quad , \quad
	\begin{tikzpicture}[string diagram , x = {(24mm , 0mm)} , y = {(0mm , -12mm)} , baseline=(align.base)]
		\coordinate (f) at (1/4 , 2/3) ;
		\coordinate (g) at (2/4 , 2/3) ;
		\coordinate (h) at (3/4 , 2/3) ;
		\draw (f |- 0 , 0) to [out = south , in = north] node [fromabove] {f} coordinate [pos = 1/3] (a) node [tobelow] {f} (f |- 1 , 1) ;
		\draw (g |- 0 , 0) to [out = south , in = north] node [fromabove] {g} coordinate [pos = 1/3] (b) node [tobelow] {g′} (g |- 1 , 1) ;
		\draw (h |- 0 , 0) to [out = south , in = north] node [fromabove] {h} coordinate [pos = 1/3] (c) node [tobelow] {h} (h |- 1 , 1) ;
		\node [bead] at (b) {β} ;
		\node at ($ (f) ! -1/2 ! (g) $) {A} ;
		\node at ($ (f) ! 1/2 ! (g) $) {B} ;
		\node at ($ (g) ! 1/2 ! (h) $) {C} ;
		\node at ($ (g) ! 3/2 ! (h) $) (align) {D} ;
	\end{tikzpicture}
$$

In a sesquicategory we can compose $2$-cells vertically, but only whisker them by $1$-cells horizontally.
A string diagram is called \define[ordered string diagram]{ordered}
if it can be decomposed into horizontal layers such that each layer contains at most one $2$-cell.
Two ordered string diagrams are equal just in case they have equal ordered layer decompositions.

For our purposes it will be convenient to characterize a $2$-category
as a sesquicategory with an additional relation that loosens the ordering restriction.

\begin{definition}[$2$-category]
	A \define[2-category]{$2$-category} $\cat{C}$ is a \refer{sesquicategory} satisfying the following additional relation:
	\begin{description}
		\item[interchange law:]
			for adjacent $2$-cells $α : \hom[\hom[\cat{C}]{A}{B}]{f}{f′}$ and $β : \hom[\hom[\cat{C}]{B}{C}]{g}{g′}$
			the relation:
			\begin{equation} \label{2-cell interchange}
				\comp{(\comp[2]{α , g}) , (\comp[2]{f′ , β})}
				\; = \;
				\comp{(\comp[2]{f , β}) , (\comp[2]{α , g′})}
			\end{equation}
	\end{description}
\end{definition}

This relation allows us to coherently define a \define{horizontal composition} operation on $2$-cells
$\comp[2]{α , β} : \hom[]{\comp{f , g}}{\comp{f′, g′}}$
that is equal to each of its perturbations into a composition of whiskerings:
\begin{equation} \label{horizontal composition}
	\begin{tikzpicture}[string diagram , x = {(16mm , 0mm)} , y = {(0mm , -12mm)} , baseline=(align.base)]
		\coordinate (f) at (1/4 , 1/2) ;
		\coordinate (g) at (3/4 , 1/2) ;
		\draw (f |- 0 , 0) to [out = south , in = north] node [fromabove] {f} coordinate [pos = 1/4] (a) node [tobelow] {f′} (f |- 1 , 1) ;
		\draw (g |- 0 , 0) to [out = south , in = north] node [fromabove] {g} coordinate [pos = 3/4] (b) node [tobelow] {g′} (g |- 1 , 1) ;
		\node [bead] at (a) {α} ;
		\node [bead] at (b) {β} ;
		\node at ($ (f) ! -1/2 ! (g) $) {A} ;
		\node at ($ (f) ! 1/2 ! (g) $) (align) {B} ;
		\node at ($ (f) ! 3/2 ! (g) $) {C} ;
	\end{tikzpicture}
	\quad ≕ \quad
	\begin{tikzpicture}[string diagram , x = {(16mm , 0mm)} , y = {(0mm , -12mm)} , baseline=(align.base)]
		\coordinate (f) at (1/4 , 1/2) ;
		\coordinate (g) at (3/4 , 1/2) ;
		\draw (f |- 0 , 0) to [out = south , in = north] node [fromabove] {f} coordinate [pos = 1/2] (a) node [tobelow] {f′} (f |- 1 , 1) ;
		\draw (g |- 0 , 0) to [out = south , in = north] node [fromabove] {g} coordinate [pos = 1/2] (b) node [tobelow] {g′} (g |- 1 , 1) ;
		\node [bead] at (a) {α} ;
		\node [bead] at (b) {β} ;
		\node at ($ (f) ! -1/2 ! (g) $) {A} ;
		\node at ($ (f) ! 1/2 ! (g) $) (align) {B} ;
		\node at ($ (f) ! 3/2 ! (g) $) {C} ;
	\end{tikzpicture}
	\quad ≔ \quad
	\begin{tikzpicture}[string diagram , x = {(16mm , 0mm)} , y = {(0mm , -12mm)} , baseline=(align.base)]
		\coordinate (f) at (1/4 , 1/2) ;
		\coordinate (g) at (3/4 , 1/2) ;
		\draw (f |- 0 , 0) to [out = south , in = north] node [fromabove] {f} coordinate [pos = 3/4] (a) node [tobelow] {f′} (f |- 1 , 1) ;
		\draw (g |- 0 , 0) to [out = south , in = north] node [fromabove] {g} coordinate [pos = 1/4] (b) node [tobelow] {g′} (g |- 1 , 1) ;
		\node [bead] at (a) {α} ;
		\node [bead] at (b) {β} ;
		\node at ($ (f) ! -1/2 ! (g) $) {A} ;
		\node at ($ (f) ! 1/2 ! (g) $) (align) {B} ;
		\node at ($ (f) ! 3/2 ! (g) $) {C} ;
	\end{tikzpicture}
\end{equation}
The horizontal composition of $2$-cells is associative
with units given by the vertical identity $2$-cell on the identity $1$-cell on each $0$-cell, $\comp[2]{} A ≔ \comp{} (\comp{} \, A)$.
\begin{equation} \label{2-category horizontal category}
	\begin{array}{l | l}
		\begin{relationalreasoning}
			\term[]
			{\comp[2]{(\comp[2]{α , β}) , γ}}
			\term[=]
			{\comp{(\comp[2]{(\comp[2]{α , β}) , h}) , (\comp[2]{(\comp{f′ , g′}) , γ})}}
			\term[=]
			{\comp{(\comp[2]{(\comp{(\comp[2]{α , g}) , (\comp[2]{f′ , β})}) , h}) , (\comp[2]{f′ , g′ , γ})}}
			\term[=]
			{\comp{(\comp[2]{α , g , h}) , (\comp[2]{f′ , β , h}) , (\comp[2]{f′ , g′ , γ})}}
			\term[=]
			{\comp{(\comp[2]{α  , g , h}) , (\comp[2]{f′ , (\comp{(\comp[2]{β , h}) , (\comp[2]{g′ , γ})})})}}
			\term[=]
			{\comp{(\comp[2]{α  , (\comp{g , h})}) , (\comp[2]{f′ , (\comp[2]{β , γ})})}}
			\term[=]
			{\comp[2]{α  , (\comp[2]{β , γ})}}
		\end{relationalreasoning}
		&
		\begin{relationalreasoning}
			\term[]
			{\comp[2]{\comp[2]{} A , α}}
			\term[=]
			{\comp{(\comp[2]{\comp[2]{} A , f}) , (\comp[2]{\comp{} \, A , α})}}
			\term[=]
			{\comp{\comp{} (\comp{\comp{} \, A , f}) , α}}
			\term[=]
			{\comp{\comp{} \, f , α}}
			\term[=]
			{α}
			\term[=]
			{\comp{α , \comp{} \, f′}}
			\term[=]
			{\comp{α , \comp{} (\comp{f′, \comp{} \, B})}}
			\term[=]
			{\comp{(\comp[2]{α , \comp{} \, B}) , (\comp[2]{f′ , \comp[2]{} B})}}
			\term[=]
			{\comp[2]{α , \comp[2]{} B}}
		\end{relationalreasoning}
		\\
	\end{array}
\end{equation}
The $0$-cells and $2$-cells of a $2$-category themselves form a category
so we adopt an unbracketed notation for horizontal composition of $2$-cells,
and each of the following string diagrams represents a unique $2$-cell.
$$
	\begin{tikzpicture}[string diagram , x = {(36mm , 0mm)} , y = {(0mm , -12mm)} , baseline=(align.base)]
		\coordinate (f) at (1/4 , 1/2) ;
		\coordinate (g) at (2/4 , 1/2) ;
		\coordinate (h) at (3/4 , 1/2) ;
		\draw (f |- 0 , 0) to [out = south , in = north] node [fromabove] {f} coordinate [pos = 1/2] (a) node [tobelow] {f′} (f |- 1 , 1) ;
		\draw (g |- 0 , 0) to [out = south , in = north] node [fromabove] {g} coordinate [pos = 1/2] (b) node [tobelow] {g′} (g |- 1 , 1) ;
		\draw (h |- 0 , 0) to [out = south , in = north] node [fromabove] {h} coordinate [pos = 1/2] (c) node [tobelow] {h′} (h |- 1 , 1) ;
		\node [bead] at (a) {α} ;
		\node [bead] at (b) {β} ;
		\node [bead] at (c) {γ} ;
		\node at ($ (f) ! -1/2 ! (g) $) {A} ;
		\node at ($ (f) ! 1/2 ! (g) $) {B} ;
		\node at ($ (g) ! 1/2 ! (h) $) {C} ;
		\node at ($ (g) ! 3/2 ! (h) $) (align) {D} ;
	\end{tikzpicture}
	\quad , \quad
	\begin{tikzpicture}[string diagram , x = {(12mm , 0mm)} , y = {(0mm , -12mm)} , baseline=(align.base)]
		\coordinate (f) at (1/2 , 1/2) ;
		\draw (f |- 0 , 0) to [out = south , in = north] node [fromabove] {f} coordinate [pos = 1/2] (a) node [tobelow] {f′} (f |- 1 , 1) ;
		\node [bead] at (a) (align) {α} ;
		\node at ( $(a) + (-1/2 , 0) $) {A} ;
		\node at ( $(a) + (1/2 , 0) $) {B} ;
	\end{tikzpicture}
$$
Moreover, the $2$-cells satisfy the \define{middle-four exchange law}
$\comp[2]{(\comp{α , γ}) , (\comp{β , δ})}  =  \comp{(\comp[2]{α , β}) , (\comp[2]{γ , δ})}$,
which is a $2$-dimensional associative law
asserting that the following string diagram represents a unique $2$-cell.
$$
	\begin{tikzpicture}[string diagram , x = {(16mm , 0mm)} , y = {(0mm , -12mm)} , baseline=(align.base)]
		\coordinate (f) at (1/4 , 1/2) ;
		\coordinate (g) at (3/4 , 1/2) ;
		\draw (f |- 0 , 0) to [out = south , in = north] node [fromabove] {f} coordinate [pos = 1/4] (a) coordinate [pos = 3/4] (c) node [tobelow] {f′′} (f |- 1 , 1) ;
		\draw (g |- 0 , 0) to [out = south , in = north] node [fromabove] {g} coordinate [pos = 1/4] (b) coordinate [pos = 3/4] (d) node [tobelow] {g′′} (g |- 1 , 1) ;
		\node [bead] at (a) {α} ;
		\node [bead] at (b) {β} ;
		\node [bead] at (c) {γ} ;
		\node [bead] at (d) {δ} ;
		\node at ($ (f) ! -1/2 ! (g) $) {A} ;
		\node at ($ (f) ! 1/2 ! (g) $) (align) {B} ;
		\node at ($ (f) ! 3/2 ! (g) $) {C} ;
	\end{tikzpicture}
$$
Because of the horizontal composition operation string diagrams of $2$-cells in $2$-categories need not be ordered.
Two such diagrams are equal just in case they are related by applications of the ternary interchange law \eqref{horizontal composition}.

\begin{definition}[box-category]
	A \define{box-category} $\cat{C}$ consists of a \refer{sesquicategory} together with the following additional structure:
	\begin{description}
		\item[$3$-cells:]
			for each parallel pair of $2$-cells, $α , α′ : \hom[\hom[\cat{C}]{A}{B}]{f}{f′}$
			a collection of $3$-cells or ``globes'' $\hom[\hom[\hom[\cat{C}]{A}{B}]{f}{f′}]{α}{α′}$,
		\item[nullary $3$-cell composition:]
			for each $2$-cell $α : \hom{f}{f′}$
			a $3$-cell $\comp{} \, α : \hom{α}{α}$,
		\item[binary $3$-cell composition:]
			for consecutive $3$-cells $Γ : \hom[\hom[]{f}{f′}]{α}{α′}$ and $Γ′ : \hom[\hom[]{f}{f′}]{α′}{α′′}$
			a $3$-cell $\comp{Γ , Γ′} : \hom[\hom[]{f}{f′}]{α}{α′′}$,
		\item[$3$-cell horizontal composition:]
			for adjacent $3$-cells $Γ : \hom[\hom[]{f}{f′}]{α}{α′}$ and $Δ : \hom[\hom[]{f′}{f′′}]{β}{β′}$
			a $3$-cell $\comp[2]{Γ , Δ} : \hom[\hom[]{f}{f′′}]{\comp{α , β}}{\comp{α′ , β′}}$,
		\item[$3$-cell right-whiskering:]
			for a $3$-cell $Γ : \hom[\hom[\hom[]{A}{B}]{f}{f′}]{α}{α′}$ and a proximate $1$-cell $g : \hom{B}{C}$
			a $3$-cell $\comp[3]{Γ , g} : \hom[\hom[\hom[]{A}{C}]{\comp{f , g}}{\comp{f′, g}}]{\comp[2]{α , g}}{\comp[2]{α′ , g}}$,
		\item[$3$-cell left-whiskering:]
			for a $1$-cell $f : \hom{A}{B}$ and a proximate $3$-cell $Φ : \hom[\hom[\hom[]{B}{C}]{g}{g′}]{γ}{γ′}$
			a $3$-cell $\comp[3]{f , Φ} : \hom[\hom[\hom[]{A}{C}]{\comp{f , g}}{\comp{f, g′}}]{\comp[2]{f , γ}}{\comp[2]{f , γ′}}$.
	\end{description}
	
	In addition to the composition laws for the underlying sesquicategory we have the following relations:
	\begin{description}
		\item[strict $3$-cell composition:]
			for parallel $1$-cells $f , f′ : \hom[\cat{C}]{A}{B}$,
			the $2$-cells and $3$-cells of $\cat{C}$ in  $\hom[\hom[\cat{C}]{A}{B}]{f}{f′}$ form a category:
			\begin{equation} \label{strict 3-cell composition}
				\comp{\comp{} \, α , Γ}  =  Γ  =  \comp{Γ , \comp{} \, α′}
				\; , \;
				\comp{(\comp{Γ , Φ}) , Ω}  =  \comp{Γ , (\comp{Φ , Ω})}
			\end{equation}
		\item[strict $3$-cell horizontal composition:]
			for $0$-cells $A , B : \cat{C}$,
			the $1$-cells and $3$-cells of $\cat{C}$ in  $\hom[\cat{C}]{A}{B}$ form a category:
			\begin{equation} \label{strict 3-cell horizontal composition}
				\comp[2]{\comp[2]{} f , Γ}  =  Γ  =  \comp[2]{Γ , \comp[2]{} f′}
				\; , \;
				\comp[2]{(\comp[2]{Γ , Δ}) , Λ}  =  \comp[2]{Γ , (\comp[2]{Δ , Λ})}
			\end{equation}
			\item[middle-four-exchange for $3$-cell composition:]
				there is a unique configuration of four $3$-cells with the same object boundary arranged in a square:
			\begin{equation} \label{3-cell exchange law}
					\comp[2]{(\comp{Γ , Φ}) , (\comp{Δ , Ψ})}
					=
					\comp{(\comp[2]{Γ , Δ}) , (\comp[2]{Φ , Ψ})}
				\end{equation}
		\item[$3$-whiskering functoriality:]
			the whiskering operations are functorial in their $3$-cell argument:
			\begin{equation} \label{3-whiskering functoriality}
				\begin{array}{c}
					\comp[3]{\comp{} \, α , g}  =  \comp{} (\comp[2]{α , g})
					\, , \,
					\comp[3]{f , \comp{} \, β}  =  \comp{} (\comp[2]{f , β})
					\, ,
					\\
					\comp[3]{(\comp{Γ , Γ′}) , g}  =  \comp{(\comp[3]{Γ , g}) , (\comp[3]{Γ′ , g})}
					\, , \,
					\comp[3]{f , (\comp{Φ , Φ′})}  =  \comp{(\comp[3]{f , Φ}) , (\comp[3]{f , Φ′})}
					\, ,
					\\
					\comp[3]{(\comp[2]{Γ , Δ}) , g}  =  \comp[2]{(\comp[3]{Γ , g}) , (\comp[3]{Δ , g})}
					\, , \,
					\comp[3]{f , (\comp[2]{Φ , Ψ})}  =  \comp[2]{(\comp[3]{f , Φ}) , (\comp[3]{f , Ψ})}
					\, \phantom{,}
					\\
				\end{array}
			\end{equation}
		\item[$3$-whiskering algebra laws:]
			the whiskering operations are compatible with $1$-cell composition:
			\begin{equation} \label{3-cell whiskering algebra laws}
				\begin{array}{c}
					\comp[3]{Γ , \comp{} \, B}  =  Γ  =  \comp[3]{\comp{} \, A , Γ}
					\; ,
					\\
					\comp[3]{Γ , (\comp{g , h})}  =  \comp[3]{(\comp[3]{Γ , g}) , h}
					\; , \;
					\comp[3]{(\comp{f , g}) , Λ}  =  \comp[3]{f , (\comp[3]{g , Λ})}
				\end{array}
			\end{equation}
		\item[$3$-whiskering bialgebra law:]
			the whiskering operations are compatible with each other:
			\begin{equation} \label{3-cell whiskering bialgebra law}
				\comp[3]{(\comp[3]{f , Δ}) , h}  =  \comp[3]{f , (\comp[3]{Δ , h})}
			\end{equation}
	\end{description}
\end{definition}

In light of the strict associativity of $3$-cell composition
along both $2$-cells and $1$-cells
we adopt unbracketed notations for these operations.
Because whiskering a $3$-cell with any number of $1$-cells on either side is also unambiguous
we adopt an unbracketed notation for this as well,
letting us write the unique $3$-cells appearing in the binary composition instances of the algebra laws and in the bialgebra law
as $\comp[3]{Γ , g , h}$, $\comp[3]{f , g , Λ}$ and $\comp[3]{f , Δ , h}$.

Because of the unit laws for $3$-cell composition in \eqref{strict 3-cell composition} and \eqref{strict 3-cell horizontal composition},
and for whiskering in \eqref{3-cell whiskering algebra laws}
the following surface diagram represents a unique $3$-cell.
$$

$$
Because of the associative laws for $3$-cell composition in \eqref{strict 3-cell composition} and \eqref{strict 3-cell horizontal composition}
each of the following surface diagrams represents a unique $3$-cell.
$$
	%
$$
Because of the middle-four-exchange law \eqref{3-cell exchange law}
the following surface diagram represents a unique $3$-cell.
$$
	%
$$
Note that the hom objects of a box category are themselves \refer[2-category]{$2$-categories}.

Because of whiskering functoriality \eqref{3-whiskering functoriality}
each of the following surface diagrams represents a unique $3$-cell.
$$
	%
	\quad \phantom{,}
$$
Because of the whiskering algebra laws for binary composites in \eqref{3-cell whiskering algebra laws}
and the biaglebra law \eqref{3-cell whiskering bialgebra law}
each of the following surface diagrams represents a unique $3$-cell.
$$
	%
$$

The top and bottom boundaries of a surface diagram in a box-category
are \refer[ordered string diagram]{ordered string diagrams} of its $0$-, $1$- and $2$-cells,
as are all horizontal slices that don't intersect a point representing a $3$-cell.
However, these surface diagrams are not themselves vertically ordered
because the hom objects are $2$-categories, which have horizontal composition and interchange.

We extended sesquicategories to $2$-categories by adding the property of $2$-cell interchange.
Gray-categories extend box-categories in a similar way,
but because box-categories have $3$-cells there is now ``room'' for the interchange of $2$-cells to be a structure rather than a property.

\begin{definition}[oplax Gray-category]
	An \define[Gray-category]{oplax Gray-category} $\cat{C}$ consists of a \refer{box-category} together with the following additional structure:
	\begin{description}
		\item[oplax interchangers:]
			for adjacent $2$-cells $α : \hom[\hom[]{A}{B}]{f}{f′}$ and $β : \hom[\hom[]{B}{C}]{g}{g′}$
			a $3$-cell
			\begin{equation} \label{gray interchanger}
				\interchange{α , β} : \hom[\hom[\hom[]{A}{C}]{\comp{f , g}}{\comp{f′ , g′}}]
				{\comp{(\comp[2]{α , g}) , (\comp[2]{f′ , β})}}
				{\comp{(\comp[2]{f  , β}) , (\comp[2]{α , g′})}}
				.
			\end{equation}
	\end{description}
	
	In addition to the composition laws for the underlying box-category we have the following relations:
	\begin{description}
		\item[composite $2$-cell interchangers:]
			interchangers respect  $2$-cell composition in each index:
			\begin{equation} \label{composite interchangers}
				\begin{array}{c}
					\interchange{\comp{} f , β}
					=
					\comp{} (\comp[2]{f , β})
					\; , \;
					\interchange{\comp{α , α′} , β}
					=
					\comp
					{
						(\comp[2]{\comp{} (\comp[2]{α , g}) , \interchange{α′ , β}}) ,
						(\comp[2]{\interchange{α , β} , \comp{} (\comp[2]{α′ , g′})})
					}
					\; ,
					\\
					\interchange{α , \comp{} g}
					=
					\comp{} (\comp[2]{α , g})
					\; , \;
					\interchange{α , \comp{β , β′}}
					=
					\comp
					{
						(\comp[2]{\interchange{α , β} , \comp{} (\comp[2]{f′ , β′})}) ,
						(\comp[2]{\comp{} (\comp[2]{f , β}) , \interchange{α , β′}})
					}
					\; \phantom{,}
					\\
				\end{array}
			\end{equation}
		\item[interchanger extremal whiskering:]
			interchangers are compatible with extremal whiskering in each index:
			\begin{equation} \label{interchanger extremal whiskering}
				\interchange{\comp[2]{f , β} , γ}
				=
				\comp[3]{f , \interchange{β , γ}}
				\; , \;
				\interchange{α , \comp[2]{β , h}}
				=
				\comp[3]{\interchange{α , β} , h}
			\end{equation}
		\item[interchanger medial whiskering:]
			interchangers are compatible with medial whiskering:
			\begin{equation} \label{interchanger medial whiskering}
				\interchange{\comp[2]{α , g} , γ}
				=
				\interchange{α , \comp[2]{g , γ}}
			\end{equation}
		\item[interchanger naturality:]
			interchangers are natural in each index,
			in the sense that for $3$-cells $Γ : \hom{α}{α′}$ and $Δ : \hom{β}{β′}$,
			\begin{equation} \label{interchanger naturality}
				\begin{array}{c}
					\comp{(\comp[2]{(\comp[3]{Γ , g}) , \comp{} (\comp[2]{f′ , β})}) , \interchange{α′ , β}}
					=
					\comp{\interchange{α , β} , (\comp[2]{\comp{} (\comp[2]{f , β}) , (\comp[3]{Γ , g′})})}
					\, ,
					\\
					\comp{(\comp[2]{\comp{} (\comp[2]{α , g}) , (\comp[3]{f′ , Δ})}) , \interchange{α , β′}}
					=
					\comp{\interchange{α , β} , (\comp[2]{(\comp[3]{f , Δ}) , \comp{} (\comp[2]{α , g′})})}
					\, \phantom{,}
					\\
				\end{array}
			\end{equation}
	\end{description}
\end{definition}

In light of their extremal and medial whiskering laws
we adopt an interchanger notation with arbitrarily many additional $1$-cell indices,
letting us write the unique $3$-cells expressed by those laws as
$\interchange{f , β , γ}$, $\interchange{α , β , h}$, and $\interchange{α , g , γ}$.

A \define{lax Gray-category} is one with interchangers oriented the other way.
A \define{pseudo Gray-category} has invertible interchangers, making it both lax and oplax.
A \define[3-category]{$3$-category} is a Gray-category whose interchangers are identity $3$-cells.

In diagrams interchanger $3$-cells are drawn as
lines on separated surfaces whose projections intersect at a crossing point:
\begin{equation} \label{interchanger elaboration}
	\begin{tikzpicture}[string diagram , x = {(16mm , -8mm)} , y = {(0mm , -16mm)} , z = {(20mm , 0mm)} , baseline=(align.base)]
		\node at (1/2 , 1/2 , 3/2) {C} ;
		\coordinate (sheet) at (1/2 , 1/2 , 1) ;
		\draw [sheet]
			($ (sheet) + (-1/2 , -1/2 , 0) $) to coordinate [pos = 3/4] (in)
			($ (sheet) + (1/2 , -1/2 , 0) $) to coordinate [pos = 3/4] (to)
			($ (sheet) + (1/2 , 1/2 , 0) $) to coordinate [pos = 3/4] (out)
			($ (sheet) + (-1/2 , 1/2 , 0) $) to coordinate [pos = 3/4] (from)
			cycle
		;
		\draw [on sheet] (in) to [out = south , in = north] node [fromabove] {β} node [tobelow] {β} (out) ;
		\node [anchor = west] at (from) {g} ;
		\node [anchor = east] at (to) {g′} ;
		\node at (1/2 , 1/2 , 1/2) {B} ;
		\coordinate (sheet) at (1/2 , 1/2 , 0) ;
		\draw [sheet]
			($ (sheet) + (-1/2 , -1/2 , 0) $) to coordinate [pos = 1/4] (in)
			($ (sheet) + (1/2 , -1/2 , 0) $) to coordinate [pos = 0] (to)
			($ (sheet) + (1/2 , 1/2 , 0) $) to coordinate [pos = 1/4] (out)
			($ (sheet) + (-1/2 , 1/2 , 0) $) to coordinate [pos = 0] (from)
			cycle
		;
		\draw [on sheet] (in) to [out = south , in = north] node [fromabove] {α} node [tobelow] {α} (out) ;
		\node [anchor = south west] at (from) {f} ;
		\node [anchor = north east] at (to) {f′} ;
		\node (align) at (1/2 , 1/2 , -1/2) {A} ;
	\end{tikzpicture}
	\quad ≔ \quad
	\begin{tikzpicture}[string diagram , x = {(18mm , -8mm)} , y = {(0mm , -18mm)} , z = {(6mm , 3mm)} , baseline=(align.base)]
		\coordinate (cell) at (1/2 , 1/2 , 1/2) ;
		\node at (1/2 , 1/2 , 3) {C} ;
		\coordinate (sheet) at (1/2 , 1/2 , 1) ;
		\draw [sheet]
			($ (sheet) + (-1/2 , -1/2 , 0) $) to coordinate [pos = 3/4] (in)
			($ (sheet) + (1/2 , -1/2 , 0) $) to coordinate [pos = 3/4] (to)
			($ (sheet) + (1/2 , 1/2 , 0) $) to coordinate [pos = 3/4] (out)
			($ (sheet) + (-1/2 , 1/2 , 0) $) to coordinate [pos = 3/4] (from)
			cycle
		;
		\draw [on sheet]
			(in) to [out = {180 + 15} , in = {90 - 15}] node [fromabove] {β}
			(cell) to [out = {270 - 15} , in = {90 - 0}] node [tobelow] {β}
			(out)
		;
		\node [anchor = west] at (from) {g} ;
		\node [anchor = east] at (to) {g′} ;
		\coordinate (sheet) at (1/2 , 1/2 , 0) ;
		\draw [sheet]
			($ (sheet) + (-1/2 , -1/2 , 0) $) to coordinate [pos = 1/4] (in)
			($ (sheet) + (1/2 , -1/2 , 0) $) to coordinate [pos = 0] (to)
			($ (sheet) + (1/2 , 1/2 , 0) $) to coordinate [pos = 1/4] (out)
			($ (sheet) + (-1/2 , 1/2 , 0) $) to coordinate [pos = 0] (from)
			cycle
		;
		\draw [on sheet]
			(in) to [out = {0 - 30} , in = {90 + 15}] node [fromabove] {α}
			(cell) to [out = {270 +15} , in = {90 + 0}] node [tobelow] {α}
			(out)
		;
		\node [anchor = south west] at (from) {f} ;
		\node [anchor = north east] at (to) {f′} ;
		\node [bead] at (cell) (align) {χ} ;
		\node at (1/2 , 1/2 , -2) {A} ;
	\end{tikzpicture}
\end{equation}
The topological properties of this representation turn out to be a good fit for the algebraic properties of interchangers.
However, this representation introduces the same potential ambiguity that we encountered in section \ref{section: natural transformations}.
We must still restrict the geometry so that the lines in the projection string diagram intersect \refer[transverse intersection]{transversely}.
We must also ensure that the line representing the second $2$-cell index of an interchanger
crosses the line representing its first $2$-cell index from the upper right to the lower left
only if the interchanger is oplax, and from the upper left to the lower right only if it is lax.

Because of the interchanger laws for composite $2$-cells \eqref{composite interchangers}
each of the following surface diagrams represents a unique $3$-cell.
$$

	\quad \phantom{,}
$$
By induction, interchangers for arbitrary $2$-cell composites in each index are well-defined.
Because of the interchanger extremal \eqref{interchanger extremal whiskering} and medial \eqref{interchanger medial whiskering} whiskering laws
each of the following surface diagrams represents a unique $3$-cell.
$$
	%
$$
Because of the interchanger naturality laws \eqref{interchanger naturality}
we have the following surface diagram equations
which, thanks to the separation convention, are homotopies.
$$
	%
	\quad \phantom{,}
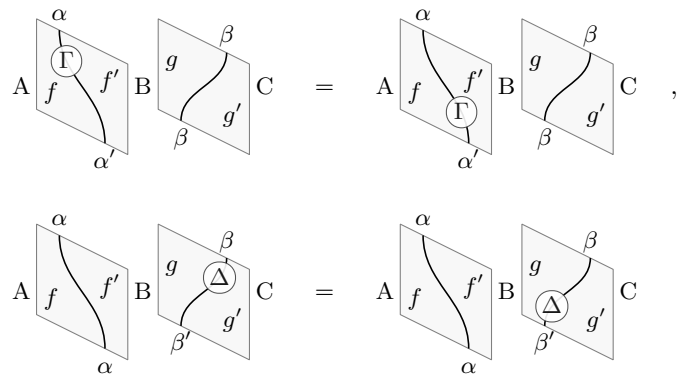
$$

A notable consequence of the composite $2$-cell interchanger laws and the interchanger naturality laws
is that permuting the composition order of any number of whiskered $2$-cells
using interchangers is coherent \cite{gray-1976-interchanger_coherence}.
The special case of reversing the order of three $2$-cells is known as the Yang-Baxter equation,
$$
	\begin{array}{c}
		\comp
		{
			(\comp[2]{\interchange{α , β , h} , \comp{} (\comp[2]{f′ , g′ , γ})}) ,
			(\comp[2]{\comp{} (\comp[2]{f , β , h}) , \interchange{α , g , γ}}) ,
			(\comp[2]{\interchange{f , β , γ} , \comp{} (\comp[2]{α , g′ , h′})})
		}
		\\
		=
		\\
		\comp
		{
			(\comp[2]{\comp{} (\comp[2]{α , g , h}) , \interchange{f′ , β , γ}}) ,
			(\comp[2]{\interchange{α , g , γ} , \comp{} (\comp[2]{f′ , β , h′})}) ,
			(\comp[2]{\comp{} (\comp[2]{f , g , γ}) , \interchange{α , β , h′}})
		}
		\\
	\end{array}
$$
which we can represent as the following equation between surface diagrams.
$$
	%
$$
These have projection string diagrams with the shapes of the two in \eqref{triple crossing} having only pairwise crossings.

Thus we could coherently define $n$-ary $2$-cell interchangers, such as $\interchange{α , β ,  γ}$,
and even $n$-ary $3$-cell interchangers, such as $\interchange{Γ , Δ , Λ}$.
The latter would be equal to each of the $54$ perturbations of the following surface diagram
containing only binary $2$-cell interchangers;
namely, those whose projection string diagrams have
a bead on one of the three segments of each of the three wires.
\begin{equation} \label{generalized interchangers}
	\begin{tikzpicture}[string diagram , x = {(12mm , -6mm)} , y = {(0mm , -12mm)} , z = {(16mm , 0mm)} , baseline=(align.base)]
		\node at (1/2 , 1/2 , 5/2) {D} ;
		\coordinate (sheet) at (1/2 , 1/2 , 2) ;
		\draw [sheet]
			($ (sheet) + (-1/2 , -1/2 , 0) $) to coordinate [pos = 3/4] (in)
			($ (sheet) + (1/2 , -1/2 , 0) $) to coordinate [pos = 3/4] (to)
			($ (sheet) + (1/2 , 1/2 , 0) $) to coordinate [pos = 3/4] (out)
			($ (sheet) + (-1/2 , 1/2 , 0) $) to coordinate [pos = 3/4] (from)
			cycle
		;
		\draw [on sheet] (in) to [out = south , in = north] node [fromabove] {γ} coordinate [pos = 1/2] (cell) node [tobelow] {γ′} (out) ;
		\node [bead] at (cell) {Λ} ;
		\node [anchor = north west] at (from) {h} ;
		\node [anchor = south east] at (to) {h′} ;
		\node at (1/2 , 1/2 , 3/2) {C} ;
		\coordinate (sheet) at (1/2 , 1/2 , 1) ;
		\draw [sheet]
			($ (sheet) + (-1/2 , -1/2 , 0) $) to coordinate [pos = 1/2] (in)
			($ (sheet) + (1/2 , -1/2 , 0) $) to coordinate [pos = 1/2] (to)
			($ (sheet) + (1/2 , 1/2 , 0) $) to coordinate [pos = 1/2] (out)
			($ (sheet) + (-1/2 , 1/2 , 0) $) to coordinate [pos = 1/2] (from)
			cycle
		;
		\draw [on sheet] (in) to [out = south , in = north] node [fromabove] {β} coordinate [pos = 1/2] (cell) node [tobelow] {β′}  (out) ;
		\node [bead] at (cell) {Δ} ;
		\node [anchor = west] at (from) {g} ;
		\node [anchor = east] at (to) {g′} ;
		\node at (1/2 , 1/2 , 1/2) {B} ;
		\coordinate (sheet) at (1/2 , 1/2 , 0) ;
		\draw [sheet]
			($ (sheet) + (-1/2 , -1/2 , 0) $) to coordinate [pos = 1/4] (in)
			($ (sheet) + (1/2 , -1/2 , 0) $) to coordinate [pos = 0] (to)
			($ (sheet) + (1/2 , 1/2 , 0) $) to coordinate [pos = 1/4] (out)
			($ (sheet) + (-1/2 , 1/2 , 0) $) to coordinate [pos = 0] (from)
			cycle
		;
		\draw [on sheet] (in) to [out = south , in = north] node [fromabove] {α} coordinate [pos = 1/2] (cell) node [tobelow] {α′} (out) ;
		\node [bead] at (cell) {Γ} ;
		\node [anchor = south west] at (from) {f} ;
		\node [anchor = north east] at (to) {f′} ;
		\node (align) at (1/2 , 1/2 , -1/2) {A} ;
	\end{tikzpicture}
\end{equation}
Even if we don't admit such \define[generalized interchanger]{generalized interchangers}
it is sometimes convenient to write an expression or draw a diagram involving them
as a shorthand for any/all of their equal admissible perturbations.

The top and bottom boundaries of a surface diagram in a Gray-category are ordered string diagrams of its $0$-, $1$- and $2$-cells,
as in the case for a box-category.
However, now there may be horizontal slices not intersecting a point representing a $3$-cell of the underlying box-category (left)
forming string diagrams that are not ordered (right).
By the transverse intersection requirement, these correspond to line crossings in the projection string diagram,
which we elaborate to interchanger $3$-cells as in \eqref{interchanger elaboration}.
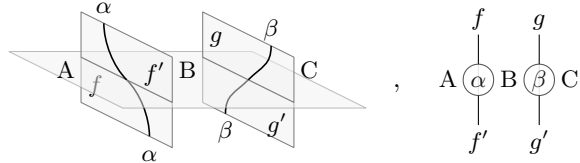
$$
	\begin{tikzpicture}[string diagram , x = {(12mm , -6mm)} , y = {(0mm , -12mm)} , z = {(16mm , 0mm)} , baseline=(align.base)]
		\node at (1/2 , 2/5 , 3/2) {C} ;
		\coordinate (sheet) at (1/2 , 1/2 , 1) ;
		\draw [sheet]
			($ (sheet) + (-1/2 , 0 , 0) $) to coordinate [pos = 1/2] (in)
			($ (sheet) + (1/2 , 0 , 0) $) to coordinate [pos = 1/2] (to)
			($ (sheet) + (1/2 , 1/2 , 0) $) to coordinate [pos = 3/4] (out)
			($ (sheet) + (-1/2 , 1/2 , 0) $) to coordinate [pos = 3/4] (from)
			cycle
		;
		\draw [on sheet] (in) to [out = south west , looseness = 1.2 , in = north] node [tobelow] {β} (out) ;
		\node [anchor = east] at (to) {g′} ;
		\node at (1/2 , 1/2 , 1/2) {} ;
		\coordinate (sheet) at (1/2 , 1/2 , 0) ;
		\draw [sheet]
			($ (sheet) + (-1/2 , 0 , 0) $) to coordinate [pos = 1/2] (in)
			($ (sheet) + (1/2 , 0 , 0) $) to coordinate [pos = 0] (to)
			($ (sheet) + (1/2 , 1/2 , 0) $) to coordinate [pos = 1/4] (out)
			($ (sheet) + (-1/2 , 1/2 , 0) $) to coordinate [pos = 0] (from)
			cycle
		;
		\draw [on sheet] (in) to [out = south east , looseness = 0.8 , in = north] node [tobelow] {α} (out) ;
		\node [anchor = south west] at (from) {f} ;
		\node (align) at (1/2 , 1/2 , -1/2) {} ;
		\draw [sheet , opacity = 0.5] (-1/8 , 1/2 , -1/2) -- (-1/8 , 1/2 , 3/2) -- (9/8 , 1/2 , 3/2) -- (9/8 , 1/2 , -1/2) -- cycle ;
		\node at (1/2 , 2/5 , 3/2) {} ;
		\coordinate (sheet) at (1/2 , 1/2 , 1) ;
		\draw [sheet]
			($ (sheet) + (-1/2 , -1/2 , 0) $) to coordinate [pos = 3/4] (in)
			($ (sheet) + (1/2 , -1/2 , 0) $) to coordinate [pos = 1/2] (to)
			($ (sheet) + (1/2 , 0 , 0) $) to coordinate [pos = 1/2] (out)
			($ (sheet) + (-1/2 , 0 , 0) $) to coordinate [pos = 1/2] (from)
			cycle
		;
		\draw [on sheet] (in) to [out = south , looseness = 1.2 , in = north east] node [fromabove] {β} (out) ;
		\node [anchor = west] at (from) {g} ;
		\node at (1/2 , 2/5 , 1/2) {B} ;
		\coordinate (sheet) at (1/2 , 1/2 , 0) ;
		\draw [sheet]
			($ (sheet) + (-1/2 , -1/2 , 0) $) to coordinate [pos = 1/4] (in)
			($ (sheet) + (1/2 , -1/2 , 0) $) to coordinate [pos = 0] (to)
			($ (sheet) + (1/2 , 0 , 0) $) to coordinate [pos = 1/2] (out)
			($ (sheet) + (-1/2 , 0 , 0) $) to coordinate [pos = 0] (from)
			cycle
		;
		\draw [on sheet] (in) to [out = south , looseness = 0.8 , in = north west] node [fromabove] {α} (out) ;
		\node [anchor = north east] at (to) {f′} ;
		\node at (1/2 , 2/5 , -1/2) {A} ;
	\end{tikzpicture}
	\quad , \quad
	\begin{tikzpicture}[string diagram , x = {(16mm , 0mm)} , y = {(0mm , -12mm)} , baseline=(align.base)]
		\coordinate (f) at (1/4 , 1/2) ;
		\coordinate (g) at (3/4 , 1/2) ;
		\draw (f |- 0 , 0) to [out = south , in = north] node [fromabove] {f} coordinate [pos = 1/2] (a) node [tobelow] {f′} (f |- 1 , 1) ;
		\draw (g |- 0 , 0) to [out = south , in = north] node [fromabove] {g} coordinate [pos = 1/2] (b) node [tobelow] {g′} (g |- 1 , 1) ;
		\node [bead] at (a) {α} ;
		\node [bead] at (b) {β} ;
		\node at ($ (f) ! -1/2 ! (g) $) {A} ;
		\node at ($ (f) ! 1/2 ! (g) $) {B} ;
		\node at ($ (f) ! 3/2 ! (g) $) {C} ;
		\node (align) at (f) {} ;
	\end{tikzpicture}
$$

\end{sect}

\begin{sect}{Diagram Semantics} \label{section: diagram semantics}
	Thus far we have used
string diagrams for $2$-dimensional categories and
surface diagrams for $3$-dimensional categories
informally as convenient notations for sometimes unwieldy expressions.
But in fact they have precise mathematical semantics.
Such diagrams can be understood as presentations of free categorical structures.
This was worked out in the case of string diagrams for monoidal categories in \cite{joyal-1991-geometry_of_tensor_calculus_1},
and extended to string diagrams for $2$-categories.
The case of surface diagrams for Gray Categories appears in \cite{hummon-2012-surface_diagrams}.

In the following we will not address the geometric aspects of dual diagrams.
While these are important, their development is rather complex, orthogonal to our present discussion, and can be found in the above references.
For our purposes it will suffice to consider diagrams as algebraic objects,
with equivalence given by an intuitive notion of configuration-preserving deformation.

A computad is an algebraic structure used to give a presentation for a free (generally higher-dimensional) category.
 The original construction was used by Street to give presentations for $2$- and $3$-categories \cite{street-1996-categorical_structures}.
Intuitively, an $n$-dimensional computad
consists of an $(n−1)$-dimensional computad specifying its lower-dimensional structure
together with a set of $n$-dimensional generator cells
and functions specifying their boundaries.
These boundaries are valued in the free $(n−1)$-dimensional category determined by the underlying $(n−1)$-dimensional computad,
which is built up from its generator cells using the composition operations of the category and quotiented by its relations.
We do not attempt to characterize the general construction here
but rather describe the computads and free categories that we will need.

\begin{itemize}
	\item
		A \define[0-computad]{$0$-computad} $G$
		consists of a set of $0$-generators $\celldim{0}{G}$.
	\item
		The \define[free 0-category]{free $0$-category} presented by a $0$-computad $\mathrm{diag}_0 G$
		has as $0$-cells the elements of the set $\celldim{0}{G}$, with no further structure or relations.
	\item
		A \define[1-computad]{$1$-computad}
		is a $0$-computad $G$
		together with a set of $1$-generators $\celldim{1}{G}$
		and boundary maps $\dom , \cod : \hom{\celldim{1}{G}}{\mathrm{diag}_0 G}$,
		in other words, a directed graph.
	\item
		The \define[free 1-category]{free $1$-category} presented by a $1$-computad $\mathrm{diag}_1 G$
		has as $1$-cells \define[dot diagram]{dot diagrams} built from the generators of $G$
		using nullary and binary composition and satisfying the categorical unit and associative laws
		so that each of the following dot diagrams represents a unique $1$-cell.
		$$
			\begin{tikzpicture}[string diagram , x = {(18mm , 0mm)} , y = {(0mm , -12mm)} , baseline=(align.base)]
				\coordinate (string) at (1/2 , 1/2) ;
				\draw
					(0 , 0 |- string) to node [fromleft] {A}
					coordinate [pos = 1/2] (f)
					node (align) [toright] {B} (1 , 1 |- string)
				;
				\node [bead] at (f) {f} ;
			\end{tikzpicture}
			\quad , \quad
			\begin{tikzpicture}[string diagram , x = {(36mm , 0mm)} , y = {(0mm , -12mm)} , baseline=(align.base)]
				\coordinate (string) at (1/2 , 1/2) ;
				\draw
					(0 , 0 |- string) to node [fromleft] {A}
					coordinate [pos = 1/6] (f) node [pos = 2/6 , auto] {B} coordinate [pos = 3/6] (g) node [pos = 4/6 , auto] {C} coordinate [pos = 5/6] (h)
					node (align) [toright] {D} (1 , 1 |- string)
				;
				\node [bead] at (f) {f} ;
				\node [bead] at (g) {g} ;
				\node [bead] at (h) {h} ;
			\end{tikzpicture}
		$$
	\item
		A \define{sesquicomputad}, which is also a \define[2-computad]{$2$-computad}, is a $1$-computad $G$
		together with a set of $2$-generators $\celldim{2}{G}$
		and boundary maps $\dom , \cod : \hom{\celldim{2}{G}}{\mathrm{diag}_1 G}$.
		These are required to satisfy globularity relations
		asserting that the boundaries of $2$-generators are parallel dot diagrams: 
		$\comp{\dom , \dom} = \comp{\cod , \dom}$ and $\comp{\dom , \cod} = \comp{\cod , \cod}$
		in $\hom{\celldim{2}{G}}{\mathrm{diag}_0 G}$,
		where the first boundary operator of each composite
		is that of the computad $G$ and the second is that of the category $\mathrm{diag}_1 G$.
	\item
		The \define{free sesquicategory} presented by a sesquicomputad $\mathrm{diag}_{\mathrm{ses}} G$
		has as $2$-cells \refer[ordered string diagram]{ordered string diagrams} built from the generators of $G$.
		The \define[free 2-category]{free $2$-category} presented by a $2$-computad $\mathrm{diag}_2 G$
		is similar, except that the string diagrams need not be ordered
		and may contain \refer[horizontal composition]{horizontal compositions} of $2$-cells.
	\item
		A \define{Gray-computad} is a sesquicomputad $G$
		together with a set of $3$-generators $\celldim{3}{G}$
		and boundary maps $\dom , \cod : \hom{\celldim{3}{G}}{\mathrm{diag}_{\mathrm{ses}} G}$.
		These are required to satisfy globularity relations
		asserting that the boundaries of $3$-generators are parallel ordered string diagrams:
		$\comp{\dom , \dom} = \comp{\cod , \dom}$ and $\comp{\dom , \cod} = \comp{\cod , \cod}$
		in $\hom{\celldim{3}{G}}{\mathrm{diag}_1 G}$.
	\item
		The \define{free Gray-category} presented by a Gray-computad $\mathrm{diag}_{\mathrm{Gray}} G$
		has as $3$-cells separated surface diagrams built from the generators of $G$.
		Interchanger $3$-cells are inferred from line crossings in the projection string diagrams.
		We want the result to be a valid diagram in a Gray-category
		so we require that in the projection string diagram:
		\begin{itemize}
			\item
				lines intersect \refer[transverse intersection]{transversely} (so that line intersections represent interchangers),
			\item
				such line crossings are oriented allowably (so that the interchangers are lax and/or oplax, as the case may be),
			\item
				only two lines cross at a point (so that interchangers have only two $2$-cell indices),
			\item
				points do not intersect points or lines (so that compositions involving generator $3$-cells along $0$-cells are whiskerings).
		\end{itemize}
		Such surface diagrams are sometimes called ``photogenic''.
		If we wish to allow $n$-ary $2$-cell interchangers and $3$-cell interchangers as in \eqref{generalized interchangers}
		then the last two conditions can be relaxed to requiring lines and line boundaries of points to obey the other conditions pairwise.
\end{itemize}

For $* ∈ \set{0 , 1 , \mathrm{ses} , 2 , \mathrm{Gray}}$ the collection of $*$-categories itself forms a category $*\Cat{Cat}$,
and the collection of $*$-computads forms a category $*\Cat{Ctd}$,
whose morphisms are the boundary-respecting functions between generator sets.
It is shown in \cite{hummon-2012-surface_diagrams} that in each case
the free $*$-category map is a functor $\mathrm{diag} : \hom{*\Cat{Ctd}}{*\Cat{Cat}}$
with a right adjoint $\mathrm{forget} : \hom{*\Cat{Cat}}{*\Cat{Ctd}}$
that regards each $n$-cell as an $n$-generator and forgets the composition structure.
The adjunction units have components known as \define[cone diagram]{cones},
which send $n$-generators to their singleton diagrams.
The adjunction counits have components known as \define[diagram evaluation]{evaluations},
which send diagrams of $n$-cells to their composites.

The right adjunction law says:
$$
	\begin{tikzpicture}[string diagram , x = {(24mm , 0mm)} , y = {(0mm , -16mm)} , baseline=(align.base)]
		\coordinate (unit) at (3/4 , 1/8) ;
		\coordinate (counit) at (1/4 , 7/8) ;
		\draw
			(0 , 0) to [out = south , in = west] node [fromabove] {\mathrm{forget}}
			(counit) to [out = east , looseness = 1.2 , in = west] node [online] (align) {\mathrm{diag}}
			(unit) to [out = east , in = north] node [tobelow] {\mathrm{forget}}
			(1 , 1)
		;
		\node [bead] at (unit) {\mathrm{cone}} ;
		\node [bead] at (counit) {\mathrm{eval}} ;
		\node at (unit |- counit) {*\Cat{Cat}} ;
		\node at (counit |- unit) {*\Cat{Ctd}} ;
	\end{tikzpicture}
	\quad = \quad
	\begin{tikzpicture}[string diagram , x = {(8mm , 0mm)} , y = {(0mm , -16mm)} , baseline=(align.base)]
		\draw
			(0 , 0) to [out = south , in = north] node [fromabove] {\mathrm{forget}} node [tobelow] {\mathrm{forget}}
			(1 , 1)
		;
		\node at (0 , 5/6) {*\Cat{Cat}} ;
		\node at (1 , 1/6) {*\Cat{Ctd}} ;
	\end{tikzpicture}
$$
Given a $*$-category $\cat{C}$,  the $*$-computad $\mathrm{forget} \, \cat{C}$ has as $n$-generators $n$-cells of $\cat{C}$,
while  the $*$-computad $\mathrm{forget} (\mathrm{diag} (\mathrm{forget} \, \cat{C}))$
has as $n$-generators diagrams of $n$-cells of $\cat{C}$.
The relation says that for any $n$-cell $φ$ of $\cat{C}$
if we form the singleton diagram comprising its cone
and evaluate it
then we get back $φ$.

The left adjunction law says:
$$
	\begin{tikzpicture}[string diagram , x = {(24mm , 0mm)} , y = {(0mm , -16mm)} , baseline=(align.base)]
		\coordinate (unit) at (1/4 , 1/8) ;
		\coordinate (counit) at (3/4 , 7/8) ;
		\draw
			(0 , 1) to [out = north , in = west] node [frombelow] {\mathrm{diag}}
			(unit) to [out = east , looseness = 1.2 , in = west] node [online] (align) {\mathrm{forget}}
			(counit) to [out = east , in = south] node [toabove] {\mathrm{diag}}
			(1 , 0)
		;
		\node [bead] at (unit) {\mathrm{cone}} ;
		\node [bead] at (counit) {\mathrm{eval}} ;
		\node at (counit |- unit) {*\Cat{Ctd}} ;
		\node at (unit |- counit) {*\Cat{Cat}} ;
	\end{tikzpicture}
	\quad = \quad
	\begin{tikzpicture}[string diagram , x = {(8mm , 0mm)} , y = {(0mm , -16mm)} , baseline=(align.base)]
		\draw
			(1 , 0) to [out = south , in = north] node [fromabove] {\mathrm{diag}} node [tobelow] {\mathrm{diag}}
			(0 , 1)
		;
		\node at (0 , 1/6) {*\Cat{Ctd}} ;
		\node at (1 , 5/6) {*\Cat{Cat}} ;
	\end{tikzpicture}
$$
Given a $*$-computad $G$, the $*$-category $\mathrm{diag}\,  G$ has as $n$-cells diagrams of $n$-generators of $G$,
while the $*$-category $\mathrm{diag} (\mathrm{forget} (\mathrm{diag} \, G))$
has as $n$-cells diagrams of diagrams of $n$-generators of $G$.
The relation says that for any diagram $δ$ of $G$'s $n$-generators
if we form the diagram of diagrams resulting from coning each $n$-generator in place
and evaluate the resulting diagram of diagrams to a diagram
then we get back $δ$.

\end{sect}

\begin{sect}{Morphisms of $2$-Categories} \label{section: 2-category theory}
	By combining Gray's classical results about the structure of the hierarchy of morphisms between $2$-categories
\cite{gray-1974-formal_category_theory}
with Hummon's results about surface diagrams for Gray-categories
\cite{hummon-2012-surface_diagrams}
we obtain an intuitive graphical calculus for $2$-category theory.
We can recover the conventional component-based presentation
by using global elements and projection.
This works because projection flattens a surface diagram in a Gray-category
into a string diagram in the hom $2$-category determined by its $0$-dimensional boundary.
By choosing the domain of this boundary to be the singleton $2$-category $\cat{1}$
we obtain a string diagram in the codomain $2$-category.

The simplest sort of $1$-dimensional morphism between $2$-categories is known as a $2$-functor \cite{kelly-1974-2_categories}.

\begin{definition}[$2$-functor between $2$-categories]
	A \define[2-functor]{$2$-functor} between \refer[2-category]{$2$-categories} $F : \hom{\cat{C}}{\cat{D}}$
	consists of the following structure:
	\begin{description}
		\item[object map:]
			a function on objects $\celldim{0}{F} : \hom{\celldim{0}{\cat{C}}}{\celldim{0}{\cat{D}}}$,
		\item[hom functors:]
			for each pair of objects $A , B : \cat{C}$ a local functor
			$\celldim{1}{F}^{\tuple*{A , B}} : \hom{\hom[\cat{C}]{A}{B}}{\hom[\cat{D}]{\celldim{0}{F} A}{\celldim{0}{F} B}}$.
	\end{description}
	The hom functors are required to strictly preserve the horizontal composition structure is the following sense\footnote
	{Except for the sake of emphasis, we will follow the convention of omitting the super- and subscript annotations on the constituent maps.}.
	\begin{description}
		\item[arrow composition preservation:]
			for each $\cat{C}$-object $A$
			and consecutive $\cat{C}$-arrows $f : \hom{A}{B}$ and $g : \hom{B}{C}$ we have
			\begin{equation}
				F (\comp{} A)  =  \comp{} (F A)
				\quad \text{and} \quad
				F (\comp{f , g})  =  \comp{F f , F g}
			\end{equation}
		\item[disk horizontal composition preservation:]
			for adjacent $\cat{C}$-disks $φ : \hom[\hom[]{A}{B}]{f}{f′}$ and $ψ : \hom[\hom[]{B}{C}]{g}{g′}$ we have
			\begin{equation}
				F (\comp[2]{φ , ψ})  =  \comp[2]{F φ , F ψ}
			\end{equation}
	\end{description}
\end{definition}

Note that $F$ necessarily preserves the vertical composition structure of disks
because the $\celldim{1}{F}$s are functors on the hom categories.
Thus by \eqref{2-category horizontal category} it preserves the horizontal composition units as well:
$$
	F (\comp[2]{} \, A)  ≔
	F (\comp{} (\comp{} \, A))  =
	\comp{} (F (\comp{} \, A))  =
	\comp{} (\comp{} (F A))  ≕
	\comp[2] {} (F A)
	.
$$
In diagrams we represent a $2$-functor as a surface separating the volumes representing its boundary $2$-categories.
This representation preserves all composition structure in the domain $2$-category
because given a string diagram $δ$ in $\cat{C}$
the string diagram $F δ$ in $\cat{D}$
looks just like $δ$, but with ``$F$'' prefixed to each label.
For example,
$$
	
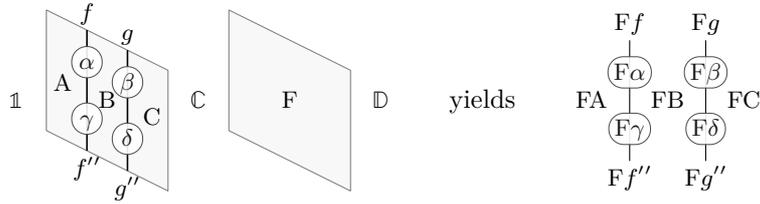
\begin{tikzpicture}[string diagram , x = {(16mm , -8mm)} , y = {(0mm , -16mm)} , z = {(24mm , 0mm)} , baseline=(align.base)]
		\node at (1/2 , 1/2 , 3/2) {\cat{D}} ;
		\coordinate (sheet) at (1/2 , 1/2 , 1) ;
		\draw [sheet]
			($ (sheet) + (-1/2 , -1/2 , 0) $) to
			($ (sheet) + (1/2 , -1/2 , 0) $) to
			($ (sheet) + (1/2 , 1/2 , 0) $) to
			($ (sheet) + (-1/2 , 1/2 , 0) $) to
			cycle
		;
		\node at (sheet) {F} ;
		\node at (1/2 , 1/2 , 1/2) {\cat{C}} ;
		\coordinate (sheet) at (1/2 , 1/2 , 0) ;
		\draw [sheet]
			($ (sheet) + (-1/2 , -1/2 , 0) $) to coordinate [pos = 1/3] (in 1) coordinate [pos = 2/3] (in 2)
			($ (sheet) + (1/2 , -1/2 , 0) $) to coordinate [pos = 1/4] (to)
			($ (sheet) + (1/2 , 1/2 , 0) $) to coordinate [pos = 1/3] (out 2) coordinate [pos = 2/3] (out 1)
			($ (sheet) + (-1/2 , 1/2 , 0) $) to coordinate [pos = 1/4] (from)
			cycle
		;
		\draw [on sheet]
			(in 1) to [out = south , in = north] node [fromabove] {f}
				coordinate [pos = 1/4] (cell 11) coordinate [pos = 3/4] (cell 12) node [tobelow] {f′′}
			(out 1)
		;
		\draw [on sheet]
			(in 2) to [out = south , in = north] node [fromabove] {g}
				coordinate [pos = 1/4] (cell 21) coordinate [pos = 3/4] (cell 22) node [tobelow] {g′′}
			(out 2)
		;
		\node [bead] at (cell 11) {α} ;
		\node [bead] at (cell 12) {γ} ;
		\node [bead] at (cell 21) {β} ;
		\node [bead] at (cell 22) {δ} ;
		\node [anchor = south west] at (from) {A} ;
		\node [anchor = center] at ($ (from) ! 1/2 ! (to) $) {B} ;
		\node [anchor = north east] at (to) {C} ;
		\node (align) at (1/2 , 1/2 , -1/2) {\cat{1}} ;
	\end{tikzpicture}
	\qquad \text{yields} \qquad
	\begin{tikzpicture}[string diagram , x = {(20mm , 0mm)} , y = {(0mm , -16mm)} , baseline=(align.base)]
		\coordinate (f) at (1/4 , 1/2) ;
		\coordinate (g) at (3/4 , 1/2) ;
		\draw (f |- 0 , 0) to [out = south , in = north] node [fromabove] {F f} coordinate [pos = 1/4] (a) coordinate [pos = 3/4] (c) node [tobelow] {F f′′} (f |- 1 , 1) ;
		\draw (g |- 0 , 0) to [out = south , in = north] node [fromabove] {F g} coordinate [pos = 1/4] (b) coordinate [pos = 3/4] (d) node [tobelow] {F g′′} (g |- 1 , 1) ;
		\node [bead] at (a) {F α} ;
		\node [bead] at (b) {F β} ;
		\node [bead] at (c) {F γ} ;
		\node [bead] at (d) {F δ} ;
		\node at ($ (f) ! -1/2 ! (g) $) {F A} ;
		\node at ($ (f) ! 1/2 ! (g) $) (align) {F B} ;
		\node at ($ (f) ! 3/2 ! (g) $) {F C} ;
	\end{tikzpicture}
$$

Like a \refer{natural transformation} between $1$-functors between $1$-categories,
a transformation between $2$-functors between $2$-categories
\cite{kelly-1974-2_categories, leinster-1998-bicategories}
has component arrows for objects.
But because there is now ``room'' for $2$-dimensional structure,
we can categorify the property of being natural for arrows into the structure of component disks.

\begin{definition}[oplax transformation between $2$-functors]
	An \define[oplax transformation between 2-functors]{oplax transformation} between $2$-functors between $2$-categories
	$α : \hom[\hom[]{\cat{C}}{\cat{D}}]{F}{G}$
	consists of the following data:
	\begin{description}
		\item[object-component arrows:]
			for each object of the domain $2$-category $A : \cat{C}$
			an arrow of the codomain $2$-category $α A : \hom[\cat{D}]{F A}{G A}$,
		\item[arrow-component disks:]
			for each arrow of the domain $2$-category $f : \hom[\cat{C}]{A}{B}$
			a disk of the codomain $2$-category $α f : \hom[\hom[\cat{D}]{F A}{G B}]{\comp{F f , α B}}{\comp{α A , G f}}$.
	\end{description}
	This data is required to satisfy the following relations:
	\begin{description}
		\item[compatibility with arrow composition:]
			for each $\cat{C}$-object $A$ and consecutive $\cat{C}$-arrows $f : \hom{A}{B}$ and $g :  \hom{B}{C}$
			we have
			\begin{equation} \label{transformation arrow composition compatibility}
				α (\comp{} \, A) = \comp{} (α A)
				\quad \text{and} \quad
				α (\comp{f , g}) = \comp{(\comp[2]{\comp{} (F f) , α g}) , (\comp[2]{α f , \comp{} (G g)})}
			\end{equation}
		\item[naturality for disks:]
			for each $\cat{C}$-disk $φ : \hom[\hom[]{A}{B}]{f}{f′}$
			we have
			\begin{equation} \label{transformation naturality for disks}
				\comp{(\comp[2]{F φ , \comp{} (α B)}) , α f′} = \comp{α f , (\comp[2]{\comp{} (α A) , G φ})}
			\end{equation}
	\end{description}
\end{definition}

A \define[lax transformation between 2-functors]{lax transformation}
has its arrow-component disks oriented the other way and satisfies suitably dualized relations.
A \define[pseudo transformation between 2-functors]{pseudo transformation}
has invertible arrow-component disks that make it lax one way and oplax the other.
A pseudo transformation is \define[strict transformation between 2-functors]{strict}
if its arrow-component disks are identities.

Transformations are $2$-cells between parallel $1$-cells,
so we represent them in surface diagrams as lines separating co-bounded surfaces.
Their arrow-component disks $α f$ are the projections of Gray-interchangers $\interchange{f , α}$.
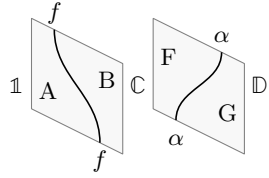
$$
	\begin{tikzpicture}[string diagram , x = {(12mm , -6mm)} , y = {(0mm , -12mm)} , z = {(16mm , 0mm)} , baseline=(align.base)]
		\node at (1/2 , 1/2 , 3/2) {\cat{D}} ;
		\coordinate (sheet) at (1/2 , 1/2 , 1) ;
		\draw [sheet]
			($ (sheet) + (-1/2 , -1/2 , 0) $) to coordinate [pos = 3/4] (in)
			($ (sheet) + (1/2 , -1/2 , 0) $) to coordinate [pos = 3/4] (to)
			($ (sheet) + (1/2 , 1/2 , 0) $) to coordinate [pos = 3/4] (out)
			($ (sheet) + (-1/2 , 1/2 , 0) $) to coordinate [pos = 3/4] (from)
			cycle
		;
		\draw [on sheet] (in) to [out = south , in = north] node [fromabove] {α} node [tobelow] {α} (out) ;
		\node [anchor = north west] at (from) {F} ;
		\node [anchor = south east] at (to) {G} ;
		\node at (1/2 , 1/2 , 1/2) {\cat{C}} ;
		\coordinate (sheet) at (1/2 , 1/2 , 0) ;
		\draw [sheet]
			($ (sheet) + (-1/2 , -1/2 , 0) $) to coordinate [pos = 1/4] (in)
			($ (sheet) + (1/2 , -1/2 , 0) $) to coordinate [pos = 0] (to)
			($ (sheet) + (1/2 , 1/2 , 0) $) to coordinate [pos = 1/4] (out)
			($ (sheet) + (-1/2 , 1/2 , 0) $) to coordinate [pos = 0] (from)
			cycle
		;
		\draw [on sheet] (in) to [out = south , in = north] node [fromabove] {f} node [tobelow] {f} (out) ;
		\node [anchor = south west] at (from) {A} ;
		\node [anchor = north east] at (to) {B} ;
		\node (align) at (1/2 , 1/2 , -1/2) {\cat{1}} ;
	\end{tikzpicture}
$$
Their compatibility with arrow composition \eqref{transformation arrow composition compatibility}
is the respect by interchangers for composition in the first index \eqref{composite interchangers},
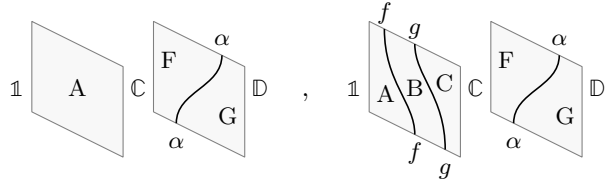
$$
	\begin{tikzpicture}[string diagram , x = {(12mm , -6mm)} , y = {(0mm , -12mm)} , z = {(16mm , 0mm)} , baseline=(align.base)]
		\node at (1/2 , 1/2 , 3/2) {\cat{D}} ;
		\coordinate (sheet) at (1/2 , 1/2 , 1) ;
		\draw [sheet]
			($ (sheet) + (-1/2 , -1/2 , 0) $) to coordinate [pos = 3/4] (in)
			($ (sheet) + (1/2 , -1/2 , 0) $) to coordinate [pos = 3/4] (to)
			($ (sheet) + (1/2 , 1/2 , 0) $) to coordinate [pos = 3/4] (out)
			($ (sheet) + (-1/2 , 1/2 , 0) $) to coordinate [pos = 3/4] (from)
			cycle
		;
		\draw [on sheet] (in) to [out = south , in = north] node [fromabove] {α} node [tobelow] {α} (out) ;
		\node [anchor = north west] at (from) {F} ;
		\node [anchor = south east] at (to) {G} ;
		\node at (1/2 , 1/2 , 1/2) {\cat{C}} ;
		\coordinate (sheet) at (1/2 , 1/2 , 0) ;
		\draw [sheet]
			($ (sheet) + (-1/2 , -1/2 , 0) $) to
			($ (sheet) + (1/2 , -1/2 , 0) $) to
			($ (sheet) + (1/2 , 1/2 , 0) $) to
			($ (sheet) + (-1/2 , 1/2 , 0) $) to
			cycle
		;
		\node at (sheet) {A} ;
		\node (align) at (1/2 , 1/2 , -1/2) {\cat{1}} ;
	\end{tikzpicture}
	\quad , \quad
	\begin{tikzpicture}[string diagram , x = {(12mm , -6mm)} , y = {(0mm , -12mm)} , z = {(16mm , 0mm)} , baseline=(align.base)]
		\node at (1/2 , 1/2 , 3/2) {\cat{D}} ;
		\coordinate (sheet) at (1/2 , 1/2 , 1) ;
		\draw [sheet]
			($ (sheet) + (-1/2 , -1/2 , 0) $) to coordinate [pos = 3/4] (in)
			($ (sheet) + (1/2 , -1/2 , 0) $) to coordinate [pos = 3/4] (to)
			($ (sheet) + (1/2 , 1/2 , 0) $) to coordinate [pos = 3/4] (out)
			($ (sheet) + (-1/2 , 1/2 , 0) $) to coordinate [pos = 3/4] (from)
			cycle
		;
		\draw [on sheet] (in) to [out = south , in = north] node [fromabove] {α} node [tobelow] {α} (out) ;
		\node [anchor = north west] at (from) {F} ;
		\node [anchor = south east] at (to) {G} ;
		\node at (1/2 , 1/2 , 1/2) {\cat{C}} ;
		\coordinate (sheet) at (1/2 , 1/2 , 0) ;
		\draw [sheet]
			($ (sheet) + (-1/2 , -1/2 , 0) $) to coordinate [pos = 1/6] (in 1) coordinate [pos = 3/6] (in 2)
			($ (sheet) + (1/2 , -1/2 , 0) $) to coordinate [pos = 0] (to)
			($ (sheet) + (1/2 , 1/2 , 0) $) to coordinate [pos = 1/6] (out 2) coordinate [pos = 3/6] (out 1)
			($ (sheet) + (-1/2 , 1/2 , 0) $) to coordinate [pos = 0] (from)
			cycle
		;
		\draw [on sheet] (in 1) to [out = south , out looseness = 1.25 , in looseness = 0.75 , in = north] node [fromabove] {f} node [tobelow] {f} (out 1) ;
		\draw [on sheet] (in 2) to [out = south , out looseness = 0.75 , in looseness = 1.25 , in = north] node [fromabove] {g} node [tobelow] {g} (out 2) ;
		\node [anchor = south west] at (from) {A} ;
		\node [anchor = center] at ($ (from) ! 1/2 ! (to) $) {B} ;
		\node [anchor = north east] at (to) {C} ;
		\node (align) at (1/2 , 1/2 , -1/2) {\cat{1}} ;
	\end{tikzpicture}
$$
and their naturality for disks \eqref{transformation naturality for disks}
is interchanger naturality in the first index \eqref{interchanger naturality}.
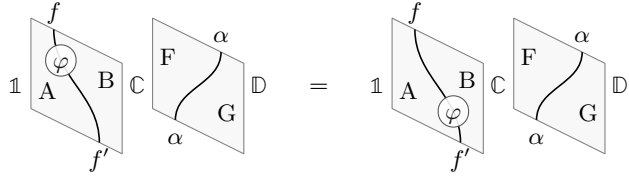
$$
	\begin{tikzpicture}[string diagram , x = {(12mm , -6mm)} , y = {(0mm , -12mm)} , z = {(16mm , 0mm)} , baseline=(align.base)]
		\node at (1/2 , 1/2 , 3/2) {\cat{D}} ;
		\coordinate (sheet) at (1/2 , 1/2 , 1) ;
		\draw [sheet]
			($ (sheet) + (-1/2 , -1/2 , 0) $) to coordinate [pos = 3/4] (in)
			($ (sheet) + (1/2 , -1/2 , 0) $) to coordinate [pos = 3/4] (to)
			($ (sheet) + (1/2 , 1/2 , 0) $) to coordinate [pos = 3/4] (out)
			($ (sheet) + (-1/2 , 1/2 , 0) $) to coordinate [pos = 3/4] (from)
			cycle
		;
		\draw [on sheet] (in) to [out = south , in = north] node [fromabove] {α} coordinate [pos = 1/2] (cell) node [tobelow] {α} (out) ;
		\node at (cell) {} ;
		\node [anchor = north west] at (from) {F} ;
		\node [anchor = south east] at (to) {G} ;
		\node at (1/2 , 1/2 , 1/2) {\cat{C}} ;
		\coordinate (sheet) at (1/2 , 1/2 , 0) ;
		\draw [sheet]
			($ (sheet) + (-1/2 , -1/2 , 0) $) to coordinate [pos = 1/4] (in)
			($ (sheet) + (1/2 , -1/2 , 0) $) to coordinate [pos = 0] (to)
			($ (sheet) + (1/2 , 1/2 , 0) $) to coordinate [pos = 1/4] (out)
			($ (sheet) + (-1/2 , 1/2 , 0) $) to coordinate [pos = 0] (from)
			cycle
		;
		\draw [on sheet] (in) to [out = south , in = north] node [fromabove] {f} coordinate [pos = 1/4] (cell) node [tobelow] {f′} (out) ;
		\node [bead] at (cell) {φ} ;
		\node [anchor = south west] at (from) {A} ;
		\node [anchor = north east] at (to) {B} ;
		\node (align) at (1/2 , 1/2 , -1/2) {\cat{1}} ;
	\end{tikzpicture}
	\quad = \quad
	\begin{tikzpicture}[string diagram , x = {(12mm , -6mm)} , y = {(0mm , -12mm)} , z = {(16mm , 0mm)} , baseline=(align.base)]
		\node at (1/2 , 1/2 , 3/2) {\cat{D}} ;
		\coordinate (sheet) at (1/2 , 1/2 , 1) ;
		\draw [sheet]
			($ (sheet) + (-1/2 , -1/2 , 0) $) to coordinate [pos = 3/4] (in)
			($ (sheet) + (1/2 , -1/2 , 0) $) to coordinate [pos = 3/4] (to)
			($ (sheet) + (1/2 , 1/2 , 0) $) to coordinate [pos = 3/4] (out)
			($ (sheet) + (-1/2 , 1/2 , 0) $) to coordinate [pos = 3/4] (from)
			cycle
		;
		\draw [on sheet] (in) to [out = south , in = north] node [fromabove] {α} coordinate [pos = 1/2] (cell) node [tobelow] {α} (out) ;
		\node at (cell) {} ;
		\node [anchor = north west] at (from) {F} ;
		\node [anchor = south east] at (to) {G} ;
		\node at (1/2 , 1/2 , 1/2) {\cat{C}} ;
		\coordinate (sheet) at (1/2 , 1/2 , 0) ;
		\draw [sheet]
			($ (sheet) + (-1/2 , -1/2 , 0) $) to coordinate [pos = 1/4] (in)
			($ (sheet) + (1/2 , -1/2 , 0) $) to coordinate [pos = 0] (to)
			($ (sheet) + (1/2 , 1/2 , 0) $) to coordinate [pos = 1/4] (out)
			($ (sheet) + (-1/2 , 1/2 , 0) $) to coordinate [pos = 0] (from)
			cycle
		;
		\draw [on sheet] (in) to [out = south , in = north] node [fromabove] {f} coordinate [pos = 3/4] (cell) node [tobelow] {f′} (out) ;
		\node [bead] at (cell) {φ} ;
		\node [anchor = south west] at (from) {A} ;
		\node [anchor = north east] at (to) {B} ;
		\node (align) at (1/2 , 1/2 , -1/2) {\cat{1}} ;
	\end{tikzpicture}
$$

The \emph{naturality} of transformations between $2$-functors
corresponds to the \emph{parametricity} of their component structure
in the sense that it is independent of the domain $2$-category.
This is what makes the separated surface notation so useful:
with no connection between the surfaces
there is no obstruction to moving a disk in the domain $2$-category
past the crossing in the projection string diagram representing an arrow-component disk,
so naturality is just a homotopy.

In a setting where objects have $2$-dimensional structure,
so too should the homs between them.
Indeed, we wouldn't be able to represent disks as global elements if this were not the case.
This motivates the definition of the following $3$-dimensional morphisms
\cite{kelly-1974-2_categories, johnson-2020-2_dimensional_categories}.

\begin{definition}[modification between oplax transformations]
	A \define[modification between transformations]{modification}
	between oplax transformations between $2$-functors between $2$-categories
	$μ : \hom[\hom[\hom[]{\cat{C}}{\cat{D}}]{F}{G}]{α}{β}$
	consists of the following data:
	\begin{description}
		\item[object-component disks:]
			for each object of the domain $2$-category $A : \cat{C}$
			a disk of the codomain $2$-category $μ A : \hom[\hom[\cat{D}]{F A}{G A}]{α A}{β A}$.
	\end{description}
	These object-component disks are required to satisfy the following relation:
	\begin{description}
		\item[naturality for arrows:]
			for each $\cat{C}$-arrow $f : \hom{A}{B}$ we have
			\begin{equation} \label{modification naturality for arrows}
				\comp{(\comp[2]{\comp{} (F f) , μ B}) , β f} = \comp{α f , (\comp[2]{μ A , \comp{} (G f)})}
			\end{equation}
	\end{description}
\end{definition}

Modifications are $3$-cells between parallel $2$-cells,
so we represent them in surface diagrams as points separating co-bounded lines.
Their naturality for arrows \eqref{modification naturality for arrows}
is interchanger naturality in the second index \eqref{interchanger naturality}.
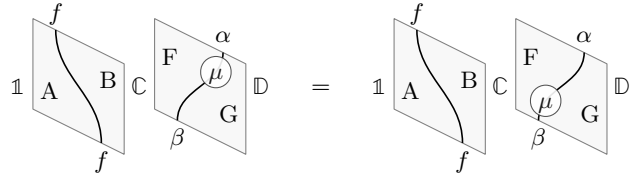
$$
	\begin{tikzpicture}[string diagram , x = {(12mm , -6mm)} , y = {(0mm , -12mm)} , z = {(16mm , 0mm)} , baseline=(align.base)]
		\node at (1/2 , 1/2 , 3/2) {\cat{D}} ;
		\coordinate (sheet) at (1/2 , 1/2 , 1) ;
		\draw [sheet]
			($ (sheet) + (-1/2 , -1/2 , 0) $) to coordinate [pos = 3/4] (in)
			($ (sheet) + (1/2 , -1/2 , 0) $) to coordinate [pos = 3/4] (to)
			($ (sheet) + (1/2 , 1/2 , 0) $) to coordinate [pos = 3/4] (out)
			($ (sheet) + (-1/2 , 1/2 , 0) $) to coordinate [pos = 3/4] (from)
			cycle
		;
		\draw [on sheet] (in) to [out = south , in = north] node [fromabove] {α} coordinate [pos = 1/4] (cell) node [tobelow] {β} (out) ;
		\node [bead] at (cell) {μ} ;
		\node [anchor = north west] at (from) {F} ;
		\node [anchor = south east] at (to) {G} ;
		\node at (1/2 , 1/2 , 1/2) {\cat{C}} ;
		\coordinate (sheet) at (1/2 , 1/2 , 0) ;
		\draw [sheet]
			($ (sheet) + (-1/2 , -1/2 , 0) $) to coordinate [pos = 1/4] (in)
			($ (sheet) + (1/2 , -1/2 , 0) $) to coordinate [pos = 0] (to)
			($ (sheet) + (1/2 , 1/2 , 0) $) to coordinate [pos = 1/4] (out)
			($ (sheet) + (-1/2 , 1/2 , 0) $) to coordinate [pos = 0] (from)
			cycle
		;
		\draw [on sheet] (in) to [out = south , in = north] node [fromabove] {f} coordinate [pos = 1/2] (cell) node [tobelow] {f} (out) ;
		\node at (cell) {} ;
		\node [anchor = south west] at (from) {A} ;
		\node [anchor = north east] at (to) {B} ;
		\node (align) at (1/2 , 1/2 , -1/2) {\cat{1}} ;
	\end{tikzpicture}
	\quad = \quad
	\begin{tikzpicture}[string diagram , x = {(12mm , -6mm)} , y = {(0mm , -12mm)} , z = {(16mm , 0mm)} , baseline=(align.base)]
		\node at (1/2 , 1/2 , 3/2) {\cat{D}} ;
		\coordinate (sheet) at (1/2 , 1/2 , 1) ;
		\draw [sheet]
			($ (sheet) + (-1/2 , -1/2 , 0) $) to coordinate [pos = 3/4] (in)
			($ (sheet) + (1/2 , -1/2 , 0) $) to coordinate [pos = 3/4] (to)
			($ (sheet) + (1/2 , 1/2 , 0) $) to coordinate [pos = 3/4] (out)
			($ (sheet) + (-1/2 , 1/2 , 0) $) to coordinate [pos = 3/4] (from)
			cycle
		;
		\draw [on sheet] (in) to [out = south , in = north] node [fromabove] {α} coordinate [pos = 3/4] (cell) node [tobelow] {β} (out) ;
		\node [bead] at (cell) {μ} ;
		\node [anchor = north west] at (from) {F} ;
		\node [anchor = south east] at (to) {G} ;
		\node at (1/2 , 1/2 , 1/2) {\cat{C}} ;
		\coordinate (sheet) at (1/2 , 1/2 , 0) ;
		\draw [sheet]
			($ (sheet) + (-1/2 , -1/2 , 0) $) to coordinate [pos = 1/4] (in)
			($ (sheet) + (1/2 , -1/2 , 0) $) to coordinate [pos = 0] (to)
			($ (sheet) + (1/2 , 1/2 , 0) $) to coordinate [pos = 1/4] (out)
			($ (sheet) + (-1/2 , 1/2 , 0) $) to coordinate [pos = 0] (from)
			cycle
		;
		\draw [on sheet] (in) to [out = south , in = north] node [fromabove] {f} coordinate [pos = 1/2] (cell) node [tobelow] {f} (out) ;
		\node at (cell) {} ;
		\node [anchor = south west] at (from) {A} ;
		\node [anchor = north east] at (to) {B} ;
		\node (align) at (1/2 , 1/2 , -1/2) {\cat{1}} ;
	\end{tikzpicture}
$$

The composition structure of transformations and modifications
is that of $2$-cells and $3$-cells in a Gray-category and,
as we have seen in section \ref{section: gray categories},
is apparent from their diagrammatic representations.

We see that in the construction of
the Gray-category of $2$-categories, $2$-functors, (lax/oplax/pseudo) transformations and modification
by explicit definitions,
the structure of interchangers is given as part of the structure of transformations, namely, as their arrow-component disks,
and the naturality of interchangers in their first index as a property of this structure;
whereas the naturality of interchangers in their second index is given as a property of modifications.
In contrast, in a general \refer{Gray-category} interchangers are structures in their own right,
and their naturality in each index is a property of this structure itself.
Arguably, the nature of the Gray-category of $2$-categories appears more clearly when viewed abstractly,
where the interchanger $3$-cells have an independent existence and properties uniform in their indices.

\end{sect}

\begin{sect}{Cones and Limits of $2$-Functors} \label{section: cones and limits}
	\subimport*{parts/}{cones}
\end{sect}

\begin{sect}{Lax Functoriality} \label{section: lax functoriality}
	In this section we consider the setting where
we require the morphisms between $2$-categories to preserve their composition
only up to a specified directed comparison structure
\cite{benabou-1967-bicategories, leinster-1998-bicategories}.
Such morphisms are called lax or oplax functors, depending on the variance of their comparison structure.

\begin{definition}[lax functor between $2$-categories]
	A \define{lax functor} between $2$-categories $F : \hom{\cat{C}}{\cat{D}}$
	has an object map and hom functors, just like a \refer[2-functor]{$2$-functor}.
	But instead of strictly preserving the composition of arrows and the horizontal composition of disks,
	a lax functor has a \define{lax comparison structure} $θ$,
	which for each path of consecutive arrows $f_1 , ⋯ , f_n$ in $\cat{C}$
	gives a comparison disk $θ \tuple*{f_1 , ⋯ , f_n} : \hom{\comp{F f_1 , … , F f_n}}{F (\comp{f_1 , … , f_n})}$ in $\cat{D}$.
	These comparison disks are required to be natural and coherent in the sense described next.
\end{definition}

A presentation for a comparison structure for a lax functor can be given as follows.
\begin{description}
	\item[nullary comparitors:]
		for each object $A : \cat{C}$ a natural transformation
		$
			θ_0^{A} : \hom{\comp{}_{\cat{D}} (\celldim{0}{F} A)}{\celldim{1}{F} (\comp{}_{\cat{C}} A)}
		$,
		$$
			\begin{tikzpicture}[math diagram , x = {(24mm , 0mm)} , y = {(0mm , -12mm)} , baseline=(align.base)]
				\node (00) at (0 , 0) {\cat{1}} ;
				\node (10) at (1 , 0) {\cat{1}} ;
				\node (01) at (0 , 1) {\hom[\cat{C}]{A}{A}} ;
				\node (11) at (1 , 1) {\hom[\cat{D}]{F A}{F A}} ;
				\draw [equals] (00) to node [auto] {} (10) ;
				\draw (01) to node [swap] {F} (11) ;
				\draw (00) to node [swap] {\comp{}_{\cat{C}} A} (01) ;
				\draw (10) to node [auto] (align) {\comp{}_{\cat{D}} (F A)} (11) ;
				\path (10) to node [pos = 2/5 , anchor = center] {θ_0^{A}} node [pos = 3/5 , anchor = center , sloped] {⇐} (01) ;
			\end{tikzpicture}
		$$
	\item[binary comparitors:]
		for objects $A , B , C : \cat{C}$ a natural transformation
		$θ_2^{A , B , C} : \hom{(\celldim{1}{F} \arg) \comp{,} _{\cat{D}} (\celldim{1}{F} \arg)}{\celldim{1}{F} (\arg \comp{,}_{\cat{C}} \arg)}$.
		$$
			\begin{tikzpicture}[math diagram , x = {(48mm , 0mm)} , y = {(0mm , -12mm)} , baseline=(align.base)]
				\node (00) at (0 , 0) {\product{\hom[\cat{C}]{A}{B} , \hom[\cat{C}]{B}{C}}} ;
				\node (10) at (1 , 0) {\product{\hom[\cat{D}]{F A}{F B} , \hom[\cat{D}]{F B}{F C}}} ;
				\node (01) at (0 , 1) {\hom[\cat{C}]{A}{C}} ;
				\node (11) at (1 , 1) {\hom[\cat{D}]{F A}{F C}} ;
				\draw (00) to node [auto] {\product{F , F}} (10) ;
				\draw (01) to node [swap] {F} (11) ;
				\draw (00) to node [swap] {\arg ⋅ _{\cat{C}} \arg} (01) ;
				\draw (10) to node [auto] (align) {\arg ⋅ _{\cat{D}} \arg} (11) ;
				\path (10) to node [pos = 2/5 , anchor = center] {θ_2^{A,B,C}} node [pos = 3/5 , anchor = center , sloped] {⇐} (01) ;
			\end{tikzpicture}
		$$
\end{description}
These, respectively,  have the following component disks in $\cat{D}$:
\[
	θ \tuple*{A} ≔ θ_0^{A} \tuple*{⋆} : \hom{\comp{} (F A)}{F (\comp{} \, A)}
	\; \text{and} \;
	θ \tuple*{f , g} ≔ θ_2^{A,B,C} \tuple*{f , g} : \hom{\comp{F f , F g}}{F (\comp{f , g})}
\]
\begin{equation} \label{lax comparitor components}
	\begin{tikzpicture}[string diagram , x = {(12mm , 0mm)} , y = {(0mm , -12mm)} , baseline=(align.base)]
		\coordinate (start) at (1/2 , 1/2) ;
		\draw (start) to node [tobelow] {F (\comp{} \, A)} (start |- 1 , 1) ;
		\node [bead] (align) at (start) {θ_0} ;
		\node [anchor = south] at ($ (start) + (0 , -1/4) $) {F A} ;
	\end{tikzpicture}
	\qquad \text{and} \qquad
	\begin{tikzpicture}[string diagram , x = {(12mm , 0mm)} , y = {(0mm , -12mm)} , baseline=(align.base)]
		\coordinate (merge) at (1/2 , 1/2) ;
		\draw ($ (merge |- 0 , 0) + (-3/8 , 0) $) [out = south , in = west] to node [fromabove] {F f} (merge) ;
		\draw ($ (merge |- 0 , 0) + (3/8 , 0) $) [out = south , in = east] to node [fromabove] {F g} (merge) ;
		\draw (merge) to node [tobelow] {F (\comp{f , g})} (merge |- 1 , 1) ;
		\node [bead] (align) at (merge) {θ_2} ;
		\node [anchor = north east] at ($ (merge) + (-1/4 , 0) $) {F A} ;
		\node [anchor = south] at ($ (merge) + (0 , -1/4) $) {F B} ;
		\node [anchor = north west] at ($ (merge) + (1/4 , 0) $) {F C} ;
	\end{tikzpicture}
\end{equation}
These are required to satisfy the following coherence laws.
\begin{description}
	\item[comparitor unit laws:]
		for $\cat{C}$-arrow $f : \hom{A}{B}$,
		the relations
		\[
			\comp{(\comp[2]{θ \tuple*{A} , F f}) , θ \tuple*{\comp{} \, A , f}}
			\quad = \quad
			\idty{F f}
			\quad = \quad
			\comp{(\comp[2]{F f , θ \tuple*{B}}) , θ \tuple*{f , \comp{} \, B}}
		\]
		\begin{equation} \label{comparitor unit law}
			\begin{tikzpicture}[string diagram , x = {(12mm , 0mm)} , y = {(0mm , -12mm)} , baseline=(current bounding box.center)]
				\coordinate (merge) at (1/2 , 1/2) ;
				\begin{scope}
					\clip [name path = clip] (0 , 1/16) rectangle (1 , 1) ;
					\draw [name path = p₀] ($ (merge |- 0 , 0) + (-1/3 , 0) $) to [out = south , in = west] (merge) ;
					\path [name intersections = {of = clip and p₀ , by = {start}}] ;
				\end{scope}
				\draw [name path = p₁] ($ (merge |- 0 , 0) + (1/3 , 0) $) to [out = south , in = east] node [fromabove] {F f} (merge) ;
				\draw (merge) to node [tobelow] {F f} (merge |- 1 , 1) ;
				\node [bead] at (start) {θ_0} ;
				\node [bead] at (merge) {θ_2} ;
			\end{tikzpicture}
			\qquad = \qquad
			\begin{tikzpicture}[string diagram , x = {(12mm , 0mm)} , y = {(0mm , -12mm)} , baseline=(current bounding box.center)]
				\coordinate (merge) at (1/2 , 1/2) ;
				\draw (merge |- 0 , 0) to [out = south , in = north] node [fromabove] {F f} node [tobelow] {F f} (merge |- 1 , 1) ;
			\end{tikzpicture}
			\qquad = \qquad
			\begin{tikzpicture}[string diagram , x = {(12mm , 0mm)} , y = {(0mm , -12mm)} , baseline=(current bounding box.center)]
				\coordinate (merge) at (1/2 , 1/2) ;
				\draw [name path = p₀] ($ (merge |- 0 , 0) + (-1/3 , 0) $) to [out = south , in = west] node [fromabove] {F f} (merge) ;
				\begin{scope}
					\clip [name path = clip] (0 , 1/16) rectangle (1 , 1) ;
					\draw [name path = p₁] ($ (merge |- 0 , 0) + (1/3 , 0) $) to [out = south , in = east] (merge) ;
					\path [name intersections = {of = clip and p₁ , by = {start}}] ;
				\end{scope}
				\draw (merge) to node [tobelow] {F f} (merge |- 1 , 1) ;
				\node [bead] at (start) {θ_0} ;
				\node [bead] at (merge) {θ_2} ;
			\end{tikzpicture}
		\end{equation}
		
	\item[comparitor associative law:]
		for consecutive $\cat{C}$-arrows $f : \hom{A}{B}$, $g : \hom{B}{C}$, and  $h : \hom{C}{D}$,
		the relation
		\[
			\comp{(\comp[2]{θ \tuple*{f , g} , F h}) , θ \tuple*{\comp{f , g} , h}}
			\quad = \quad
			\comp{(\comp[2]{F f , θ \tuple*{g , h}}) , θ \tuple*{f , \comp{g , h}}}
		\]
		\begin{equation} \label{comparitor associative law}
			\begin{tikzpicture}[string diagram , x = {(20mm , 0mm)} , y = {(0mm , -16mm)} , baseline=(current bounding box.center)]
				\coordinate (merge) at (2/3 , 2/3) ;
				\begin{scope}
					\clip [name path = clip] (0 , 1/3) rectangle (1 , 1) ;
					\draw [name path = p₀] ($ (merge |- 0 , 0) + (-1/3 , 0) $) to [out = south , in = west] node [fromabove] {F f} (merge) ;
					\path [name intersections = {of = clip and p₀ , by = {merge'}}] ;
				\end{scope}
				\draw ($ (merge' |- 0 , 0) + (-1/5 , 0) $) to [out = south , in = west] node [fromabove] {F f} (merge') ;
				\draw ($ (merge' |- 0 , 0) + (1/5 , 0) $) to [out = south , in = east] node [fromabove] {F g} (merge') ;
				\draw [name path = p₁] ($ (merge |- 0 , 0) + (1/3 , 0) $) to [out = south , in = east] node [fromabove] {F h} (merge) ;
				\draw (merge) to node [tobelow] {F (\comp{f , g , h})} (merge |- 1 , 1) ;
				\node [bead] at (merge') {θ_2} ;
				\node [bead] at (merge) {θ_2} ;
			\end{tikzpicture}
			\qquad = \qquad
			\begin{tikzpicture}[string diagram , x = {(20mm , 0mm)} , y = {(0mm , -16mm)} , baseline=(current bounding box.center)]
				\coordinate (merge) at (1/3 , 2/3) ;
				\draw [name path = p₀] ($ (merge |- 0 , 0) + (-1/3 , 0) $) to [out = south , in = west] node [fromabove] {F f} (merge) ;
				\begin{scope}
					\clip [name path = clip] (0 , 1/3) rectangle (1 , 1) ;
					\draw [name path = p₁] ($ (merge |- 0 , 0) + (1/3 , 0) $) to [out = south , in = east] node [fromabove] {F g} (merge) ;
					\path [name intersections = {of = clip and p₁ , by = {merge'}}] ;
				\end{scope}
				\draw ($ (merge' |- 0 , 0) + (-1/5 , 0) $) to [out = south , in = west] node [fromabove] (Fg) {F g} (merge') ;
				\draw ($ (merge' |- 0 , 0) + (1/5 , 0) $) to [out = south , in = east] node [fromabove] {F h} (merge') ;
				\draw (merge) to node [tobelow] {F (\comp{f , g , h})} (merge |- 1 , 1) ;
				\node [bead] at (merge') {θ_2} ;
				\node [bead] at (merge) {θ_2} ;
			\end{tikzpicture}
		\end{equation}
\end{description}

The naturality of the lax comparitors amount to the following relations on their components.
\begin{description}	
	\item[comparitor unit naturality:]
		for object $A : \cat{C}$,
		we have
		$
			\comp{θ \tuple*{A} , F (\comp[2]{} A)}
			=
			θ \tuple*{A}
		$:
		\begin{equation} \label{comparitor unit naturality}
			\begin{tikzpicture}[string diagram , x = {(16mm , 0mm)} , y = {(0mm , -16mm)} , baseline=(align.base)]
				\coordinate (start) at (1/2 , 1/3) ;
				\coordinate (id) at (1/2 , 3/4) ;
				\draw [name path = F id] (start) to node [tobelow] {F (\comp{} \, A)} (start |- 1 , 1) ;
				\path [name path = cell height] (0 , 0 |- id) to (1 , 1 |- id) ;
				\node [bead] at (start) {θ_0} ;
				\path [name intersections={of = F id and cell height , by = {cell}}] ;
				\node [bead] (align) at (cell) {F (\comp[2]{} A)} ;
			\end{tikzpicture}
			\quad = \quad
			\begin{tikzpicture}[string diagram , x = {(16mm , 0mm)} , y = {(0mm , -16mm)} , baseline=(align.base)]
				\coordinate (id) at (1/2 , 1/3) ;
				\coordinate (start) at (1/2 , 3/4) ;
				\draw [name path = F id] (start) to node [tobelow] {F (\comp{} \, A)} (start |- 1 , 1) ;
				\node [bead] (align) at (start) {θ_0} ;
			\end{tikzpicture}
		\end{equation}
		But because $F (\comp[2]{} A) ≔ F (\comp{} (\comp{} \, A)) = \comp{} (F (\comp{} \, A))$ by the local functoriality of $F$,
		this condition is redundant.
		
	\item[comparitor composition naturality:]
		for adjacent disks $φ : \hom[\hom[\cat{C}]{A}{B}]{f}{f′}$ and  $ψ : \hom[\hom[\cat{C}]{B}{C}]{g}{g′}$,
		we have
		$
			\comp{θ \tuple*{f , g} , F (\comp[2]{φ , ψ})}
			=
			\comp{(\comp[2]{F φ , F ψ}) , θ \tuple*{f′ , g′}}
		$:
		\begin{equation} \label{comparitor composition naturality}
			\begin{tikzpicture}[string diagram , x = {(16mm , 0mm)} , y = {(0mm , -16mm)} , baseline=(align.base)]
				\coordinate (merge) at (1/2 , 1/3) ;
				\coordinate (disk) at (1/2 , 3/4) ;
				\draw [name path = f] ($ (merge |- 0 , 0) + (-1/3 , 0) $) to [out = south , in = west] node [fromabove] {F f} (merge) ;
				\draw [name path = g] ($ (merge |- 0 , 0) + (1/3 , 0) $) to [out = south , in = east] node [fromabove] {F g} (merge) ;
				\draw [name path = fg] (merge) to node [tobelow] {F (\comp{f′ , g′})} (merge |- 1 , 1) ;
				\path [name path = disk height] (0 , 0 |- disk) to (1 , 1 |- disk) ;
				\node [bead] at (merge) {θ_2} ;
				\path [name intersections = {of = disk height and fg , by = {disk}}] ;
				\node [bead] (align) at (disk) {F (\comp[2]{φ , ψ})} ;
			\end{tikzpicture}
			\quad = \quad
			\begin{tikzpicture}[string diagram , x = {(16mm , 0mm)} , y = {(0mm , -16mm)} , baseline=(align.base)]
				\coordinate (merge) at (1/2 , 3/4) ;
				\coordinate (disk) at (1/2 , 1/3) ;
				\draw [name path = f] ($ (merge |- 0 , 0) + (-1/3 , 0) $) to [out = south , in = west] node [fromabove] {F f} (merge) ;
				\draw [name path = g] ($ (merge |- 0 , 0) + (1/3 , 0) $) to [out = south , in = east] node [fromabove] {F g} (merge) ;
				\draw [name path = fg] (merge) to node [tobelow] {F (\comp{f′ , g′})} (merge |- 1 , 1) ;
				\path [name path = disk height] (0 , 0 |- disk) to (1 , 1 |- disk) ;
				\node [bead] (align) at (merge) {θ_2} ;
				\path [name intersections = {of = disk height and f , by = {disk1}}] ;
				\node [bead] at (disk1) {F φ} ;
				\path [name intersections = {of = disk height and g , by = {disk2}}] ;
				\node [bead] at (disk2) {F ψ} ;
			\end{tikzpicture}
		\end{equation}
\end{description}

An \define{oplax functor} is one where the comparitors are oriented the other way and satisfy dual conditions.
A \define{pseudo functor} is both lax and oplax with invertible comparitors,
and a \refer[2-functor]{$2$-functor} is a pseudo functor with identity comparitors.

We would like to incorporate lax functors into our graphical calculus,
but in trying to do so we encounter a problem.
For a $2$-functor $F : \hom{\cat{C}}{\cat{D}}$ it doesn't matter
whether we compose a diagram of cells in $\cat{C}$ and then send the composite cell to $\cat{D}$ by $F$
or whether we send the cells over individually by $F$ and then compose them in $\cat{D}$, the result is the same.
But if $F$ is only a lax functor then \emph{where} we do the composition makes a difference:
the composite in $\cat{D}$ of $F$-images is generally not the same as the $F$-image of the composite in $\cat{C}$
because a lax functor preserves the composition structure of its domain $2$-category only up to its comparitor's components.

In fact, if we apply a lax functor to the cells of a diagram that is composable in its domain $2$-category,
the resulting collection of cells may not even be composable in its codomain $2$-category.
For example, for the diagram shown below,
the vertical codomain of the $F$-image composite of its top half $\comp[2]{F φ , F h}$
and the vertical domain of the $F$-image composite of its bottom half $\comp[2]{F f , F ψ}$
coincide only if $\comp{F (\comp{f , g}) , F h}  =  \comp{F f , F (\comp{g , h})}$.
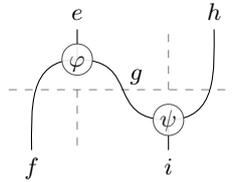
\[
	\begin{tikzpicture}[string diagram , x = {(24mm , 0mm)} , y = {(0mm , -16mm)} , baseline=(align.base)]
		\coordinate (split) at (1/4 , 1/4) ;
		\coordinate (merge) at (3/4 , 3/4) ;
		\coordinate (in 1) at (split |- 0 , 0) ;
		\coordinate (in 2) at (1 , 0) ;
		\coordinate (out 1) at (0 , 1) ;
		\coordinate (out 2) at (merge |- 1 , 1) ;
		\draw
			(out 1) to [out = north , in looseness = 1.25 , in = west] node [frombelow] {f}
			(split) to [out = east , looseness = 1.5 , in = west] node [auto] {g}
			(merge) to [out = east , out looseness = 1.25 , in = south] node [toabove] {h}
			(in 2)
		;
		\draw (in 1) to [out = south , in = north] node [fromabove] {e} (split) ;
		\draw (merge) to [out = south , in = north] node [tobelow] {i} (out 2) ;
		\node [bead] at (split) {φ} ;
		\node [bead] at (merge) {ψ} ;
		\draw [dashed , opacity = 0.5] (-1/8 , 1/2) to coordinate [pos = 1/2] (mid) (9/8 , 1/2) ;
		\draw [dashed , opacity = 0.5] (split |- mid) to (split |- out 1) ;
		\draw [dashed , opacity = 0.5] (merge |- mid) to (merge |- in 2) ;
	\end{tikzpicture}
\]

The lax functor images of the cells of a composable diagram themselves form a composable diagram
if the diagram can be decomposed into a horizontal composition of vertical compositions of cells.
We will call such diagrams \define[stratified diagram]{stratified}.

In order to accommodate lax functors
we augment our surface diagrams with annotations to indicate where compositions take place.
A \define{compositor} sheet is a surface drawn within a volume representing a $2$-category.
It is intended to represent the composition of the portion of the diagram preceding it into a single cell.
A compositor sheet is significant only when it occurs before or after
a surface representing a lax or oplax functor that is not a strict $2$-functor.
If the compositor sheet precedes the functor surface in the composition order
then the diagram is first composed in the functor's domain
and then the functor is applied to the resulting composite cell.
If the compositor sheet follows the functor surface in the composition order
then the functor is first applied to the constituent cells of a diagram in the functor's domain
and then the resulting diagram is composed in the functor's codomain.
The latter configuration requires the side-condition that the diagram's cells' functor images be composable,
and is satisfied whenever the diagram is \refer[stratified diagram]{stratified}.
For example, the left diagram below represents $F (\comp[2]{φ , ψ}) : \hom[\cat{D}]{F (\comp{f , g})}{F (\comp{f′ , g′})}$,
whereas the right diagram represents $\comp[2]{F φ , F ψ} : \hom[\cat{D}]{\comp{F f , F g}}{\comp{F f′ , F g′}}$.
\[
	\begin{tikzpicture}[string diagram , x = {(16mm , -8mm)} , y = {(0mm , -16mm)} , z = {(24mm , 0mm)} , baseline=(align.base)]
		\node at (1/2 , 1/2 , 3/2) {\cat{D}} ;
		\coordinate (sheet) at (1/2 , 1/2 , 1) ;
		\draw [sheet]
			($ (sheet) + (-1/2 , -1/2 , 0) $) to
			($ (sheet) + (1/2 , -1/2 , 0) $) to
			($ (sheet) + (1/2 , 1/2 , 0) $) to
			($ (sheet) + (-1/2 , 1/2 , 0) $) to
			cycle
		;
		\node at ($ (sheet) + (1/3 , 0 , 0) $) {F} ;
		\coordinate (sheet) at (1/2 , 1/2 , 3/4) ;
		\draw [compositor]
			($ (sheet) + (-1/2 , -1/2 , 0) $) to
			($ (sheet) + (1/2 , -1/2 , 0) $) to
			($ (sheet) + (1/2 , 1/2 , 0) $) to
			($ (sheet) + (-1/2 , 1/2 , 0) $) to
			cycle
		;
		\node at (1/2 , 1/2 , 1/2) {\cat{C}} ;
		\coordinate (sheet) at (1/2 , 1/2 , 0) ;
		\draw [sheet]
			($ (sheet) + (-1/2 , -1/2 , 0) $) to coordinate [pos = 1/6] (in 1) coordinate [pos = 3/6] (in 2)
			($ (sheet) + (1/2 , -1/2 , 0) $) to coordinate [pos = 0] (to)
			($ (sheet) + (1/2 , 1/2 , 0) $) to coordinate [pos = 1/6] (out 2) coordinate [pos = 3/6] (out 1)
			($ (sheet) + (-1/2 , 1/2 , 0) $) to coordinate [pos = 0] (from)
			cycle
		;
		\draw [on sheet , name path = path 1]
			(in 1) to [out = south , out looseness = 1.25 , in looseness = 0.75 , in = north]
				node [fromabove] {f} node [tobelow] {f′}
			(out 1)
		;
		\draw [on sheet , name path = path 2]
			(in 2) to [out = south , out looseness = 0.75 , in looseness = 1.25 , in = north]
				node [fromabove] {g} node [tobelow] {g′}
			(out 2)
		;
		\path [name path = cell height] ($ (sheet) + (-1/2 , 0 , 0) $) to ($ (sheet) + (1/2 , 0 , 0) $) ;
		\path [name intersections = {of = cell height and path 1 , by = {cell 1}}] ;
		\path [name intersections = {of = cell height and path 2 , by = {cell 2}}] ;
		\node [bead] at (cell 1) {φ} ;
		\node [bead] at (cell 2) {ψ} ;
		\node [anchor = south west] at (from) {A} ;
		\node [anchor = center] at ($ (from) ! 1/2 ! (to) $) {B} ;
		\node [anchor = north east] at (to) {C} ;
		\node (align) at (1/2 , 1/2 , -1/2) {\cat{1}} ;
	\end{tikzpicture}
	\quad , \quad
	\begin{tikzpicture}[string diagram , x = {(16mm , -8mm)} , y = {(0mm , -16mm)} , z = {(24mm , 0mm)} , baseline=(align.base)]
		\node at (1/2 , 1/2 , 3/2) {\cat{D}} ;
		\coordinate (sheet) at (1/2 , 1/2 , 5/4) ;
		\draw [compositor]
			($ (sheet) + (-1/2 , -1/2 , 0) $) to
			($ (sheet) + (1/2 , -1/2 , 0) $) to
			($ (sheet) + (1/2 , 1/2 , 0) $) to
			($ (sheet) + (-1/2 , 1/2 , 0) $) to
			cycle
		;
		\coordinate (sheet) at (1/2 , 1/2 , 1) ;
		\draw [sheet]
			($ (sheet) + (-1/2 , -1/2 , 0) $) to
			($ (sheet) + (1/2 , -1/2 , 0) $) to
			($ (sheet) + (1/2 , 1/2 , 0) $) to
			($ (sheet) + (-1/2 , 1/2 , 0) $) to
			cycle
		;
		\node at ($ (sheet) + (-1/3 , 0 , 0) $) {F} ;
		\node at (1/2 , 1/2 , 1/2) {\cat{C}} ;
		\coordinate (sheet) at (1/2 , 1/2 , 0) ;
		\draw [sheet]
			($ (sheet) + (-1/2 , -1/2 , 0) $) to coordinate [pos = 1/6] (in 1) coordinate [pos = 3/6] (in 2)
			($ (sheet) + (1/2 , -1/2 , 0) $) to coordinate [pos = 0] (to)
			($ (sheet) + (1/2 , 1/2 , 0) $) to coordinate [pos = 1/6] (out 2) coordinate [pos = 3/6] (out 1)
			($ (sheet) + (-1/2 , 1/2 , 0) $) to coordinate [pos = 0] (from)
			cycle
		;
		\draw [on sheet , name path = path 1]
			(in 1) to [out = south , out looseness = 1.25 , in looseness = 0.75 , in = north]
				node [fromabove] {f} node [tobelow] {f′}
			(out 1)
		;
		\draw [on sheet , name path = path 2]
			(in 2) to [out = south , out looseness = 0.75 , in looseness = 1.25 , in = north]
				node [fromabove] {g} node [tobelow] {g′}
			(out 2)
		;
		\path [name path = cell height] ($ (sheet) + (-1/2 , 0 , 0) $) to ($ (sheet) + (1/2 , 0 , 0) $) ;
		\path [name intersections = {of = cell height and path 1 , by = {cell 1}}] ;
		\path [name intersections = {of = cell height and path 2 , by = {cell 2}}] ;
		\node [bead] at (cell 1) {φ} ;
		\node [bead] at (cell 2) {ψ} ;
		\node [anchor = south west] at (from) {A} ;
		\node [anchor = center] at ($ (from) ! 1/2 ! (to) $) {B} ;
		\node [anchor = north east] at (to) {C} ;
		\node (align) at (1/2 , 1/2 , -1/2) {\cat{1}} ;
	\end{tikzpicture}
\]
The global elements $A , B , C : \hom{\cat{\product{}}}{\cat{C}}$
should themselves be strict \refer[2-functor]{2-functors},
otherwise their lax comparison structures would induce nontrivial monads on the corresponding $\cat{C}$-objects,
as was observed by Bénabou \cite{benabou-1967-bicategories}.

The comparison structure $θ$ for a lax functor $F : \hom{\cat{C}}{\cat{D}}$ provides a map
from the composite in $\cat{D}$ of $F$-images to the $F$-image of the composite in $\cat{C}$.
We can represent this in surface diagrams
as the passing of a compositor sheet backward in the composition order through the lax functor surface.
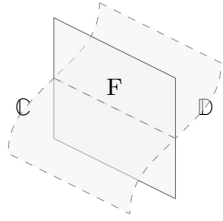
\[

\]
If we compose this with a diagram in $\cat{C}$ consisting of a consecutive pair arrows
then the projection gives the corresponding binary comparitor component.
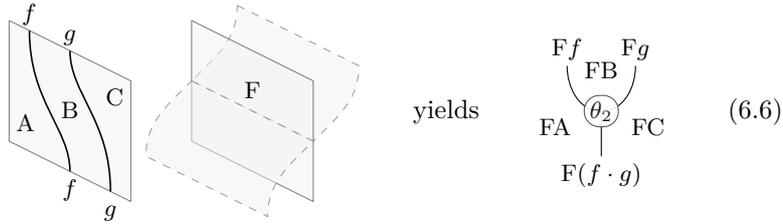
\begin{equation} \label{binary lax comparitor component}
	%
\end{equation}
If we compose it with a diagram consisting of no arrows
then the projection gives the corresponding nullary comparitor component.
\begin{equation} \label{nullary lax comparitor component}
	%
\end{equation}

The comparitor unit \eqref{comparitor unit law} and associative \eqref{comparitor associative law} laws
together imply that there is a unique $θ_n^{A_0,⋯,A_n}$ for each $n ∈ \Nat$,
with $θ \tuple*{f} ≔ θ_1^{A,B} \tuple*{f} : \hom{\comp{F f}}{F \comp{f}} = \idty{F f}$.

For a diagram consisting of a single disk $φ : \hom[\hom[\cat{C}]{A}{B}]{f}{f′}$,
the image of the diagram composite and the composite of the diagram comprising the disk's image
are both just the image of the disk.
So we have $\comp{θ \tuple*{f} , F φ} = F φ = \comp{F φ , θ \tuple*{f′}}$.
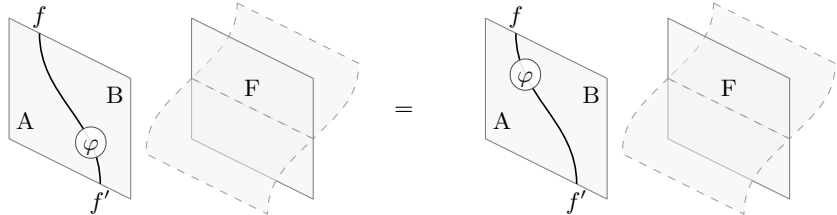
\[
	%
\]
Thus vertical composites of disks can freely slide past unary comparitors.
The same is true for horizontal composites of disks and binary comparitors
due to the comparitor naturality condition \eqref{comparitor composition naturality}.
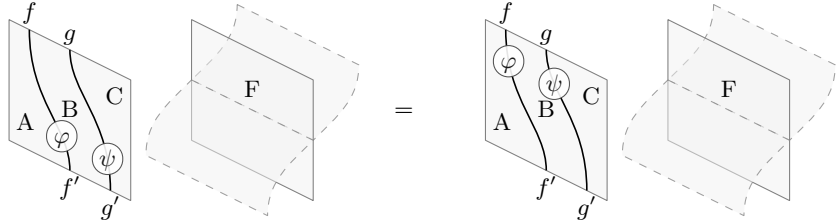
\[
	%
\]
Moreover, because identity disks are strict units of vertical composition
and their lax functor images are again identity disks,
we may slide any disk of a horizontal composition past a comparitor independently of the others.
So the next two diagrams are also equal to the previous two.
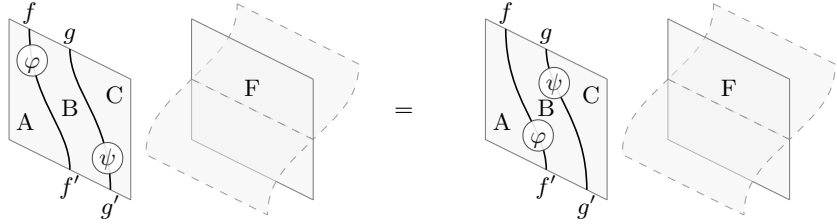
\[
	%
\]
That is to say, disks in different strata can move past a comparitor independently.
Thus any \refer{stratified diagram} in $\cat{C}$ may be positioned arbitrarily with respect to the comparitor for a lax functor $F$,
and all of the resulting projections into $\cat{D}$ will be equal.

Consecutive lax functors $F : \hom{\cat{C}}{\cat{D}}$ and $G : \hom{\cat{D}}{\cat{E}}$ compose strictly,
with composite comparitor components given by
\begin{equation} \label{lax functor composition}
	%
\end{equation}
In surface diagrams we represent these as
the passage of a compositor sheet through the two lax functor surfaces sequentially.
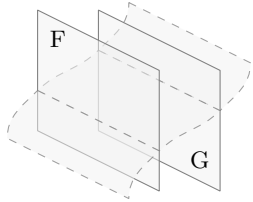
\[
	%
\]
This composition of lax functors is associative with identity $2$-functors acting as composition units.

We can define transformations between lax functors as well \cite{kelly-1974-2_categories, leinster-1998-bicategories}.

\begin{definition}[oplax transformation between lax functors]
	An \define[oplax transformation between lax functors]{oplax transformation} between lax functors between $2$-categories
	$α : \hom[\hom[]{\cat{C}}{\cat{D}}]{F}{G}$
	has object-component arrows and arrow-component disks
	that are natural for disks in the domain $2$-category \eqref{transformation naturality for disks}
	just like an \refer{oplax transformation between 2-functors}.
	
	However, now the arrow composition compatibility conditions \eqref{transformation arrow composition compatibility}
	do not make sense because for object $A : \cat{C}$ we have
	$α (\comp{} \, A) : \hom{\comp{F (\comp{} \, A) , α A}}{\comp{α A , G (\comp{} \, A)}}$ whereas
	$\comp{} (α A) : \hom{α A}{α A}$,
	and for consecutive arrows $f : \hom[\cat{C}]{A}{B}$ and $g : \hom[\cat{C}]{B}{C}$ we have
	$α (\comp{f , g}) : \hom{\comp{F (\comp{f , g}) , α C}}{\comp{α A , G (\comp{f , g})}}$ whereas
	$\comp{(\comp[2]{F f , α g}) , (\comp[2]{α f , G g})} : \hom{\comp{F f , F g , α C}}{\comp{α A , G f , G g}}$,
	but the lax functors $F$ and $G$ preserve arrow composition structure only up to their respective comparitors.
	
	So instead we require the arrow-component disks to be compatible with these comparitors
	in the sense that
	\begin{equation} \label{nullary comparitor compatibility}
		\comp{(\comp[2]{θ^F \tuple*{A} , α A}) , α (\comp{} \, A)}
		\; = \;
		\comp[2]{α A , θ^G \tuple*{A}}
	\end{equation}
	and
	\begin{equation} \label{binary comparitor compatibility}
		\comp{(\comp[2]{θ^F \tuple*{f , g} , α C}) , α (\comp{f , g})}
		\; = \;
		\comp{(\comp[2]{F f , α g}) , (\comp[2]{α f , G g}) , (\comp[2]{α A , θ^G \tuple*{f , g}})}
	\end{equation}
\end{definition}

We can represent these lax comparitor compatibility conditions in string diagrams as
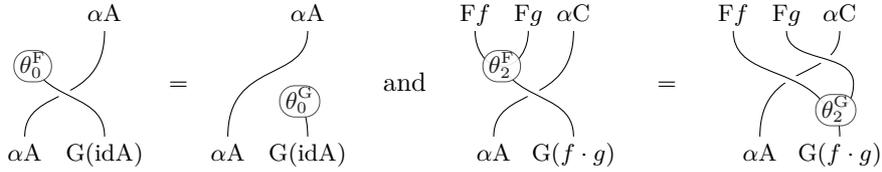
\[
		\begin{tikzpicture}[string diagram , x = {(14mm , 0mm)} , y = {(0mm , -14mm)} , baseline=(align.base)]
		\coordinate (in 1) at (1/8 , 0) ;
		\coordinate (in 2) at (7/8 , 0) ;
		\coordinate (out 1) at (1/8 , 1) ;
		\coordinate (out 2) at (7/8 , 1) ;
		\draw [name path = back , overcross]
			(in 2) to [out = south , out looseness = 1.5 , in = north]
				node [fromabove] {α A} node [tobelow] {α A}
			(out 1)
		;
		\begin{scope}
			\clip [name path = clip] (0 , 1/3) rectangle (4/3 , 4/3) ;
			\draw [name path = front , overcross]
				(in 1) to [out = south , out looseness = 1.5 , in = north]
					node [fromabove] {F (\comp{} A)} node [tobelow] {G (\comp{} A)}
				(out 2)
			;
			\path [name intersections={of = clip and front , by = {disk}}] ;
			\path [name intersections={of = front and back , by = {cell}}] ;
		\end{scope}
		\node [bead] at (disk) {θ_0^F} ;
		\node (align) at (1/2 , 1/2) {} ;
	\end{tikzpicture}
	\; = \;
		\begin{tikzpicture}[string diagram , x = {(14mm , 0mm)} , y = {(0mm , -14mm)} , baseline=(align.base)]
		\coordinate (in 1) at (1/8 , 0) ;
		\coordinate (in 2) at (7/8 , 0) ;
		\coordinate (out 1) at (1/8 , 1) ;
		\coordinate (out 2) at (7/8 , 1) ;
		\draw [name path = back , overcross]
			(in 2) to [out = south , in looseness = 1.5 , in = north]
				node [fromabove] {α A} node [tobelow] {α A}
			(out 1)
		;
		\begin{scope}
			\clip [name path = clip] (0 , 2/3) rectangle (4/3 , 4/3) ;
			\draw [name path = front , overcross]
				(in 1) to [out = south , in looseness = 1.5 , in = north]
					node [fromabove] {F (\comp{} A)} node [tobelow] {G (\comp{} A)}
				(out 2)
			;
			\path [name intersections={of = clip and front , by = {disk}}] ;
			\path [name intersections={of = front and back , by = {cell}}] ;
		\end{scope}
		\node [bead] at (disk) {θ_0^G} ;
		\node (align) at (1/2 , 1/2) {} ;
	\end{tikzpicture}
	\quad \text{and} \quad
		\begin{tikzpicture}[string diagram , x = {(14mm , 0mm)} , y = {(0mm , -14mm)} , baseline=(align.base)]
			\coordinate (in 1) at (1/8 , 0) ;
			\coordinate (in 2) at (7/8 , 0) ;
			\coordinate (out 1) at (1/8 , 1) ;
			\coordinate (out 2) at (7/8 , 1) ;
			\draw [name path = back , overcross]
				(in 2) to [out = south , out looseness = 1.5 , in = north]
					node [fromabove] {α C} node [tobelow] {α A}
				(out 1)
			;
			\begin{scope}
				\clip [name path = clip] (0 , 1/3) rectangle (4/3 , 4/3) ;
				\draw [name path = front , overcross]
					(in 1) to [out = south , out looseness = 1.5 , in = north]
						node [fromabove] {F (\comp{f , g})} node [tobelow] {G (\comp{f , g})}
					(out 2)
				;
				\path [name intersections={of = clip and front , by = {disk}}] ;
				\path [name intersections={of = front and back , by = {cell}}] ;
			\end{scope}
			\draw ($ (disk |- 0 , 0) + (-1/4 , 0) $) to [out = south , in = west] node [fromabove] {F f} (disk) ;
			\draw ($ (disk |- 0 , 0) + (1/4 , 0) $) to [out = south , in = east] node [fromabove] {F g} (disk) ;
			\node [bead] at (disk) {θ_2^F} ;
		\end{tikzpicture}
		\quad = \quad
		\begin{tikzpicture}[string diagram , x = {(14mm , 0mm)} , y = {(0mm , -14mm)} , baseline=(align.base)]
			\coordinate (in 1) at (1/8 , 0) ;
			\coordinate (in 2) at (7/8 , 0) ;
			\coordinate (out 1) at (1/8 , 1) ;
			\coordinate (out 2) at (7/8 , 1) ;
			\draw [name path = back , overcross]
				(in 2) to [out = south , in looseness = 1.5 , in = north]
					node [fromabove] {α C} node [tobelow] {α A}
				(out 1)
			;
			\begin{scope}
				\clip [name path = clip] (0 , 3/4) rectangle (4/3 , 4/3) ;
				\draw [name path = front , overcross]
					(in 1) to [out = south , in looseness = 1.5 , in = north]
						node [fromabove] {F (\comp{f , g})} node [tobelow] {G (\comp{f , g})}
					(out 2)
				;
				\path [name intersections={of = clip and front , by = {disk}}] ;
				\path [name intersections={of = front and back , by = {cell}}] ;
			\end{scope}
			\draw [name path = front 1 , overcross]
				($ (in 1) + (-1/4 , 0) $) to [out = south , in = north west] node [fromabove] {F f}
				($ (disk) + (-1/6 , -1/8) $) to [out = south east , in = west]
				(disk)
			;
			\draw [name path = front 2 , overcross]
				($ (in 1) + (1/4 , 0) $) to [out = south , in = north] node [fromabove] {F g}
				($ (disk) + (1/6 , -1/4) $) to [out = south , in = east]
				(disk)
			;
			\node [bead] at (disk) {θ_2^G} ;
		\end{tikzpicture}
\]
In surface diagrams this corresponds to the movement of the line representing the oplax transformation
past the intersection of a compositor sheet with a lax functor surface.
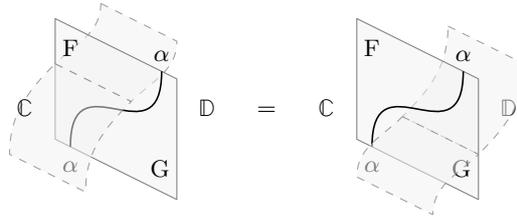
\[
	\begin{tikzpicture}[string diagram , x = {(16mm , -8mm)} , y = {(0mm , -16mm)} , z = {(24mm , 0mm)} , baseline=(align.base)]
		\node at (1/2 , 1/2 , 3/2) {\cat{D}} ;
		\coordinate (sheet) at (1/2 , 1/2 , 1) ;
		\draw [compositor]
			($ (sheet) + (-1/2 , -5/8 , 1/4) $) to
			($ (sheet) + (1/8 , -5/8 , 1/4) $) to [out = south , in = north east]
			($ (sheet) + (1/8 , -1/8 , 0) $) to
			($ (sheet) + (-1/2 , -1/8 , 0) $) to [out = north east , in = south]
			cycle
		;
		\draw [sheet]
			($ (sheet) + (-1/2 , -1/2 , 0) $) to coordinate [pos = 7/8] (in)
			($ (sheet) + (1/2 , -1/2 , 0) $) to coordinate [pos = 3/4] (to)
			($ (sheet) + (1/2 , 1/2 , 0) $) to coordinate [pos = 7/8] (out)
			($ (sheet) + (-1/2 , 1/2 , 0) $) to coordinate [pos = 3/4] (from)
			cycle
		;
		\draw [on sheet] (in) to [out = south , looseness = 2 , in = north] node [fromabove] {α} coordinate [pos = 1/2] (cell) node [tobelow] {α} (out) ;
		\node [anchor = west] at (from) {F} ;
		\node [anchor = east] at (to) {G} ;
		\draw [compositor]
			($ (sheet) + (-1/2 , -1/8 , 0) $) to
			($ (sheet) + (1/8 , -1/8 , 0) $) to [out = south west , in = north]
			($ (sheet) + (1/8 , 5/8 , -1/4) $) to
			($ (sheet) + (-1/2 , 5/8 , -1/4) $) to [out = north , in = south west]
			cycle
		;
		\node (align) at (1/2 , 1/2 , 1/2) {\cat{C}} ;
	\end{tikzpicture}
	\quad = \quad
	\begin{tikzpicture}[string diagram , x = {(16mm , -8mm)} , y = {(0mm , -16mm)} , z = {(24mm , 0mm)} , baseline=(align.base)]
		\node at (1/2 , 1/2 , 3/2) {\cat{D}} ;
		\coordinate (sheet) at (1/2 , 1/2 , 1) ;
		\draw [compositor]
			($ (sheet) + (-1/8 , -5/8 , 1/4) $) to
			($ (sheet) + (1/2 , -5/8 , 1/4) $) to [out = south , in = north east]
			($ (sheet) + (1/2 , 1/8 , 0) $) to
			($ (sheet) + (-1/8 , 1/8 , 0) $) to [out = north east , in = south]
			cycle
		;
		\draw [sheet]
			($ (sheet) + (-1/2 , -1/2 , 0) $) to coordinate [pos = 7/8] (in)
			($ (sheet) + (1/2 , -1/2 , 0) $) to coordinate [pos = 3/4] (to)
			($ (sheet) + (1/2 , 1/2 , 0) $) to coordinate [pos = 7/8] (out)
			($ (sheet) + (-1/2 , 1/2 , 0) $) to coordinate [pos = 3/4] (from)
			cycle
		;
		\draw [on sheet] (in) to [out = south , looseness = 2 , in = north] node [fromabove] {α} coordinate [pos = 1/2] (cell) node [tobelow] {α} (out) ;
		\node [anchor = west] at (from) {F} ;
		\node [anchor = east] at (to) {G} ;
		\draw [compositor]
			($ (sheet) + (-1/8 , 1/8 , 0) $) to
			($ (sheet) + (1/2 , 1/8 , 0) $) to [out = south west , in = north]
			($ (sheet) + (1/2 , 5/8 , -1/4) $) to
			($ (sheet) + (-1/8 , 5/8 , -1/4) $) to [out = north , in = south west]
			cycle
		;
		\node (align) at (1/2 , 1/2 , 1/2) {\cat{C}} ;
	\end{tikzpicture}
\]

A \define[modification between oplax transformations between lax functors]{modification}
between oplax transformations between lax functors
is defined just like a \refer[modification between oplax transformations]{modification}
between oplax transformations between $2$-functors,
with no additional compatibility relations involving the lax comparison structure
\cite{leinster-1998-bicategories, johnson-2020-2_dimensional_categories}.

Compositor sheets give us pretty pictures,
but we would like to understand them in terms of the diagram semantics of section \ref{section: diagram semantics}.
A compositor is intended to take a string diagram of disks of a $2$-category
and compose it to a single disk of that $2$-category.
For a $2$-category $\cat{C}$, a string diagram of $\cat{C}$-disks
is a $2$-generator of the $2$-computad
$\mathrm{forget} (\mathrm{diag} (\mathrm{forget} \, \cat{C}))$.
And in general, the operation that forms diagrams of cells of a $*$-category is the functor
$\comp{\mathrm{forget} , \mathrm{diag} , \mathrm{forget}} : \hom{*\Cat{Cat}}{*\Cat{Ctd}}$.
For brevity we write this functor as ``$\mathrm{fdf}$''.

Recall that the adjunction $\adjoint{\mathrm{diag} , \mathrm{forget}}$
has as counit $\mathrm{eval} : \hom{\comp{\mathrm{forget} , \mathrm{diag}}}{\comp{} (*\Cat{Cat})}$
whose components evaluate diagrams of cells to their composites,
and as unit $\mathrm{cone} : \hom{\comp{} (*\Cat{Ctd})}{\comp{\mathrm{diag} , \mathrm{forget}}}$
whose components create singleton diagrams.
From these we can form the following natural transformation:
$$
	κ
	\quad ≔ \quad
	\begin{tikzpicture}[string diagram , x = {(12mm , 0mm)} , y = {(0mm , -16mm)} , baseline=(align.base)]
		\coordinate (in1) at (0 , 0) ;
		\coordinate (in2) at (1 , 0) ;
		\coordinate (in3) at (2 , 0) ;
		\coordinate (out1) at (1 , 1) ;
		\coordinate (out2) at (2 , 1) ;
		\coordinate (out3) at (3 , 1) ;
		\coordinate (cup) at ({$ (in1) ! 1/2 ! (in2) $} |- {$ (in3) ! 1/3 ! (out1) $}) ;
		\coordinate (cap) at ({$ (out2) ! 1/2 ! (out3) $} |- {$ (in3) ! 2/3 ! (out1) $}) ;
		\draw
			(in1) to [out = south , in = west] node [fromabove] {\mathrm{forget}}
			(cup) to [out = east , in = south] node [toabove] {\mathrm{diag}}
			(in2)
		;
		\draw (in3) to [out = south , in = north]
			node [fromabove] {\mathrm{forget}} coordinate [pos = 1/2] (forget) node [tobelow] {\mathrm{forget}}
			(out1)
		;
		\draw
			(out2) to [out = north , in = west] node [frombelow] {\mathrm{diag}}
			(cap) to [out = east , in = north] node [tobelow] {\mathrm{forget}}
			(out3)
		;
		\node [bead] at (cup) {\mathrm{eval}} ;
		\node [bead] at (cap) {\mathrm{cone}} ;
		\node (align) at (forget) {} ;
	\end{tikzpicture}
$$

Given a $*$-category $\cat{C}$,
$κ \cat{C}$ transforms a diagram of $\cat{C}$-cells
into the singleton diagram of its composite $\cat{C}$-cell.
This is just how we want compositors to behave.
Moreover, $κ$ is idempotent by the right adjunction law.
This matches our expectation that if we compose an already-composed diagram then we leave it unchanged.

The comparison structure of a \refer{lax functor} $F : \hom{\cat{C}}{\cat{D}}$
allows us to pass a compositor sheet backward through the functor's surface.
This gives us a morphism
$θ : \hom{\comp{(\comp[2]{F , \mathrm{fdf}}) , (\comp[2]{\cat{D} , κ})}}{\comp{(\comp[2]{\cat{C} , κ}) , (\comp[2]{F , \mathrm{fdf}})}}$,
which we can think of as an \emph{oplax} interchanger of $F$ and $κ$:
$$
	\begin{tikzpicture}[string diagram , x = {(16mm , -8mm)} , y = {(0mm , -16mm)} , z = {(24mm , 0mm)} , baseline=(align.base)]
		\node at (1/2 , 1/2 , 3/2) {2\Cat{Ctd}} ;
		\coordinate (sheet) at (1/2 , 1/2 , 1) ;
		\draw [sheet]
			($ (sheet) + (-1/2 , -1/2 , 0) $) to coordinate [pos = 3/4] (in)
			($ (sheet) + (1/2 , -1/2 , 0) $) to coordinate [pos = 3/4] (to)
			($ (sheet) + (1/2 , 1/2 , 0) $) to coordinate [pos = 3/4] (out)
			($ (sheet) + (-1/2 , 1/2 , 0) $) to coordinate [pos = 3/4] (from)
			cycle
		;
		\draw [on sheet] (in) to [out = south , in = north] node [fromabove] {κ} node [tobelow] {κ} (out) ;
		\node [anchor = west] at (from) {\mathrm{fdf}} ;
		\node [anchor = east] at (to) {\mathrm{fdf}} ;
		\node at (1/2 , 1/2 , 1/2) {2\Cat{Cat}} ;
		\coordinate (sheet) at (1/2 , 1/2 , 0) ;
		\draw [sheet]
			($ (sheet) + (-1/2 , -1/2 , 0) $) to coordinate [pos = 1/4] (in)
			($ (sheet) + (1/2 , -1/2 , 0) $) to coordinate [pos = 0] (to)
			($ (sheet) + (1/2 , 1/2 , 0) $) to coordinate [pos = 1/4] (out)
			($ (sheet) + (-1/2 , 1/2 , 0) $) to coordinate [pos = 0] (from)
			cycle
		;
		\draw [on sheet] (in) to [out = south , in = north] node [fromabove] {F} node [tobelow] {F} (out) ;
		\node [anchor = south west] at (from) {\cat{C}} ;
		\node [anchor = north east] at (to) {\cat{D}} ;
		\node (align) at (1/2 , 1/2 , -1/2) {\cat{1}} ;
	\end{tikzpicture}
	\quad ≔ \quad
	\begin{tikzpicture}[string diagram , x = {(18mm , -8mm)} , y = {(0mm , -18mm)} , z = {(6mm , 3mm)} , baseline=(current bounding box.center)]
		\coordinate (cell) at (1/2 , 1/2 , 1/2) ;
		\node at (1/2 , 1/2 , 3) {2\Cat{Ctd}} ;
		\coordinate (sheet) at (1/2 , 1/2 , 1) ;
		\draw [sheet]
			($ (sheet) + (-1/2 , -1/2 , 0) $) to coordinate [pos = 3/4] (in)
			($ (sheet) + (1/2 , -1/2 , 0) $) to coordinate [pos = 3/4] (to)
			($ (sheet) + (1/2 , 1/2 , 0) $) to coordinate [pos = 3/4] (out)
			($ (sheet) + (-1/2 , 1/2 , 0) $) to coordinate [pos = 3/4] (from)
			cycle
		;
		\draw [on sheet]
			(in) to [out = {180 + 15} , in = {90 - 15}] node [fromabove] {κ}
			(cell) to [out = {270 - 15} , in = {90 - 0}] node [tobelow] {κ}
			(out)
		;
		\node [anchor = west] at (from) {\mathrm{fdf}} ;
		\node [anchor = east] at (to) {\mathrm{fdf}} ;
		\coordinate (sheet) at (1/2 , 1/2 , 0) ;
		\draw [sheet]
			($ (sheet) + (-1/2 , -1/2 , 0) $) to coordinate [pos = 1/4] (in)
			($ (sheet) + (1/2 , -1/2 , 0) $) to coordinate [pos = 0] (to)
			($ (sheet) + (1/2 , 1/2 , 0) $) to coordinate [pos = 1/4] (out)
			($ (sheet) + (-1/2 , 1/2 , 0) $) to coordinate [pos = 0] (from)
			cycle
		;
		\draw [on sheet]
			(in) to [out = {0 - 30} , in = {90 + 15}] node [fromabove] {F}
			(cell) to [out = {270 +15} , in = {90 + 0}] node [tobelow] {F}
			(out)
		;
		\node [anchor = south west] at (from) {\cat{C}} ;
		\node [anchor = north east] at (to) {\cat{D}} ;
		\node [bead] at (cell) {θ} ;
		\node at (1/2 , 1/2 , -2) {\cat{1}} ;
	\end{tikzpicture}
$$

The composition structure of lax functors \eqref{lax functor composition}
says that these interchangers respect composition in their first index \eqref{composite interchangers}:
$$
	\begin{tikzpicture}[string diagram , x = {(12mm , -6mm)} , y = {(0mm , -12mm)} , z = {(16mm , 0mm)} , baseline=(align.base)]
		\node at (1/2 , 1/2 , 3/2) {} ;
		\coordinate (sheet) at (1/2 , 1/2 , 1) ;
		\draw [sheet]
			($ (sheet) + (-1/2 , -1/2 , 0) $) to coordinate [pos = 3/4] (in)
			($ (sheet) + (1/2 , -1/2 , 0) $) to coordinate [pos = 3/4] (to)
			($ (sheet) + (1/2 , 1/2 , 0) $) to coordinate [pos = 3/4] (out)
			($ (sheet) + (-1/2 , 1/2 , 0) $) to coordinate [pos = 3/4] (from)
			cycle
		;
		\draw [on sheet] (in) to [out = south , in = north] node [fromabove] {κ} node [tobelow] {κ} (out) ;
		\node [anchor = north west] at (from) {\mathrm{fdf}} ;
		\node [anchor = south east] at (to) {\mathrm{fdf}} ;
		\node at (1/2 , 1/2 , 1/2) {} ;
		\coordinate (sheet) at (1/2 , 1/2 , 0) ;
		\draw [sheet]
			($ (sheet) + (-1/2 , -1/2 , 0) $) to coordinate [pos = 1/6] (in 1) coordinate [pos = 3/6] (in 2)
			($ (sheet) + (1/2 , -1/2 , 0) $) to coordinate [pos = 0] (to)
			($ (sheet) + (1/2 , 1/2 , 0) $) to coordinate [pos = 1/6] (out 2) coordinate [pos = 3/6] (out 1)
			($ (sheet) + (-1/2 , 1/2 , 0) $) to coordinate [pos = 0] (from)
			cycle
		;
		\draw [on sheet] (in 1) to [out = south , out looseness = 1.25 , in looseness = 0.75 , in = north] node [fromabove] {F} node [tobelow] {F} (out 1) ;
		\draw [on sheet] (in 2) to [out = south , out looseness = 0.75 , in looseness = 1.25 , in = north] node [fromabove] {G} node [tobelow] {G} (out 2) ;
		\node [anchor = south west] at (from) {\cat{C}} ;
		\node [anchor = center] at ($ (from) ! 1/2 ! (to) $) {\cat{D}} ;
		\node [anchor = north east] at (to) {\cat{E}} ;
		\node (align) at (1/2 , 1/2 , -1/2) {} ;
	\end{tikzpicture}
$$
The comparitor compatibility requirement of an \refer{oplax transformation between lax functors}
\eqref{nullary comparitor compatibility} and \eqref{binary comparitor compatibility}
say that these interchangers are natural in their first index \eqref{interchanger naturality}:
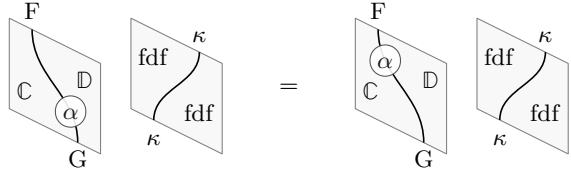
$$
	\begin{tikzpicture}[string diagram , x = {(12mm , -6mm)} , y = {(0mm , -12mm)} , z = {(16mm , 0mm)} , baseline=(align.base)]
		\node at (1/2 , 1/2 , 3/2) {} ;
		\coordinate (sheet) at (1/2 , 1/2 , 1) ;
		\draw [sheet]
			($ (sheet) + (-1/2 , -1/2 , 0) $) to coordinate [pos = 3/4] (in)
			($ (sheet) + (1/2 , -1/2 , 0) $) to coordinate [pos = 3/4] (to)
			($ (sheet) + (1/2 , 1/2 , 0) $) to coordinate [pos = 3/4] (out)
			($ (sheet) + (-1/2 , 1/2 , 0) $) to coordinate [pos = 3/4] (from)
			cycle
		;
		\draw [on sheet] (in) to [out = south , in = north] node [fromabove] {κ} node [tobelow] {κ} (out) ;
		\node [anchor = north west] at (from) {\mathrm{fdf}} ;
		\node [anchor = south east] at (to) {\mathrm{fdf}} ;
		\node at (1/2 , 1/2 , 1/2) {} ;
		\coordinate (sheet) at (1/2 , 1/2 , 0) ;
		\draw [sheet]
			($ (sheet) + (-1/2 , -1/2 , 0) $) to coordinate [pos = 1/4] (in)
			($ (sheet) + (1/2 , -1/2 , 0) $) to coordinate [pos = 0] (to)
			($ (sheet) + (1/2 , 1/2 , 0) $) to coordinate [pos = 1/4] (out)
			($ (sheet) + (-1/2 , 1/2 , 0) $) to coordinate [pos = 0] (from)
			cycle
		;
		\draw [on sheet] (in) to [out = south , in = north] node [fromabove] {F} coordinate [pos = 3/4] (cell) node [tobelow] {G} (out) ;
		\node [bead] at (cell) {α} ;
		\node [anchor = south west] at (from) {\cat{C}} ;
		\node [anchor = north east] at (to) {\cat{D}} ;
		\node (align) at (1/2 , 1/2 , -1/2) {} ;
	\end{tikzpicture}
	\quad = \quad
	\begin{tikzpicture}[string diagram , x = {(12mm , -6mm)} , y = {(0mm , -12mm)} , z = {(16mm , 0mm)} , baseline=(align.base)]
		\node at (1/2 , 1/2 , 3/2) {} ;
		\coordinate (sheet) at (1/2 , 1/2 , 1) ;
		\draw [sheet]
			($ (sheet) + (-1/2 , -1/2 , 0) $) to coordinate [pos = 3/4] (in)
			($ (sheet) + (1/2 , -1/2 , 0) $) to coordinate [pos = 3/4] (to)
			($ (sheet) + (1/2 , 1/2 , 0) $) to coordinate [pos = 3/4] (out)
			($ (sheet) + (-1/2 , 1/2 , 0) $) to coordinate [pos = 3/4] (from)
			cycle
		;
		\draw [on sheet] (in) to [out = south , in = north] node [fromabove] {κ} node [tobelow] {κ} (out) ;
		\node [anchor = north west] at (from) {\mathrm{fdf}} ;
		\node [anchor = south east] at (to) {\mathrm{fdf}} ;
		\node at (1/2 , 1/2 , 1/2) {} ;
		\coordinate (sheet) at (1/2 , 1/2 , 0) ;
		\draw [sheet]
			($ (sheet) + (-1/2 , -1/2 , 0) $) to coordinate [pos = 1/4] (in)
			($ (sheet) + (1/2 , -1/2 , 0) $) to coordinate [pos = 0] (to)
			($ (sheet) + (1/2 , 1/2 , 0) $) to coordinate [pos = 1/4] (out)
			($ (sheet) + (-1/2 , 1/2 , 0) $) to coordinate [pos = 0] (from)
			cycle
		;
		\draw [on sheet] (in) to [out = south , in = north] node [fromabove] {F} coordinate [pos = 1/4] (cell) node [tobelow] {G} (out) ;
		\node [bead] at (cell) {α} ;
		\node [anchor = south west] at (from) {\cat{C}} ;
		\node [anchor = north east] at (to) {\cat{D}} ;
		\node (align) at (1/2 , 1/2 , -1/2) {} ;
	\end{tikzpicture}
$$

The graphical calculus of compositor sheets for surface diagrams can be seen as a ``de-projection'' of
the graphical calculus of \define[functor box]{functor boxes} for string diagrams
\cite{seely-1999-linearly_distributive_functors, mellies-2006-string_diagrams}.
In the latter, parts of a string diagram may be annotated with boxes
labeled by a lax or oplax functor of $2$-categories $F : \hom{\cat{C}}{\cat{D}}$.
The part of a diagram outside of an $F$-box is a diagram with a ``hole'' in it in $F$'s codomain $2$-category $\cat{D}$.
The part inside is a diagram in $F$'s domain $2$-category $\cat{C}$
whose $F$-image is a $\cat{D}$-disk that fits in the hole to complete the diagram in $\cat{D}$.

In accordance with globularity,
a side boundary of an $F$-box may bound only a single region corresponding to a $\cat{C}$-object.
Outside of the box this region is labeled by the $F$-image of that $\cat{C}$-object,
while inside it is labeled by that $\cat{C}$-object itself.
Because $F$ is applied to everything inside the box these both represent the same $\cat{D}$-object
and this part of the box boundary acts as an identity.

Regardless of whether $F$ is a lax or oplax functor
a single wire corresponding to a $\cat{C}$-arrow may pass through the top or bottom boundary of an $F$-box.
Outside of the box this wire is labeled by the $F$-image of that $\cat{C}$-arrow,
while inside it is labeled by that $\cat{C}$-arrow itself.
Because $F$ is applied to everything inside the box these both represent the same $\cat{D}$-arrow
and this part of the box boundary also acts as an identity.
For example, drawing a $\cat{C}$-disk $φ : \hom{f}{g}$ within an $F$-box represents the $\cat{D}$-disk $F φ : \hom{F f}{F g}$.
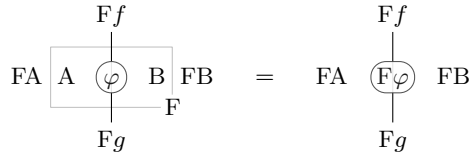
$$
	\begin{tikzpicture}[string diagram , x = {(24mm , 0mm)} , y = {(0mm , -12mm)} , baseline=(align.base)]
		\node (box) [callout , minimum width = 16mm , minimum height = 8mm] at (1/2 , 1/2) {} ;
		\node (align) [online] at (box.south east) {F} ;
		\coordinate (f) at ($ (box.west) ! 1/2 ! (box.east) $) ;
		\draw (f |- 0 , 0) to node [fromabove] {F f} coordinate [pos = 1/2] (cell) node (f |- 1 , 1) [tobelow] {F g} (f |- 1 , 1)
		;
		\node [bead] (align) at (cell) {φ} ;
		\node [anchor = east] at (box.west) {F A} ;
		\node [anchor = west] at (box.west) {A} ;
		\node [anchor = east] at (box.east) {B} ;
		\node [anchor = west] at (box.east) {F B} ;
	\end{tikzpicture}
	\quad = \quad
	\begin{tikzpicture}[string diagram , x = {(24mm , 0mm)} , y = {(0mm , -12mm)} , baseline=(align.base)]
		\coordinate (f) at (1/2 , 1/2) ;
		\draw (f |- 0 , 0) to node [fromabove] {F f}
			coordinate [pos = 1/2] (cell)
			node (f |- 1 , 1) [tobelow] {F g} (f |- 1 , 1)
		;
		\node [bead] (align) at (cell) {F φ} ;
		\node (FA) at ($ (f) + (-1/3 , 0) $) {F A} ;
		\node (FB) at ($ (f) + (1/3 , 0) $) {F B} ;
	\end{tikzpicture}
$$

If $F$ is a lax functor then the passage of any number of wires through the top boundary of an $F$-box
represents the corresponding lax comparison disk
$θ \tuple*{f_1 , ⋯ , f_n} : \hom{\comp{F f_1 , … , F f_n}}{F (\comp{f_1 , … , f_n})}$.
Note that when $n = 1$ we have $θ \tuple*{f} = \comp{} (F f)$, agreeing with the case above.
A $\cat{C}$-diagram drawn within an $F$-box represents the $F$-image of the composite of that diagram in $\cat{D}$.
For example, drawing a $\cat{C}$-disk $φ : \hom{\comp{f_1 , … , f_n}}{g}$ within an $F$-box
represents the $\cat{D}$-diagram $\comp{θ \tuple*{f_1 , ⋯ , f_n} , F φ}$.
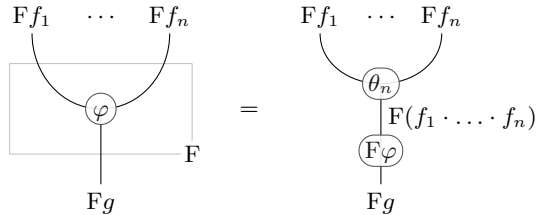
$$
	\begin{tikzpicture}[string diagram , x = {(16mm , 0mm)} , y = {(0mm , -20mm)} , baseline=(align.base)]
		\coordinate (merge) at (1/2 , 1/2) ;
		\node (box) [callout , minimum width = 24mm , minimum height = 12mm] at (merge) {} ;
		\node (align) [online] at (box.south east) {F} ;
		\coordinate (f1) at ($ (box.west) ! 1/8 ! (box.east) $) ;
		\coordinate (fn) at ($ (box.west) ! 7/8 ! (box.east) $) ;
		\coordinate (g) at ($ (box.west) ! 1/2 ! (box.east) $) ;
		\draw (f1 |- 0 , 0) to [out = south , in = west] node  (F f1) [fromabove] {F f_1} (merge) ;
		\draw (fn |- 0 , 0) to [out = south , in = east] node (F fn) [fromabove] {F f_n} (merge) ;
		\draw (merge) to node [tobelow] {F g} (box |- 1 , 1) ;
		\node [bead] (align) at (box) {φ} ;
		\node at ($ (F f1) ! 1/2 ! (F fn) $) {⋯} ;
	\end{tikzpicture}
	\quad = \quad
	\begin{tikzpicture}[string diagram , x = {(16mm , 0mm)} , y = {(0mm , -20mm)} , baseline=(align.base)]
		\coordinate (merge) at (1/2 , 1/3) ;
		\draw ($ (merge |- 0 , 0) + (-1/2 , 0) $) to [out = south , in = west] node (F f1) [fromabove] {F f_1} coordinate [pos = 1/6] (f1) (merge) ;
		\draw ($ (merge |- 0 , 0) + (1/2 , 0) $) to [out = south , in = east] node (F fn) [fromabove] {F f_n} coordinate [pos = 1/6] (fn) (merge) ;
		\draw (merge) to node [pos = 1/3] {F (\comp{f_1 , … , f_n})} coordinate [pos = 2/3] (d) node [tobelow] {F g} (merge |- 1 , 1) ;
		\node [bead] at (merge) {θ_n} ;
		\node [bead] at (d) {F φ} ;
		\node at ($ (F f1) ! 1/2 ! (F fn) $) {⋯} ;
		\node (align) at (1/2 , 1/2) {} ;
	\end{tikzpicture}
$$

Unless $F$ is oplax as well as lax, we do not generally have a morphism $\hom{F (\comp{f_1 , … , f_n})}{\comp{F f_1 , … , F f_n}}$.
In this case the bottom boundary of an $F$-box represents the identity on $F (\comp{f_1 , … , f_n})$.
Graphically, this has $n$ $\cat{C}$-wires coming in from above
and a single $\cat{D}$-wire labeled by $F$ of their composite leaving from below.
For example, the nullary and binary lax comparitor components \eqref{lax comparitor components} are represented as follows.
$$
	
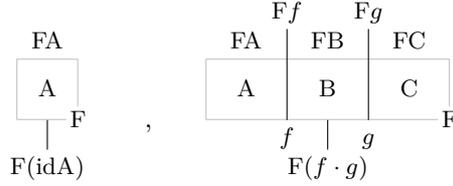
\begin{tikzpicture}[string diagram , x = {(36mm , 0mm)} , y = {(0mm , -16mm)} , baseline=(align.base)]
		\node (start) [callout , minimum width = 8mm , minimum height = 8mm] at (1/2 , 1/2) {} ;
		\node (align) [online] at (start.south east) {F} ;
		\draw (start.south) to node [tobelow] {F (\comp{} A)} (start.south |- 1 , 1) ;
		\node at (start) {A} ;
		\node [anchor = south] at (start.north) {F A} ;
	\end{tikzpicture}
	\qquad \text{,} \qquad
	\begin{tikzpicture}[string diagram , x = {(36mm , 0mm)} , y = {(0mm , -16mm)} , baseline=(align.base)]
		\node (merge) [callout , minimum width = 32mm , minimum height = 8mm] at (1/2 , 1/2) {} ;
		\node (align) [online] at (merge.south east) {F} ;
		\coordinate (f) at ($ (merge.west) ! 1/3 ! (merge.east) $) ;
		\coordinate (g) at ($ (merge.west) ! 2/3 ! (merge.east) $) ;
		\draw (f |- 0 , 0) to node [fromabove] {F f} node [tobelow] {f} (f |- merge.south) ;
		\draw (g |- 0 , 0) to node [fromabove] {F g} node [tobelow] {g} (g |- merge.south) ;
		\draw (merge.south) to node [tobelow] {F (\comp{f , g})} (merge.south |- 1 , 1) ;
		\node (A) at ($ (merge.west) ! 1/2 ! (f) $) {A} ;
		\node (B) at ($ (f) ! 1/2 ! (g) $) {B} ;
		\node (C) at ($ (g) ! 1/2 ! (merge.east) $) {C} ;
		\node [anchor = south] at (A |- merge.north) {F A} ;
		\node [anchor = south] at (B |- merge.north) {F B} ;
		\node [anchor = south] at (C |- merge.north) {F C} ;
	\end{tikzpicture}
$$

The lax comparitor unit laws \eqref{comparitor unit law} are drawn as follows.
$$
	
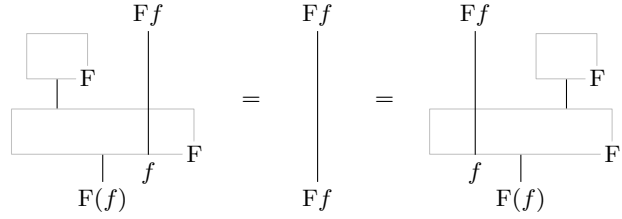
\begin{tikzpicture}[string diagram , x = {(36mm , 0mm)} , y = {(0mm , -20mm)} , baseline=(current bounding box.center)]
		\coordinate (f) at (1/3 , 1/2) ;
		\coordinate (g) at (2/3 , 1/2) ;
		\node (start) [callout , minimum width = 8mm , minimum height = 6mm] at (f |- 1/4 , 1/6) {} ;
		\node [online] at (start.south east) {F} ;
		\node (merge) [callout , minimum width = 24mm , minimum height = 6mm] at (1/2 , 2/3) {} ;
		\node [online] at (merge.south east) {F} ;
		\draw (start.south) to (start.south |- merge.north) ;
		\draw (g |- 0 , 0) to node [fromabove] {F f} node [tobelow] {f} (g |- merge.south) ;
		\draw (merge.south) to node [tobelow] {F (f)} (merge.south |- 1 , 1) ;
	\end{tikzpicture}
	\quad = \quad
	\begin{tikzpicture}[string diagram , x = {(12mm , 0mm)} , y = {(0mm , -20mm)} , baseline=(current bounding box.center)]
		\coordinate (merge) at (1/2 , 1/2) ;
		\draw (merge |- 0 , 0) to [out = south , in = north] node [fromabove] {F f} node [tobelow] {F f} (merge |- 1 , 1) ;
	\end{tikzpicture}
	\quad = \quad
	\begin{tikzpicture}[string diagram , x = {(36mm , 0mm)} , y = {(0mm , -20mm)} , baseline=(current bounding box.center)]
		\coordinate (f) at (1/3 , 1/2) ;
		\coordinate (g) at (2/3 , 1/2) ;
		\node (start) [callout , minimum width = 8mm , minimum height = 6mm] at (g |- 1/4 , 1/6) {} ;
		\node [online] at (start.south east) {F} ;
		\node (merge) [callout , minimum width = 24mm , minimum height = 6mm] at (1/2 , 2/3) {} ;
		\node [online] at (merge.south east) {F} ;
		\draw (start.south) to (start.south |- merge.north) ;
		\draw (f |- 0 , 0) to node [fromabove] {F f} node [tobelow] {f} (f |- merge.south) ;
		\draw (merge.south) to node [tobelow] {F (f)} (merge.south |- 1 , 1) ;
	\end{tikzpicture}
$$
The lax comparitor associative law \eqref{comparitor associative law} is drawn like this.
$$
	
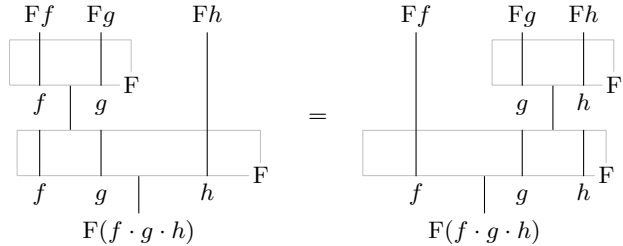
\begin{tikzpicture}[string diagram , x = {(36mm , 0mm)} , y = {(0mm , -24mm)} , baseline=(current bounding box.center)]
		\coordinate (f) at (1/4 , 1/2) ;
		\coordinate (g) at (3/4 , 1/2) ;
		\node (merge') [callout , minimum width = 16mm , minimum height = 6mm] at (f |- 1/4 , 1/6) {} ;
		\node [online] at (merge'.south east) {F} ;
		\node (merge) [callout , minimum width = 32mm , minimum height = 6mm] at (1/2 , 2/3) {} ;
		\node [online] at (merge.south east) {F} ;
		\coordinate (f0) at ({$ (merge'.west) ! 1/4 ! (merge'.east) $} |- f) ;
		\coordinate (f1) at ({$ (merge'.west) ! 3/4 ! (merge'.east) $} |- f) ;
		\draw (f0 |- 0 , 0) to node [fromabove] {F f} node [tobelow] {f} (f0 |- merge'.south) ;
		\draw (f1 |- 0 , 0) to node [fromabove] {F g} node [tobelow] {g} (f1 |- merge'.south) ;
		\draw (merge'.south) to (merge'.south |- merge.north) ;
		\draw (f0 |- merge.north) to node [tobelow] {f} (f0 |- merge.south) ;
		\draw (f1 |- merge.north) to node [tobelow] {g} (f1 |- merge.south) ;
		\draw (g |- 0 , 0) to node [fromabove] {F h} node [tobelow] {h} (g |- merge.south) ;
		\draw (merge.south) to node [tobelow] {F (\comp{f , g , h})} (merge.south |- 1 , 1) ;
	\end{tikzpicture}
	\quad = \quad
	\begin{tikzpicture}[string diagram , x = {(36mm , 0mm)} , y = {(0mm , -24mm)} , baseline=(current bounding box.center)]
		\coordinate (f) at (1/4 , 1/2) ;
		\coordinate (g) at (3/4 , 1/2) ;
		\node (merge') [callout , minimum width = 16mm , minimum height = 6mm] at (g |- 1/4 , 1/6) {} ;
		\node [online] at (merge'.south east) {F} ;
		\node (merge) [callout , minimum width = 32mm , minimum height = 6mm] at (1/2 , 2/3) {} ;
		\node [online] at (merge.south east) {F} ;
		\coordinate (g0) at ({$ (merge'.west) ! 1/4 ! (merge'.east) $} |- g) ;
		\coordinate (g1) at ({$ (merge'.west) ! 3/4 ! (merge'.east) $} |- g) ;
		\draw (g0 |- 0 , 0) to node [fromabove] {F g} node [tobelow] {g} (g0 |- merge'.south) ;
		\draw (g1 |- 0 , 0) to node [fromabove] {F h} node [tobelow] {h} (g1 |- merge'.south) ;
		\draw (merge'.south) to (merge'.south |- merge.north) ;
		\draw (f |- 0 , 0) to node [fromabove] {F f} node [tobelow] {f} (f |- merge.south) ;
		\draw (g0 |- merge.north) to node [tobelow] {g} (g0 |- merge.south) ;
		\draw (g1 |- merge.north) to node [tobelow] {h} (g1 |- merge.south) ;
		\draw (merge.south) to node [tobelow] {F (\comp{f , g , h})} (merge.south |- 1 , 1) ;
	\end{tikzpicture}
$$

Together with comparitor composition naturality \eqref{comparitor composition naturality},
these imply that boxes for lax functors support tree composition.
For example, given $\cat{C}$-disks $φ : \hom{\comp{f , g}}{\comp{f′ , g′}}$, $ψ : \hom{\comp{h , i}}{\comp{h′ , i′}}$,
and $σ : \hom{\comp{g′ , h′}}{\comp{g′′ , h′′}}$ we have the equation of string diagrams,
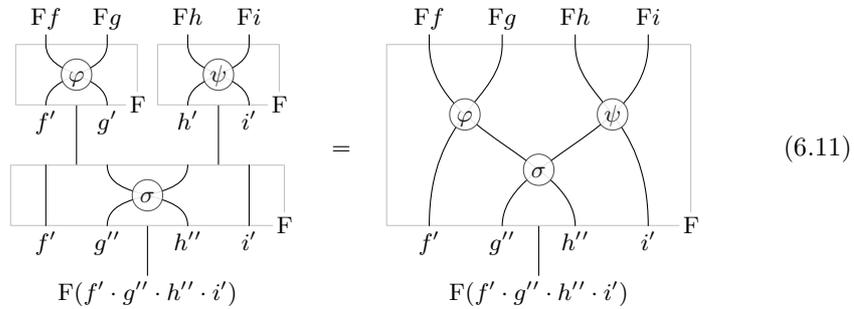
\begin{equation} \label{functor box example}
	\begin{tikzpicture}[string diagram , x = {(56mm , 0mm)} , y = {(0mm , -32mm)} , baseline=(current bounding box.center)]
		\coordinate (f) at (1/3 , 1/2) ;
		\coordinate (g) at (2/3 , 1/2) ;
		\node (merge0) [callout , minimum width = 16mm , minimum height = 8mm] at (f |- 1/4 , 1/6) {} ;
		\node [online] at (merge0.south east) {F} ;
		\node (merge1) [callout , minimum width = 16mm , minimum height = 8mm] at (g |- 1/4 , 1/6) {} ;
		\node [online] at (merge1.south east) {F} ;
		\node (merge) [callout , minimum width = 36mm , minimum height = 8mm] at (1/2 , 2/3) {} ;
		\node [online] at (merge.south east) {F} ;
		\coordinate (f0) at ({$ (merge0.west) ! 1/4 ! (merge0.east) $} |- f) ;
		\coordinate (f1) at ({$ (merge0.west) ! 3/4 ! (merge0.east) $} |- f) ;
		\draw [name path = g]
			(f0 |- 0 , 0) to node [fromabove] {F f}
			(f0 |- merge0.north) to [out = south , in = north] node [tobelow] {g′}
			(f1 |- merge0.south)
		;
		\draw [name path = f]
			(f1 |- 0 , 0) to node [fromabove] {F g}
			(f1 |- merge0.north) to [out = south , in = north] node [tobelow] {f′}
			(f0 |- merge0.south)
		;
		\path [name intersections={of = f and g , by = {cell0}}] ;
		\draw (merge0.south) to (merge0.south |- merge.north) ;
		\coordinate (g0) at ({$ (merge1.west) ! 1/4 ! (merge1.east) $} |- g) ;
		\coordinate (g1) at ({$ (merge1.west) ! 3/4 ! (merge1.east) $} |- g) ;
		\draw [name path = i]
			(g0 |- 0 , 0) to node [fromabove] {F h}
			(g0 |- merge1.north) to [out = south , in = north] node [tobelow] {i′}
			(g1 |- merge1.south)
		;
		\draw [name path = h]
			(g1 |- 0 , 0) to node [fromabove] {F i}
			(g1 |- merge1.north) to [out = south , in = north] node [tobelow] {h′}
			(g0 |- merge1.south)
		;
		\path [name intersections={of = h and i , by = {cell1}}] ;
		\draw (merge1.south) to (merge1.south |- merge.north) ;
		\draw (f0 |- merge.north) to [out = south , in = north] node [tobelow] {f′} (f0 |- merge.south) ;
		\draw [name path = h']
			(f1 |- merge.north) to [out = south , in = north] node [tobelow] {h′′}
			(g0 |- merge.south)
		;
		\draw [name path = g']
			(g0 |- merge.north) to [out = south , in = north] node [tobelow] {g′′}
			(f1 |- merge.south)
		;
		\path [name intersections={of = g' and h' , by = {cell}}] ;
		\draw (g1 |- merge.north) to [out = south , in = north] node [tobelow] {i′} (g1 |- merge.south) ;
		\draw (merge.south) to node [tobelow] {F (\comp{f′ , g′′ , h′′ , i′})} (merge.south |- 1 , 1) ;
		\node [bead] at (cell0) {φ} ;
		\node [bead] at (cell1) {ψ} ;
		\node [bead] at (cell) {σ} ;
	\end{tikzpicture}
	\quad = \quad
	\begin{tikzpicture}[string diagram , x = {(48mm , 0mm)} , y = {(0mm , -32mm)} , baseline=(current bounding box.center)]
		\coordinate (f) at (1/5 , 1/2) ;
		\coordinate (g) at (2/5 , 1/2) ;
		\coordinate (h) at (3/5 , 1/2) ;
		\coordinate (i) at (4/5 , 1/2) ;
		\node (box) [callout , minimum width = 40mm , minimum height = 24mm] at (1/2 , 5/12) {} ;
		\node [online] at (box.south east) {F} ;
		\draw [name path = f]
			(f |- 0 , 0) to node [fromabove] {F f}
			(f |- box.north) to [out = south , in looseness = 0.75 , in = north] node [tobelow] {h′′}
			(h |- box.south)
		;
		\draw [name path = g]
			(g |- 0 , 0) to node [fromabove] {F g}
			(g |- box.north) to [out = south , out looseness = 0.75 , in looseness = 1.5 , in = north] node [tobelow] {f′}
			(f |- box.south)
		;
		\draw [name path = h]
			(h |- 0 , 0) to node [fromabove] {F h}
			(h |- box.north) to [out = south , out looseness = 0.75 , in looseness = 1.5 , in = north] node [tobelow] {i′}
			(i |- box.south)
		;
		\draw [name path = i]
			(i |- 0 , 0) to node [fromabove] {F i}
			(i |- box.north) to [out = south , in looseness = 0.75 , in = north] node [tobelow] {g′′}
			(g |- box.south)
		;
		\path [name intersections={of = f and g , by = {cell0}}] ;
		\path [name intersections={of = h and i , by = {cell1}}] ;
		\path [name intersections={of = f and i , by = {cell}}] ;
		\draw (box.south) to node [tobelow] {F (\comp{f′ , g′′ , h′′ , i′})} (box.south |- 1 , 1) ;
		\node [bead] at (cell0) {φ} ;
		\node [bead] at (cell1) {ψ} ;
		\node [bead] at (cell) {σ} ;
	\end{tikzpicture}
\end{equation}
corresponding to the calculation
$$
	\begin{relationalreasoning}
		\term[]
		{\comp{[\comp[2]{(\comp{θ \tuple*{f , g} , F φ}) , (\comp{θ \tuple*{h , i} , F ψ})}] , [\comp{θ \tuple*{\comp{f′, g′} , \comp{h′, i′}} , F (\comp[2]{f′ , σ , i′})}]}}
		\term[=]
		{\comp{(\comp[2]{θ \tuple*{f , g} , θ \tuple*{h , i}}) , (\comp[2]{F φ , F ψ}) , θ \tuple*{\comp{f′, g′} , \comp{h′, i′}} , F (\comp[2]{f′ , σ , i′})}}
		\term[=]
		{\comp{(\comp[2]{θ \tuple*{f , g} , θ \tuple*{h , i}}) , θ \tuple*{\comp{f, g} , \comp{h, i}} , F (\comp[2]{φ , ψ}) , F (\comp[2]{f′ , σ , i′})}}
		\term[=]
		{\comp{θ \tuple*{f , g , h , i} , F (\comp{(\comp[2]{φ , ψ}) , (\comp[2]{f′ , σ , i′})})} .}
	\end{relationalreasoning}
$$

The string diagrammatics for oplax functors is simply the vertical reflection of that for lax functors.
In general
we can tell whether the top or bottom boundary of a functor box represents a comparison disk
by the number of wires connected to it on the outside.

The string diagram calculus for functor boxes can be seen as a projection of the surface diagram calculus for compositor sheets
where string diagrams are projected onto the compositors as a movie is projected onto a cinema screen.
The top boundary of a box for a lax functor $F$ represents its lax comparitor component.
This corresponds to the emergence of the compositor sheet from behind the surface representing $F$.
Unless $F$ is oplax as well, the compositor sheet remains in front of $F$'s surface for the duration of the surface diagram,
while in the string diagram the corresponding $F$-box seems to end at its bottom boundary.
But this is just a cinematic illusion: the bottom boundary of an $F$-box
represents the identity on the $F$-image of the composite of the arrows connected to it from the inside.
Indeed, some authors prefer to depict lax functors in string diagrams using ``tubes'' rather than ``boxes''.
These tubes simply extend the bottom boundaries of functor boxes down to top boundary of the next functor box.

The equation between string diagrams depicted in \eqref{functor box example} arises by projection from the following surface diagram,
where in the string diagram on the left the disks $φ$ and $ψ$ are positioned vertically between the comparitor components
and in the one on the right they are positioned below the lower comparitor component.
\begin{equation} \label{compositor sheet example}
	\begin{tikzpicture}[string diagram , x = {(20mm , -10mm)} , y = {(0mm , -20mm)} , z = {(32mm , 0mm)} , baseline=(align.base)]
		\node at (1/2 , 1/2 , 3/2) {} ;
		\coordinate (sheet) at (1/2 , 1/2 , 1) ;
		\draw [compositor]
			($ (sheet) + (-1/2 , -5/8 , 1/3) $) to
			($ (sheet) + (1/2 , -5/8 , 1/3) $) to [out = south , in = north east]
			($ (sheet) + (1/2 , 1/8 , 0) $) to
			($ (sheet) + (-1/2 , 1/8 , 0) $) to [out = north east , in = south]
			cycle
		;
		\draw [compositor]
			($ (sheet) + (-1/2 , -5/8 , 1/6) $) to
			($ (sheet) + (-1/8 , -5/8 , 1/6) $) to [out = south , in = north east]
			($ (sheet) + (-1/8 , -1/4 , 0) $) to
			($ (sheet) + (-1/2 , -1/4 , 0) $) to [out = north east , in = south]
			cycle
		;
		\draw [compositor]
			($ (sheet) + (1/8 , -5/8 , 1/6) $) to
			($ (sheet) + (1/2 , -5/8 , 1/6) $) to [out = south , in = north east]
			($ (sheet) + (1/2 , -1/4 , 0) $) to
			($ (sheet) + (1/8 , -1/4 , 0) $) to [out = north east , in = south]
			cycle
		;
		\draw [sheet]
			($ (sheet) + (-1/2 , -1/2 , 0) $) to
			($ (sheet) + (1/2 , -1/2 , 0) $) to
			($ (sheet) + (1/2 , 1/2 , 0) $) to
			($ (sheet) + (-1/2 , 1/2 , 0) $) to
			cycle
		;
		\node at ($ (sheet) + (0 , -1/6 , 0) $) {F} ;
		\draw [compositor]
			($ (sheet) + (-1/2 , 1/8 , 0) $) to
			($ (sheet) + (1/2 , 1/8 , 0) $) to [out = south west , in = north]
			($ (sheet) + (1/2 , 5/8 , -1/6) $) to
			($ (sheet) + (-1/2 , 5/8 , -1/6) $) to [out = north , in = south west]
			cycle
		;
		\draw [compositor]
			($ (sheet) + (-1/2 , -1/4 , 0) $) to
			($ (sheet) + (-1/8 , -1/4 , 0) $) to [out = south west , in = north]
			($ (sheet) + (-1/8 , 5/8 , -1/3) $) to
			($ (sheet) + (-1/2 , 5/8 , -1/3) $) to [out = north , in = south west]
			cycle
		;
		\draw [compositor]
			($ (sheet) + (1/8 , -1/4 , 0) $) to
			($ (sheet) + (1/2 , -1/4 , 0) $) to [out = south west , in = north]
			($ (sheet) + (1/2 , 5/8 , -1/3) $) to
			($ (sheet) + (1/8 , 5/8 , -1/3) $) to [out = north , in = south west]
			cycle
		;
		\node at (1/2 , 1/2 , 1/2) {} ;
		\coordinate (sheet) at (1/2 , 1/2 , 0) ;
		\coordinate (cell 1) at ($ (sheet) + (-1/4 , -1/8 , 0) $) ;
		\coordinate (cell 2) at ($ (sheet) + (1/4 , -1/8 , 0) $) ;
		\coordinate (cell) at ($ (sheet) + (0 , 1/4 , 0) $) ;
		\draw [sheet]
			($ (sheet) + (-1/2 , -1/2 , 0) $) to coordinate [pos = 1/8] (in 1) coordinate [pos = 3/8] (in 2) coordinate [pos = 5/8] (in 3) coordinate [pos = 7/8] (in 4)
			($ (sheet) + (1/2 , -1/2 , 0) $) to
			($ (sheet) + (1/2 , 1/2 , 0) $)  to coordinate [pos = 1/8] (out 4) coordinate [pos = 3/8] (out 3) coordinate [pos = 5/8] (out 2) coordinate [pos = 7/8] (out 1)
			($ (sheet) + (-1/2 , 1/2 , 0) $) to
			cycle
		;
		\draw [on sheet]
			(in 1) to [out = south , in = north west] node [fromabove] {f}
			(cell 1) to [out = south east , in = north west]
			(cell) to [out = south east , in = north] node [tobelow] {h′′}
			(out 3)
		;
		\draw [on sheet]
			(in 2) to [out = south , in = north east] node [fromabove] {g}
			(cell 1) to [out = south west , in = north] node [tobelow] {f′}
			(out 1)
		;
		\draw [on sheet]
			(in 3) to [out = south , in = north west] node [fromabove] {h}
			(cell 2) to [out = south east , in = north] node [tobelow] {i′}
			(out 4)
		;
		\draw [on sheet]
			(in 4) to [out = south , in = north east] node [fromabove] {i}
			(cell 2) to [out = south west , in = north east]
			(cell) to [out = south west , in = north] node [tobelow] {g′′}
			(out 2)
		;
		\node [bead] at (cell 1) {φ} ;
		\node [bead] at (cell 2) {ψ} ;
		\node [bead] at (cell) {σ} ;
		\node (align) at (1/2 , 1/2 , -1/2) {} ;
	\end{tikzpicture}
\end{equation}

Note that in \eqref{compositor sheet example} we can't move $φ$ and $ψ$ up above the upper comparitor components
without changing the top boundary of the diagram to $\comp{F (\comp{f , g}) , F (\comp{h , i})}$.
Doing so would correspond to a string diagram in which the upper $F$-boxes extend beyond the top of the diagram as follows.
$$
	
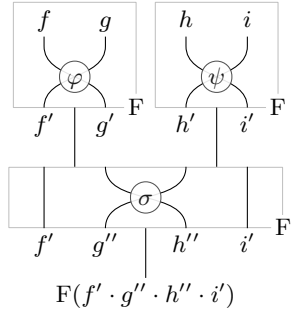
\begin{tikzpicture}[string diagram , x = {(56mm , 0mm)} , y = {(0mm , -32mm)} , baseline=(current bounding box.center)]
		\coordinate (f) at (1/3 , 1/2) ;
		\coordinate (g) at (2/3 , 1/2) ;
		\node (merge0) [minimum width = 16mm , minimum height = 8mm] at (f |- 1/4 , 1/6) {} ;
		\node (merge1) [minimum width = 16mm , minimum height = 8mm] at (g |- 1/4 , 1/6) {} ;
		\node (merge) [callout , minimum width = 36mm , minimum height = 8mm] at (1/2 , 2/3) {} ;
		\node [online] at (merge.south east) {F} ;
		\coordinate (f0) at ({$ (merge0.west) ! 1/4 ! (merge0.east) $} |- f) ;
		\coordinate (f1) at ({$ (merge0.west) ! 3/4 ! (merge0.east) $} |- f) ;
		\draw [name path = g]
			(f0 |- 0 , 0) to node (f in) [fromabove] {f}
			(f0 |- merge0.north) to [out = south , in = north] node [tobelow] {g′}
			(f1 |- merge0.south)
		;
		\draw [name path = f]
			(f1 |- 0 , 0) to node (g in) [fromabove] {g}
			(f1 |- merge0.north) to [out = south , in = north] node [tobelow] {f′}
			(f0 |- merge0.south)
		;
		\path [name intersections={of = f and g , by = {cell0}}] ;
		\node [callout , fit=(merge0)(f in)(g in)] {} ;
		\node [online] at (merge0.south east) {F} ;
		\draw (merge0.south) to (merge0.south |- merge.north) ;
		\coordinate (g0) at ({$ (merge1.west) ! 1/4 ! (merge1.east) $} |- g) ;
		\coordinate (g1) at ({$ (merge1.west) ! 3/4 ! (merge1.east) $} |- g) ;
		\draw [name path = i]
			(g0 |- 0 , 0) to node (h in) [fromabove] {h}
			(g0 |- merge1.north) to [out = south , in = north] node [tobelow] {i′}
			(g1 |- merge1.south)
		;
		\draw [name path = h]
			(g1 |- 0 , 0) to node (i in) [fromabove] {i}
			(g1 |- merge1.north) to [out = south , in = north] node [tobelow] {h′}
			(g0 |- merge1.south)
		;
		\path [name intersections={of = h and i , by = {cell1}}] ;
		\node [callout , fit=(merge1)(h in)(i in)] {} ;
		\node [online] at (merge1.south east) {F} ;
		\draw (merge1.south) to (merge1.south |- merge.north) ;
		\draw (f0 |- merge.north) to [out = south , in = north] node [tobelow] {f′} (f0 |- merge.south) ;
		\draw [name path = h']
			(f1 |- merge.north) to [out = south , in = north] node [tobelow] {h′′}
			(g0 |- merge.south)
		;
		\draw [name path = g']
			(g0 |- merge.north) to [out = south , in = north] node [tobelow] {g′′}
			(f1 |- merge.south)
		;
		\path [name intersections={of = g' and h' , by = {cell}}] ;
		\draw (g1 |- merge.north) to [out = south , in = north] node [tobelow] {i′} (g1 |- merge.south) ;
		\draw (merge.south) to node [tobelow] {F (\comp{f′ , g′′ , h′′ , i′})} (merge.south |- 1 , 1) ;
		\node [bead] at (cell0) {φ} ;
		\node [bead] at (cell1) {ψ} ;
		\node [bead] at (cell) {σ} ;
	\end{tikzpicture}
$$

\end{sect}

\begin{sect}{Conclusion} \label{section: conclusion}
	The preceding has essentially been an exercise in connecting the work of Gray \cite{gray-1974-formal_category_theory}
with that of Hummon \cite{hummon-2012-surface_diagrams}
in order to present a surface diagram calculus for $2$-category theory.
By using the canonical adjunction between Street's $2$-computads \cite{street-1996-categorical_structures} and $2$-categories
we are able to incorporate lax functors between $2$-categories as well as strict ones.
The result is a framework to represent and reason about $2$-categorical constructions
in a manner that is intended to be natural, in both the homotopical and psychological sense. 

There is of course a great deal of literature on $2$-, $3$-, $n$-, and $∞$-dimensional categorical structures.
Reaching arbitrarily high dimensions generally requires some sort of inductive or coinductive approach,
and each stage of categorification requires making choices about what structures to keep strict or make weak or directed.

$2$-categories inhabit a particularly interesting point on this spectrum
because they arise from the structure of category theory itself.
Similarly, Gray-categories are what $2$-categories and $2$-functors naturally assemble into.
But these structures are also interesting for the reason that they have
coherence theorems providing equivalences with their fully-weak counterparts,
bicategories \cite{mac_lane-1985-coherence_for_bicategories} and
tricategories \cite{street-1995-coherence_for_tricategories, gurski-2013-tricategories},
as well as compositional graphical calculi
that provide non-experts a means to help tame the daunting bureaucracy
involved in trying to do algebra in $2$ or $3$ dimensions.
Tools such as the proof assistant \emph{Globular} \cite{vicary-2018-globular} can be a great aid in this regard,
helping users discover the cells they seek by constructing diagrams incrementally
while verifying automatically that the constituent parts fit together properly along their boundaries.

\emph{Acknowledgement:} 
I'd like to express my gratitude to the members of the Compositional Systems and Methods Group at TalTech
for many enjoyable and enlightening discussions about diagrams and their semantics,
and especially to Mario Román, who encouraged me to clarify the relationship between compositors and boxes for lax functors.

\end{sect}

\begin{sect*}{References}
	\printbibliography[heading=none]
\end{sect*}

\end{document}